\documentstyle[11pt,graphicx]{article}
\def\O{{\,\rm. O\,}}

\def\E{{\cal E}}
\def\U{{\cal U}}

\def\E{{\cal E}}

\def\T{{\cal T}}

\def\U{{\cal U}}

\def\V{{\cal V}}

\def\O{{\cal O}}
\def\R{{\cal R}}
\def\A{{\cal A}}

\def\K{{\cal K}}
\def\V{{\cal V}}
\def\N{{\cal N}}
\def\ns{\noalign{\smallskip} }

\def\ds{\displaystyle}
\newfont{\Blackboard}{msbm10 scaled 1200}

\newfont{\roma}{cmr10 scaled 1200}

\def\Im{\em Im}

\newtheorem{thm}{{}\hskip\parindent Theorem}[section]
\newtheorem{lem}{{}\hskip\parindent Lemma}[section]
\newtheorem{pro}{{}\hskip\parindent Proposition}[section]
\newtheorem{exl}{{}\hskip\parindent Example}[section]

\newtheorem{rem}{{}\hskip\parindent Remark}[section]

\def\be{\begin{equation}}
\def\ee{\end{equation}}
\def\beq{\arraycolsep=1.5pt\begin{eqnarray}}
\def\eeq{\end{eqnarray}}

\def\R{I\!\!R}
\def\N{I\!\!N}
\large
\title{Second Order Optimality  Conditions for Optimal Control Problems on Riemannian Manifolds\thanks{This work is  supported by the National Science
Foundation of China under grants 11401491 and 11231007, the Fundamental research funds for the Central Universities under grant 2682014CX052, and the grant MTM2014-52347 from the Spanish
Science and Innovation Ministry.
}
\date{}
\author{Qing Cui\thanks{School of Mathematics,  Southwest Jiaotong University, Chengdu 611756, Sichuan Province, China. {\small\it E-mail:} {\small\tt
cuiqing@home.swjtu.edu.cn}.}, \ \ \ Li Deng\thanks{School of Mathematics,   Southwest Jiaotong University, Chengdu 611756, Sichuan Province, China. {\small\it E-mail:} {\small\tt
dengli@home.swjtu.edu.cn}.} \ \ \ and \ \ \ Xu Zhang\thanks{School of Mathematics, Sichuan University, Chengdu 610064, Sichuan Province, China. {\small\it E-mail:} {\small\tt
zhang$\_$xu@scu.edu.cn}.}}
}

\begin{document}
\maketitle

\begin{quote}
\begin{small}
{\bf Abstract} \,\,\,This work is concerned with an optimal control problem on a Riemannian manifold, for which two typical cases are considered.  The first case is when the endpoint is free. For this case,  the control set is assumed to be a separable metric  space. By introducing suitable dual equations, which depend on the curvature tensor of the manifold, we establish the second order necessary and sufficient optimality conditions of integral form. In particular, when the control set is a Polish space, the second order necessary condition is reduced to a pointwise form. As a key preliminary result and also an interesting byproduct, we derive a geometric
lemma, which may have some independent interest. The second case is when the endpoint is fixed. For
this more difficult case, the control set is assumed  to be open in an Euclidian space. We obtain the second order necessary and sufficient optimality conditions, in which the curvature tensor also appears explicitly. Our optimality conditions can be used to recover the following famous geometry result: Any geodesic connecting two fixed points on a Riemannian manifold satisfies the second variation of energy; while the existing optimality conditions in control literatures fail to give the same result.
\\[3mm]
\\[3mm]
{\bf MSC (2010) \,\,\,
49K15, 49K30, 93C15, 58E25, 70Q05} \\[3mm]
\end{small}
\end{quote}

\setcounter{equation}{0}

\section{Introduction }
\def\theequation{1.\arabic{equation}}
\hskip\parindent
 Let $n\in\N$ and $M$ be a complete simply connected, $n$-dimensional manifold with Riemannian metric $g$. Let $\nabla$ be the Levi-Civita connection on $M$ related to $g$, $\rho(\cdot,\cdot)$ be the distance function on $M$, $T_xM$ be the tangent space of  $M$ at $x\in M$, and $T^*_x M$ be the  cotangent space. Denote by $\langle\cdot,\cdot\rangle$ and $|\cdot|$ the inner product and the norm over $T_xM$ related to $g$, respectively. Also, denote by  $ T M \equiv \bigcup\limits_{x\in M} T_xM$, $T^* M \equiv \bigcup\limits_{x\in M} T^*_x M$ and $C^\infty(M)$ the tangent bundle, the cotangent bundle and the set of smooth functions on $M$, respectively.

Let $T>0$, $U$ be a metric  space, and  $f:[0,T]\times M\times U\to T M$ and $f^0:[0,T]\times M\times U\to \R$ be two functions (satisfying suitable assumptions to be given later). Given $y_0\in M$, let us consider the following control system
\beq\label{25}
\cases{\dot{y}(t)=f(t,y(t),u(t)),\quad a.e. \,t\in[0,T],\cr
y(0)=y_0,}
\eeq where $\dot{y}(t)=\frac{d }{dt}y(t)
$ for $t\in [0,T]$,  and $y(\cdot)$ and $u(\cdot)$ are the state and control variables valued in $M$ and $U$, respectively.
The cost functional associated with (\ref{25}) is
\beq\label{26}
J(u(\cdot))=\int_0^Tf^0(t,y(t),u(t))dt.
\eeq
In (\ref{26}), $u(\cdot)$ belongs to the following admissible control set
\begin{equation}\label{436}
\U_{ad}\equiv \big\{u(\cdot):[0,T]\to U;\ \  u(\cdot) \mbox{ is measurable}\big\}.
\end{equation}
Clearly, this is the situation without endpoint constraints.

We shall also consider the control system (\ref{25}) with the following endpoint constraint
\begin{equation}\label{503}
y(T)=y_1,
\end{equation}
for some given $y_1\in M$.  The corresponding admissible control set is then given by
$$\V_{ad}\equiv \{u(\cdot)\in\U_{ad}; \ \  y_u(T)=y_1\}, $$ where $y_u(\cdot)$ is the solution to (\ref{25}) associated to the control $u(\cdot)(\in\U_{ad})$.

In this work, we shall consider the following two optimal control problems:
\begin{description}
\item[Problem I]

{\it To find a $\bar{u}(\cdot)\in\U_{ad}$ such that}
 \begin{equation}\label{69}
J(\bar{u}(\cdot))=\inf\limits_{u(\cdot)\in\U_{ad}}J(u(\cdot));
\end{equation}
\item[Problem II] {\it To find a $\bar u(\cdot)\in\V_{ad}$ such that}
\begin{equation}\label{400}
J(\bar u(\cdot))=\inf\limits_{u(\cdot)\in\V_{ad}}J(u(\cdot)).
\end{equation}
\end{description}

For each of the above problems, we call $\bar{u}(\cdot)$ an optimal control, the corresponding solution $\bar{y}(\cdot)$ to (\ref{25}) an optimal trajectory, and $(\bar{y}(\cdot),\bar{u}(\cdot))$ an optimal pair.
Clearly, each of the above two optimal control problems can be viewed as an optimal control problem  with state constrained on a submanifold of the Euclidean space.

One of the central topics in control theory is to establish
necessary and sufficient conditions for optimal controls. As that in calculus,  one can derive the first-order necessary condition for optimal controls, as done in the classical monograph \cite{Pontryagin}, even for some situation of state constraints. Nevertheless, for some optimal control problems, it may well happen that the first-order necessary conditions
turn out to be trivial. In this case, the
first-order necessary condition cannot provide enough information for the theoretical
analysis and numerical computation, and therefore one needs to study the second (or even higher) order
optimality conditions for optimal controls.

There are many  works addressing the second order necessary and sufficient  conditions for optimal control problems in Euclidean spaces, such as \cite{bellj, bh, Cle-And, ft, gk, Gabasov73, g1, g, h2, Knobloc, k, l, m, mmp, Osmolovskii, w} (in which \cite{bellj, Cle-And,
Gabasov73, Knobloc, Osmolovskii} are five research monographs) and the references therein.  Among them, we mention several works, in which the controls are subjected to some restrictions:  Warga \cite{w}  considers a control system with the control set in a compact metric space. Frankowska and  Tonon \cite{ft} investigate the second order necessary conditions when a smooth endpoint constraint is presented and the control set is a closed subset of an Euclidean space. For the case that the control set is a general  metric space, Lou \cite{l} considers a control system  without state constraints.  For the case that the state (or mixed state-control) satisfies inequality or equality constraints, we refer to  \cite{fh, h, pz} and so on.

For the control system whose state is constrained to a manifold, there are also some literatures (e.g. \cite{as,bct,bChyb,UBBP,sl} and the references therein) devoted to the second order necessary and sufficient  conditions for optimal controls when the control set is an open subset of some manifold. Agrachev and  Sachkov  show a Legendre-type second order necessary condition for optimal controls, i.e., the Hessian of the corresponding Hamiltonian function with respect to the control variable is semi-negative definite (see \cite[Theorem 20.6, p. 300 and Proposition 20.11, p. 310]{as}); also, they give a strong Legendre-type second order sufficient condition for an optimal  controls for sufficiently short time intervals, i.e.  the Hessian of the corresponding Hamiltonian function with respect to the control variable is negative definite (see \cite[Proposition 20.12, p. 311 and Theorem 21.8, p. 347]{as}). On the basis of the second order optimality conditions, Bonnard,  Caillau and  Tr\'elat \cite{bct} present an algorithm  to compute the first conjugate time along
a smooth extremal curve, at which the trajectory fails to be optimal. In \cite[Chapter 4]{sl}, Sch\"attler and  Ledzewicz provide a ``state of art" account on the second order necessary and sufficient optimality conditions for a class of control-affine systems (Also, we refer to \cite{bChyb,UBBP} for relatively earlier works). Nevertheless, compared to the case of flat spaces, the picture of second optimality conditions in the curved spaces is, in our opinion, quite incomplete.  For example, it would be quite interesting to extend the related results in the flat spaces (say, that in the research monographs \cite{bellj, Cle-And,
Gabasov73, Knobloc, Osmolovskii} on second order optimality conditions) to the setting of curved spaces but many things remain to be done.

The main purpose of this paper is to investigate the second order necessary and sufficient  conditions for optimal control problems evolved on Riemannian manifolds. We distinguish the problems into two cases. The first case is when the endpoint is free, i.e., the case without the condition (\ref{503}). For this case, we consider the control systems with considerably general control sets.
By means of the needle variation technique, we obtain the second order Taylor's expansion of the cost functional with respect to the control variable  via the Riemannian geometric tools,  introduce a second order dual equation which depends on the curvature tensor of the manifold, and obtain  second order necessary conditions for optimal controls of integral and pointwise form respectively. In addition, we obtain the second order sufficient condition of integral form for local optimal controls. The second case is when the endpoint is fixed, i.e., the condition (\ref{503}) is imposed. For this more difficult case, as in \cite{as}, we assume that the control set is open, and  consider an extended control system involving an extra unknown variable associated to the cost functional. There is a well known relation between  the optimal controls of the original optimal control  problem and the above extended control  system (e.g. \cite[p. 179, Section 12.4]{as}): the terminal value  of the trajectory of the extended control  system corresponding to an optimal control  must be at the boundary of the attainable set of this control  system at the terminal time. Thanks to this  relation, Agrachev and Sachkov  \cite{as} observe that, the key point of finding second
 order optimality  conditions for the original optimal controls is to analyze the Hessian of the endpoint mapping of the extended control system with respect to the control variable. Following  this idea and employing again some tools from the global geometry, we compute an explicit expression of this Hessian (see (\ref{541})),  and obtain the desired second order optimality conditions, which, as the case of free endpoint, contain also the curvature tensor along the optimal trajectory. Note that the above mentioned explicit expression is absent in \cite{as}.

From the viewpoint of Riemannian geometry, it is quite natural that the curvature tensor appears in the second order optimality condition, see Synge's second variation formula for geodesics (e.g., \cite[Theorem  21, p. 158]{p1}). Nevertheless, one of the main difficulties for our control problems is to compute the second order Taylor's expansions (of the perturbed  trajectories), i.e., (\ref{32})--(\ref{37}) and (\ref{427}), with convenient and explicit expressions of the second order terms, which, in principle, should  involve the curvature tensor. To do this, quite different from the flat spaces, we need to establish a key geometric lemma, i.e, Lemma \ref{311}. This lemma enables us to obtain the second order variational equations, i.e., (\ref{477}) and (\ref{528}), whose solutions are exactly the desired expressions of the above mentioned second order terms. Clearly, the curvature tensor does appear in these two equations. It seems to us that Lemma \ref{311} has some independent interest and may be applied in other places. Another difficulty in our \textbf{Problem I} is
 how to introduce suitable dual equations such that the second order Taylor's expansion (\ref{54}) (of the cost functional) can be rewritten in terms of the dual variables. This is the key to obtain the corresponding second order optimality condition. It is not difficult to single out the first order dual equation  (\ref{64}) (because, at least formally, it is the same as that for the flat spaces), which is employed  to  rewrite the first order term in (\ref{54}) (see (\ref{498})). For the second order term (\ref{499}), there are some quadratic terms of the first order variation $X^\epsilon(\cdot)$ (see (\ref{85})).  Similar to the setting of flat spaces (\cite{l}), one needs to find some second order dual equation to ``cancel" these quadratic terms. Different from the flat setting and interestingly, the curvature tensor appears in our second order dual equation (\ref{70}).

It is worth mentioning that, compared to the corresponding results in the Euclidean spaces (\cite{bellj, bh, Cle-And, ft, gk, g1, g, h2, Knobloc, k, l, m, mmp, w}), the curvature tensor along the optimal trajectory  appears precisely in our results (Theorems \ref{66}--\ref{453}). This shows the very difference between the curved and the flat spaces from the viewpoint of optimal control theory.
Also, our results for \textbf{Problem II} (i.e., Theorems \ref{504}--\ref{453}) are different from \cite[Theorem 20.6, p. 300, Proposition 20.11, p. 310 and Theorem 20.16, p. 317]{as} and \cite[Proposition 20.12, p. 311 and Theorem 21.8, p. 347]{as}) at least in two aspects.
Firstly, our second order optimality  conditions depend on the curvature tensor, which does not appear explicitly in  \cite{as}.
    Secondly, our Theorem \ref{453} guarantees the local optimality of a control in any time interval, while the corresponding result in \cite{as} works only for a sufficiently short time interval.

In order to show the differences between our main results and that in the previous works, we shall provide two illustrative examples, i.e. Examples \ref{exl}--\ref{510}. Example \ref{exl} is about an optimal control problem on a non-compact manifold with negative curvature and with a discrete control set. We show that this problem fulfils all the assumptions in our main results. In Example \ref{510}, we apply one of our main results, i.e., Theorem \ref{504}, to the following famous geometric problem:
Given any two points on a Riemannian manifold,  find the necessary conditions of a locally shortest curve connecting these two points. We  conclude that this curve must be a geodesic and satisfies the second variation of energy, which is consistent with the corresponding results in Riemannian geometry (see \cite[p. 159]{p1}). However, the Legendre-type second order
optimality condition in \cite[Theorem 20.6, p. 300 and Proposition 20.11, p. 310]{as} becomes trivial for the same problem (see Remark \ref{rem4.1} for more detailed analysis).

Though geometric control theory is a huge subject and there exist many works on control theory via differential geometry, it seems to us that the more delicate tool, i.e. Riemannian geometry, has rarely been used to the study of ordinary differential equation control problems. As far as we know, this work is the first one which links explicitly the curvature tensor of the manifold to the optimality conditions. It seems us that this may stimulate further work  employing Riemannian geometry to the study of control problems.

The rest of this paper is organized as follows. In Section $2$, for the reader's convenience we recall first some basic notions and results on Riemannian manifolds, and also we show two lemmas which will be useful later. The main results of this paper are stated in Section 3. In Section 4, we give the above mentioned two examples. In Section 5, we obtain the first and second order variations of a trajectory of the system (\ref{25}), by means of two class of variation techniques, i.e., the needle variation and the classical variation.  Sections 6 is  devoted to proving  our main results in this work.

One of the main results in this paper, i.e., Theorem \ref{66}, has been announced in \cite{[57]} without proof.

\setcounter{equation}{0}
\section{Some preliminaries }
\def\theequation{2.\arabic{equation}}
\subsection{Exponential map}
For this part, we refer the readers to \cite[Chapter 3]{c} and \cite[Chapter 3]{wsy}.

A differentiable curve $\gamma(t)$ on $ M$ with $t\in[0,\alpha)$ (for some $\alpha>0$) is called a geodesic if it satisfies
\begin{equation}\label{102}
\nabla_{\dot{\gamma}(t)}\dot{\gamma}(t)=0,\quad t\in[0,\alpha).
\end{equation}
Let $x\in  M$ be fixed. For any $v\in T_x M$, there exists a unique geodesic $\gamma_v(\cdot)$ satisfying that $\gamma_v(0)=x$ and $\dot{\gamma}_v(0)=v$. Let $[0, \ell_v)$ be  the maximal  interval on which $\gamma_v(\cdot)$ is defined. Let $O_x\subset T_x  M$ be the set of vectors $v$ such that $\ell_v>1$. Then one can define  the  exponential  map as follows
$$\exp_x: O_x\to  M,\quad  \exp_x v=\gamma_v(1).$$
It can be shown that $O_x$ is a neighborhood of the origin $O\in T_x M$, and $\exp_x$ maps straight line
segments  in $T_x  M$ passing through the origin $O\in T_x M$  to geodesic
segments  in $ M$ passing through $x$. For any $v\in T_x M$, the differential of $\exp_x$ at $v$ is a linear map, denoted by
\begin{equation}\label{202}
d\exp_{x}|_v: \;\ T_vT_x M
\to T_{\exp_xv} M,\end{equation} where $T_vT_x M $ denotes the tangent space of the manifold $T_x M$ at the point $v\in T_x M$.

Given an $\epsilon>0$, write
 \begin{equation}\label{409}
B(O,\epsilon)\equiv\{v\in T_x M;\ \  |v|< \epsilon\}
\ \
\textrm{and} \ \
B_x(\epsilon)\equiv\{y\in  M;\ \  \rho(x,y)<\epsilon\}.
\end{equation}
We call
$
i(x)\equiv\sup\{\epsilon>0;\ \  \hbox{The map }\exp_x:B(O,\epsilon)\to B_x(\epsilon) \textrm{ is diffeomorphic}\}
$
the injectivity radius  at the point $x$ (e.g., \cite[p. 142]{p1}).

We list the following properties of the exponential map, which can be found in many books on Riemannian geometry (e.g. the proof of \cite[Proposition 2.9, p. 65]{c}).
\begin{lem}\label{317} For any $x\in  M$, the map $\exp_x$ is a local
diffeomorphism, whose  differential at the origin $O\in T_x M$ satisfies
\begin{equation}\begin{array}{c}\label{203}
d\exp_x|_O=d\exp_x^{-1}\Big|_x=\textrm{the identity operator on} \,\,T_x M.
\end{array}\end{equation}
Furthermore, for any $y\in  M$ with $\rho(x,y)<i(x)$, there exists a unique shortest piecewise smooth curve which is also a geodesic in $M$, connecting $x$ and $y$.
\end{lem}

\subsection{Parallel translation and tensors}
For the details of this part, we refer the readers to \cite[Chapter I and Chapter III ]{kn}, \cite[Chapter 2]{p1}, \cite[Chapter 1]{wsy} and \cite[Chaper 1]{h1}.

For any $x\in M$ and $r,s\in\N$, a multilinear map
 $$
F: \ \ \underbrace{T_x^*  M\times\cdots\times T_x^*  M}_{r\;\mbox{times}}\times \underbrace{T_x M\times\cdots\times T_x M}_{s\;\mbox{times}}\to \R
$$
is called a tensor of order $(r,s)$ at $x$. Denote by $\T_s^r(x)$ the tensor space of type $(r,s)$  at $x$.
A tensor field $\T$
of type $(r,s)$ on  $M$ is  a smooth assignment of a tensor $ \T(x)\in T_s^r(x)$ to each point $x$ of $ M$.
The norm of $\T$ at $x\in M$ is defined as follows:
 \begin{equation}\label{270}
\begin{array}{r}|\T(x)|=\sup\big\{\T(x)(Y_1,\cdots,Y_r,\lambda_1,\cdots,\lambda_s);\ \ Y_j\in T^*_x M,\lambda_l\in T_x M,\\[3mm]
 |Y_j|\leq 1,\, |\lambda_l|\leq 1, j=1,\cdots,r, l=1,\cdots,s\big\},\quad x\in M.\end{array}
\end{equation}
Denote by $\T_s^r(M)$ the set of all tensor fields of type $(r,s)$ over $M$.

Let $\gamma: [0,\ell]\to  M$ be a differentiable  curve  with $\gamma(0)=x\in M$, $\gamma(\ell)=y\in M$ and $\ell>0$.
Given a vector $v\in T_xM$, there exists a unique vector field $X$ along $\gamma$ satisfying
$$
\nabla_{\dot{\gamma}(s)}X=0,\qquad\forall\; s\in[0,\ell], \quad X(0)=v.
$$
The mapping $T_x M\ni v\mapsto X(\gamma(\ell))\in T_{y} M$ is a linear isometry between $T_x M$ and $T_{y} M
$. We call this map the parallel translation along the curve $\gamma$, and denote it by $L^{\gamma}_{xy}v$.
The parallel translation along the curve $\gamma$ enjoys the following property:
\begin{equation}\label{205}
\langle L_{xy}^{\gamma}X_1,L_{xy}^{\gamma}X_2\rangle=\langle X_1,X_2\rangle,\qquad\forall\; X_1, X_2\in T_x M.
\end{equation}
For any $\eta\in T_x^*M$, we define $L_{xy}^\gamma\eta\in T_y^*M$ by $L_{xy}^\gamma\eta(X)=\eta((L_{xy}^\gamma)^{-1}X)$ for any $X\in T_y M$.
One can extend the parallel translation of a vector  at $x\in M$ along the  curve $\gamma$ to a tensor $\T\in\T_s^r(x)$ by
$$
L_{xy}^\gamma \T(v_1,\cdots,v_r,\eta_1,\cdots,\eta_s)=\T((L_{xy}^\gamma)^{-1}v_1,\cdots,(L_{xy}^\gamma)^{-1}v_r,(L_{xy}^\gamma)^{-1}\eta_1,
\cdots,(L_{xy}^\gamma)^{-1}\eta_s),
$$
 for all $v_1,\cdots,v_r\in T_y^*M$ and $\eta_1,\cdots,\eta_s\in T_y M$.
From the above formula, one can get
\begin{equation}\label{82} \T(v_1,\cdots, v_r,\eta_1,\cdots,\eta_s)=L_{xy}^{\gamma}\T\Big(L^{\gamma}_{xy}v_1,\cdots, L^{\gamma}_{xy}v_r, L^{\gamma}_{xy}\eta_1,\cdots, L^{\gamma}_{xy}\eta_s\Big), \end{equation}
for any $v_1,\cdots,v_r\in T_x^* M$ and $\eta_1,\cdots,\eta_s\in T_x M$.  Especially, from  (\ref{270}), (\ref{205}) and (\ref{82}), it follows that
\begin{equation}
\label{268} |X|=|L_{xy}^{\gamma}X|,\quad |\lambda|=|L_{xy}^{\gamma}\lambda |,\qquad\forall\; X\in T_x M,\quad\lambda\in T_x^*  M.
\end{equation}

In particular, if $\rho(x,y)<\min\{i(x),i(y)\}$, according to Lemma \ref{317}, there is a unique shortest geodesic  $\gamma$ connecting  $x$ and $y$. In this case, we use $L_{xy}$ instead of $L_{xy}^{\gamma}$ for abbreviation.

Let $\T$ be a tensor field. Take any $v\in T_xM$. Let $\gamma$ be a smooth curve  such that $\gamma(0)=x$ and $\dot{\gamma}(0)=v$. Then the covariant derivative of a tensor field   (in terms  of parallel translation) is defined as follows (see \cite[p. 42]{h1}):
\begin{equation}\label{224}
\nabla_v\T=\lim_{t\to 0}\frac{1}{t}\Big(L_{\gamma(t)x}^\gamma \T(\gamma(t))-\T(x)\Big).
\end{equation}

Denote by $\nabla \T$ the covariant differential of $\T$, which is a tensor field of order $(r,s+1)$, and    is defined by (see \cite[p. 124]{kn})
\begin{equation}\label{100}
\nabla \T(Y_1,\cdots,Y_r,\lambda_1,\cdots,\lambda_s,Z)=\nabla_Z\T(Y_1,\cdots,Y_r,\lambda_1,\cdots,\lambda_s),
\end{equation} for all $Y_1,\cdots,Y_r\in T^*  M$ and $\lambda_1,\cdots,\lambda_s,Z\in T M$.
Applying \cite[Proposition 2.7, p. 123 and Proposition 2.10, p. 124 ]{kn}, one can get the following property: For each $Z\in T  M$, one has
\begin{equation}\label{99}\begin{array}{ll}
\nabla_Z \T(Y_1,\cdots,Y_r,\lambda_1,\cdots,\lambda_s)\\[2mm] =Z\Big(\T(Y_1,\cdots,Y_r,\lambda_1,\cdots,\lambda_s)\Big)-\T(\nabla_ZY_1,Y_2,\cdots,Y_r,\lambda_1,\cdots,\lambda_s)\\[2mm]
\quad-\cdots-\T(Y_1,\cdots,\nabla_ZY_r,\lambda_1,\cdots,\lambda_s)
-\T(Y_1,\cdots,Y_r,\nabla_Z\lambda_1,\lambda_2,\cdots,\lambda_s)-\cdots
\\[2mm] \quad-\T(Y_1,\cdots,Y_r,\lambda_1,\cdots,\nabla_Z\lambda_s),\quad\forall\; Y_1,\cdots,Y_r\in T^*  M,\,\,\lambda_1,\cdots,\lambda_s\in T M.
\end{array}
\end{equation}

In particular, a smooth function $f\in C^\infty( M)$ is a tensor of order $(0,0)$. $\nabla^2f$ is a tensor of order $(0,2)$. We call this tensor  the Hessian of the function $f$, which  is a symmetric tensor, and  can be computed by
\begin{equation}\label{223}
\nabla^2f(x)(X,Y)=Y(x) (X f)-(\nabla_{Y(x)}X)f,\qquad x\in  M,\ X,Y\in T M.
\end{equation} For a smooth function $h: M\times M\to \R$ of two arguments, we denote by $\nabla_i h$ the covariant derivative of $h$ with respect to  the $i^{th}$ argument with $i=1,2$, i.e., for $X\in T M$ and $(x_1,x_2)\in M\times M$,
\begin{equation}\label{476}
\langle \nabla_i h(x_1,x_2),X(x_i)\rangle= X(x_i)h(x_1,x_2).
\end{equation}
 Moreover, we define higher order derivatives of $h$ as follows: For $i,j=1,2$, $i\neq j$, any $(x_1,x_2)\in M\times M$ and $X,Y,Z\in T M$,
\begin{equation}\label{108}\begin{array}{ll}
\displaystyle\nabla_i\nabla_jh(x_1,x_2)(X,Y)\equiv Y(x_i)\Big(X(x_j)(h(x_1,x_2))\Big)=Y(x_i)(\langle \nabla_j h(x_1,x_2),X(x_j)\rangle);
\\[3mm] \displaystyle\nabla_i^2h(x_1,x_2)(X,Y)\equiv Y(x_i)\Big(X(x_i)h(x_1,x_2)\Big)-\nabla_{Y(x_i)}Xh(x_1,x_2);
\\[3mm] \displaystyle\nabla_i^2\nabla_j h(x_1,x_2)(X,Y,Z)\equiv\nabla_i^2(\langle X(x_j),\nabla_j h(x_1,x_2)\rangle)(Y,Z);
\\[3mm] \displaystyle\nabla_i\nabla_j^2h(x_1,x_2)(X,Y,Z)\equiv Z(x_i)\Big(\nabla_j^2h(x_1,x_2)(X,Y)\Big).
\end{array}\end{equation}

Let us recall the definition of the product of tensors (see \cite[p.22]{kn}):
 \begin{equation}\label{380}
 \T\otimes \K=X_1\otimes \cdots\otimes X_r\otimes Y_1\otimes\cdots\otimes Y_p\otimes \omega_1\otimes\cdots\otimes\omega_s\otimes \eta_1\otimes\cdots\otimes\eta_q\in \T_{s+q}^{r+p}(M),
 \end{equation}
 where $\T=X_1\otimes\cdots\otimes X_r\otimes\omega_1\otimes\cdots\otimes \omega_s\in \T_s^r(M)$, $\K=Y_1\otimes\cdots\otimes Y_p\otimes\eta_1\otimes\cdots\otimes\eta_q\in \T_q^p(M)$ and $r,s, p,q\in\N$. For any $i,j\in{I\!\!N}$, denote by ${\E}^{ij}$ the contraction of  the $i$th contravariant  index and the $j$th covariant  index, which is a linear mapping from $\T_s^r(M)$ to $\T_{s-1}^{r-1}(M)$ with $1\leq i\leq r$ and $1\leq j\leq s$, and is defined by ( see \cite[p.17]{h1})
\begin{equation}\label{76}\begin{array}{l}
{\E}^{ij}(X_1\otimes\cdots\otimes X_r\otimes \omega_1\otimes\cdots\otimes \omega_s)
\\=\omega_j(X_i)X_1\otimes\cdots\otimes X_{i-1}\otimes X_{i+1} \otimes\cdots\otimes X_r\otimes  \omega_1\otimes\cdots \otimes\omega_{j-1}\otimes\omega_{j+1}\otimes\cdots \otimes \omega_s,
\end{array}\end{equation}
 for all $X_1,\cdots X_r\in T M$ and $\omega_1,\cdots,\omega_s\in T^*M$.

\subsection{Two useful lemmas}

We begin with the following technical result (in which, (\ref{81})--(\ref{18}) and the first equality in (\ref{80}) can be found in \cite{d}).

\begin{lem}\label{17}
For any $x,y\in  M$ with $\rho(x,y)<\min\{i(x),i(y)\}$, $X,X_1,X_2\in T_x M$ and $Y\in T_y  M$, it holds that
\begin{eqnarray}
&&\label{80}|\displaystyle\exp_x^{-1}y|=|\exp_y^{-1}x|=\rho(x,y), \qquad \nabla_{X_1}L_{x\cdot}X=0,
\\[2mm] &&\displaystyle\label{81}\nabla_1\rho^2(x,y)=-2\exp_x^{-1}y,\qquad\nabla_2\rho^2(x,y)=-2\exp_y^{-1}x,
\\[2mm] &&\label{18}\displaystyle
L_{xy}exp_x^{-1}y=-exp_y^{-1}x, \qquad L_{xy}d_1\rho^2(x,y)=-d_1\rho^2(y,x),
\\[2mm] &&\label{92}\displaystyle
\nabla_1\nabla_2\rho^2(x,y)(Y,X)=-2\langle d\exp_y^{-1}|_xX,Y\rangle,
\\[2mm] &&\label{56}\displaystyle\langle d\exp_x^{-1}\Big|_yY,X\rangle=\langle d\exp_y^{-1}\Big|_x X,Y\rangle,
\\[2mm] &&\label{96}\displaystyle\nabla_1\nabla_2\rho^2(x,y)(Y,X)=-\nabla_1^2\rho^2(x,y)(L_{yx}Y,X)-\langle\nabla_1\rho^2(x,y),\nabla_X L_{y\cdot}Y\rangle,
\\[2mm] && \label{472}\displaystyle \nabla_1^2\rho^2(x,x)(X_1,X_2)=\nabla_2^2\rho^2(x,x)(X_1,X_2)=2\langle X_1,X_2\rangle,
\\[2mm] &&\label{260} \displaystyle\nabla_i\nabla_j^2\rho^2(x,x)=\nabla_i^2\nabla_j\rho^2(x,x)= \nabla_i^3\rho^2(x,x)=0,\quad i,j=1,2,\,\,i\neq j,
\end{eqnarray} where the notions $\nabla_1\nabla_2\rho^2$, $\nabla_i^2\nabla_j\rho^2$ and  $\nabla_i\nabla_j^2\rho^2$ with $i,j=1,2$ and $i\neq j$ are defined in (\ref{108}),  $\nabla_i^2\rho^2$ is the Hessian of $\rho^2$ with respect to the $i^{th}$ argument, $\nabla_i^3\rho^2$ is the covariant derivative of the Hessian $\nabla_i^2\rho^2(x,x)$ with respect to  the $i^{th}$ argument (see (\ref{100})),   and $d_i$ stands for the exterior derivative of a function on $M\times M$ with respect to the $i^{th}$ argument for $i=1,2$.
\end{lem}

\textbf{Proof:} \, We only prove the second equality in (\ref{80}), and (\ref{92})--(\ref{260}).

Let $\gamma_1$ be a radial geodesic satisfying  $\gamma_1(0)=x$  and $\dot{\gamma}_1(0)=X_1$. Then, $\gamma_1$ is the shortest geodesic connecting $\gamma_1(s)$ and $x$, provided that $s>0$ is small enough. By (\ref{224}), we have
$$\nabla_{X_1}L_{x\cdot}X=\lim_{s\to 0^+}\frac{1}{s}\Big(L_{\gamma_1(s)x}^{\gamma_1}L_{x\gamma_1(s)}X-X\Big)=0,
$$ which gives the second equality of (\ref{80}).

 Let $\gamma$ and $\beta$ be the curves satisfying
\begin{equation}\label{422}
\gamma(0)=x,\,\,\dot{\gamma}(0)=X,\qquad\beta(0)=y,\qquad\dot{\beta}(0)=Y.
\end{equation} By (\ref{108}) and noting (\ref{81}), we have
$$\begin{array}{lll}
\displaystyle\nabla_1\nabla_2\rho^2(x,y)(Y,X)&=&\displaystyle\frac{\partial}{\partial \theta}\frac{\partial}{\partial\tau}\rho^2(\gamma(\theta),\beta(\tau))|_{\theta=\tau=0}
= \displaystyle\frac{\partial}{\partial \theta}\langle \nabla_2\rho^2(\gamma(\theta),\beta(0)),\dot\beta(0)\rangle\Big|_{\theta=0}
\\[3mm] &=&\displaystyle-2\frac{\partial}{\partial\theta}\langle\exp_{y}^{-1}\gamma(\theta),Y\rangle\Big|_{\theta=0}
=\displaystyle-2\langle d\exp_y^{-1}|_xX,Y\rangle,
\end{array}$$ which implies (\ref{92}). Based on the above identity, (\ref{81}), (\ref{18}) and (\ref{205}), we have
$$\begin{array}{lll}
\displaystyle\nabla_1\nabla_2\rho^2(x,y)(Y,X)&=&\displaystyle-2\frac{\partial}{\partial\theta}\langle L_{y\gamma(\theta)}\exp_y^{-1}\gamma(\theta), L_{y\gamma(\theta)}Y\rangle\Big|_{\theta=0}
\\[3mm] &=&\displaystyle -\frac{\partial}{\partial\theta}\langle \nabla_1\rho^2(\gamma(\theta),y),L_{y\gamma(\theta)}Y\rangle\Big|_{\theta=0}
\\[3mm] &=&\displaystyle-\langle\nabla_X\nabla_1\rho^2(\cdot,y), L_{yx}Y\rangle-\langle\nabla_1\rho^2(x,y),\nabla_XL_{y\cdot}Y\rangle,
\end{array}
$$which, together with (\ref{100}), implies (\ref{96}).

To  prove (\ref{56}), we recall (\ref{422}) and the definition of the differential  of a differentiable  map, apply  (\ref{81}) and get
 $$\begin{array}{ll}
\displaystyle \langle d\exp_x^{-1}\Big|_yY,X\rangle&\displaystyle =\langle d\exp_x^{-1}\Big|_y\dot{\beta}(0),X\rangle
  =\langle\frac{d}{d\tau}\exp_x^{-1}\beta(\tau)\Big|_{\tau=0},X\rangle
 \\[3mm] &\displaystyle=-\frac{1}{2}\frac{d}{d\tau}\langle \nabla_1\rho^2(x,\beta(\tau)),\dot{\gamma}(0))\rangle\Big|_{\tau=0}
\\[3mm] &\displaystyle=-\frac{1}{2}\frac{d}{d\tau} \dot{\gamma}(0)\Big(\rho^2(\cdot,\beta(\tau))\Big)\Big|_{\tau=0}
=-\frac{1}{2}\frac{\partial}{\partial \tau}\frac{\partial}{\partial\theta}\rho^2(\gamma(\theta),\beta(\tau))\Big|_{\theta=\tau=0}.\end{array}
 $$
Similarly, we get
$
\langle d\exp_y^{-1}\Big|_x X,Y\rangle=-\frac{1}{2}\frac{\partial}{\partial \tau}\frac{\partial}{\partial\theta}\rho^2(\gamma(\theta),\beta(\tau))\Big|_{\theta=\tau=0}.
$ These two equalities imply (\ref{56}).

To prove  (\ref{472}), we choose the normal coordinates $\{x_1,\cdots, x_n\}$ at $x$ such that
\begin{equation}\label{261}
\Big\langle \frac{\partial}{\partial x_i}(x),\frac{\partial}{\partial x_j}(x)\Big\rangle =\delta_{ij},\qquad
\nabla_{\frac{\partial}{\partial x_i}(x)}\frac{\partial}{\partial x_j}=0,
\end{equation}
for $i,j=1,\cdots,n$.  From the above property, we can further deduce that
\begin{equation}\label{440}
\frac{\partial}{\partial x_i}\Big|_xg_{kl}(\cdot)=0,
\end{equation} with $g_{kl}(\cdot)\equiv \langle \frac{\partial}{\partial x_k}(\cdot),\frac{\partial}{\partial x_l}(\cdot)\rangle$, where $ i,k,l=1,\cdots,n$. By some computation in this system of local coordinates, we obtain (\ref{472}).

Then we go to the proof  of (\ref{260}). For any $y$ which is sufficiently close to $x$ and
 any vector field $F$, recalling (\ref{108}), we have
$$\begin{array}{ll}
\nabla_1^2\nabla_2\rho^2(y,x)(F,\cdot,\cdot)&=\nabla_1^2
\Big(\nabla_2\rho^2(y,x)(F)\Big)(\cdot,\cdot)
=-\nabla_1^2\Big(\nabla_1\rho^2(y,x)(L_{xy}F)\Big)(\cdot,
\cdot),
\end{array}$$
where we have used (\ref{82}), (\ref{81}) and (\ref{18}). Let us compute the right hand side of the above identity in local coordinates. By the property of the normal coordinates $\{x_1,\cdots, x_n\}$, we have $\rho^2(x,y)=|\exp_x^{-1}y|^2=\sum_{i=1}^nx_i^2(y)$. For any $z$ which is sufficiently close to $x$, denote
\begin{equation}\label{475}L_{zy}F=
\sum_{i=1}^na_i(z,y)\frac{\partial}{\partial x_i}(y),
\end{equation}
 where $a_i(z,y)$ is a function depending on $z$ and $y$.
Then,
$$
\nabla_1\rho^2(y,x)(L_{xy}F)=L_{xy}F\Big(\rho^2(\cdot,x)
\Big)=\sum_{i,j=1}^na_j(x,y)\frac{\partial}{\partial x_j}(y)(x_i^2(\cdot))=2\sum_{i=1}^na_i(x,y)x_i(y).
$$
For any vector fields $X_j(\cdot)=\sum_{k=1}^nX_j^k(\cdot)\frac{\partial}{\partial x_k}(\cdot)$ with $j=1,2$, noting (\ref{223}), (\ref{18}) and (\ref{261}), we obtain that
$$
\nabla_1^2\nabla_2\rho^2(x,x)(F,X_1,X_2)=-2\sum_{k,l,m=1}^n
X_2^k(x)X_1^l(x)\frac{\partial^2}{\partial x_k\partial x_l}\Big|_x\Big(a_m(x,\cdot)x_m(\cdot)\Big).
$$
We claim that
\begin{equation}\label{263}\sum_{k,l,m=1}^nX_2^k(x)X_1^l(x)\frac{\partial^2}{\partial x_k\partial x_l}\Big|_x\Big(a_m(x,\cdot)x_m(\cdot)\Big)=0.\end{equation}
In fact, according to (\ref{261}),  (\ref{440}) and  (\ref{475}) we have
$$\begin{array}{lll}
\frac{\partial }{\partial x_k}\Big|_x a_m(x,\cdot)
&=&\sum_{i=1}^n \frac{\partial}{\partial x_k}\Big|_x\Big(a_i(x,\cdot)g_{im}(\cdot)\Big)=\frac{\partial}{\partial x_k}\Big|_x\langle L_{x\cdot}F,\frac{\partial}{\partial x_m}(\cdot)\rangle
\\&=&\Big\langle \nabla_{\frac{\partial}{\partial x_k}|_x}L_{x\cdot}F,\frac{\partial}{\partial x_m}\Big|_x\Big\rangle,
\end{array}$$
 where $k=1,\cdots, n$.
By the second equality of (\ref{80}), we have  $\frac{\partial }{\partial x_k}\Big|_x a_m(x,\cdot)=0 $, which implies $\nabla_1^2\nabla_2\rho^2(x,x)=0$. Using the same method, we can show that $\nabla_2\nabla_1^2\rho^2(x,x)=0$, $\nabla_1\nabla_2^2\rho^2(x,x)=0$ and $\nabla_2^2\nabla_1\rho^2(x,x)=0$.
We now prove that $\nabla_1^3\rho^2(x,x)=0$. For this purpose, we take any $y(\in M)$ which is sufficiently close to $x$. For any $Y\in T_y M$ and $j,k=1,\cdots,n$, letting $\frac{\partial }{\partial x_k}(x)$ act on both sides of
 (\ref{96}) with $X=\frac{\partial}{\partial x_j}$, we obtain that
$$\begin{array}{lll}
&&\displaystyle\nabla_1^2\nabla_2\rho^2(x,y)(Y,\frac{\partial}{\partial x_j},\frac{\partial}{\partial x_k})
\\[3mm] &=&\displaystyle-\nabla_1^3\rho^2(x,y)(L_{yx}Y,\frac{\partial}{\partial x_j},\frac{\partial}{\partial x_k})-\nabla_1^2\rho^2(x,y)(\nabla_{\frac{\partial}{\partial x_k}}L_{y\cdot}Y,\frac{\partial}{\partial x_j})
\\[3mm]&&\displaystyle-\nabla_1^2\rho^2(x,y)(\nabla_{\frac{\partial}{\partial x_j}}L_{y\cdot}Y,\frac{\partial}{\partial x_k})
-\langle \nabla_1\rho^2(x,y),\nabla_{\frac{\partial}{\partial x_k}}\nabla_{\frac{\partial}{\partial x_j}}L_{y\cdot}Y\rangle,
\end{array}$$
where we have used  (\ref{99}) and  (\ref{261}). Letting  $y$ approach to $x$, and using  (\ref{80}), (\ref{81}) and $\nabla_1^2\nabla_2\rho^2(x,x)=0$, we conclude that $\nabla_1^3\rho^2(x,x)=0$. Similarly, we can get $\nabla_2^3\rho^2(x,x)=0$ in the same way.
$\Box$

\medskip

Denote by $[X,Y]\equiv XY-YX$ the Lie bracket of vector fields $X$ and $Y$.
Denote by  $R$  the curvature tensor (of $(M,g)$), which is a correspondence that associates to every pair $X,Y\in T M$ a mapping $R(X,Y): T M\to T M$ given by $$R(X,Y)Z=\nabla_X\nabla_Y Z-\nabla_Y\nabla_X Z-\nabla_{[X,Y]}Z,\qquad \forall\; Z\in T M.$$
We write
$$
R(X,Y,Z,W)=\langle R(X,Y)Z,W\rangle,\qquad \forall\; X,Y,Z,W\in T M.
$$
For any $x\in M$ and $X(x),Y(x)\in T_xM$ with $X(x)\not\parallel Y(x)$, the sectional curvature $sec(X(x),Y(x))$ (of the plane spanned by $X(x)$ and $Y(x)$) is given by
\begin{equation}\label{415}
sec(X(x),Y(x))\equiv \frac{R(X,Y,Y,X)(x)}{|X(x)\wedge Y(x)|^2}.
\end{equation}

The following result will play a key role in the sequel.

\begin{lem}\label{311} Let $X,V, F\in T M$. Then, for any  $x\in M$, it holds that \footnote{${}^1$The left hand side of (\ref{312}) is defined as follows:
For $y\in M$ closing enough to $x$, we view $L_{\cdot y}F(\cdot): B_x(i(x))\to T_yM$ as a map, where $B_x(i(x))$ is defined in (\ref{409}). The differential of $L_{\cdot y}F(\cdot)$ at $x$ is a linear map, denoted by $d_x(L_{xy}F(x)): T_xM\to T_{L_{xy}F(x)}T_yM$, where $T_{L_{xy}F(x)}T_yM$  is isomorphic to $T_yM$. Therefore, we have  $d_x(L_{xy}F(x))V(x)\!\in \! T_yM$. Letting $y$ vary around $x$, we get a vector field $d_x(L_{x\cdot}F(x))\!V(x)$ around $x$. Hence, $\nabla_{X(x)}\Big(d_x(L_{x\cdot}F(x))V(x)\Big)$ is the covariant derivative of vector field $d_x(L_{x\cdot}F(x))V(x)$ relative to $X(x)$.}
\begin{eqnarray}\label{312}
\nabla_{X(x)}\Big(d_x(L_{x\cdot}F(x))V(x)\Big)=\frac{1}{2}  R(X,V) F(x)\,{}^1;
\\\label{496} \nabla_2^3\nabla_1\rho^2(x,x)(X,F,V,V)+\nabla_2^2\nabla_1^2
\rho^2(x,x)(X,F,V,V)=2R(X,V,F,V)(x).
\end{eqnarray}
\end{lem}

\textbf{Proof:}  \, We will prove (\ref{312}) and (\ref{496}) by means of suitable local coordinates.
Let $\{x_1,\cdots, x_n\}$ be the normal coordinates in a neighborhood $O\subset M$ around $x$, which satisfies (\ref{261}) and (\ref{440}). Denote by $(g_{ij})(\cdot)\equiv \Big(\langle \frac{\partial}{\partial x_i},\frac{\partial}{\partial x_j}\rangle(\cdot)\Big)$ and $\Gamma_{ij}^k$ ($i,j,k=1,\cdots,n$) the metric matrix and the Christoffel symbols, respectively.
Take any $m,j \in \{1,\cdots, n\}$. Let $x_j(\cdot)$ be the coordinate curve such that   $\frac{\partial}{\partial x_j}(x)=\dot x_j(0)$ (See \cite[p. 8]{c}). Recalling  (\ref{476}), (\ref{475}) and the definition of the differential of a differentiable map, for any $y\in O$, we have
$$
d_x(L_{xy}F(x))\frac{\partial}{\partial x_j}(x)=\frac{\partial }{\partial s}\Big|_0L_{x_j(s)y}F(x_j(s))=\sum_{i=1}^n\langle \nabla_1a_i(x,y),\frac{\partial}{\partial x_j}(x)\rangle\frac{\partial}{\partial x_i}(y),
$$
where $\nabla_la_i(\cdot,\cdot)$ (with $l=1,2$ and $i=1,\cdots, n$) is defined in (\ref{476}). Then, recalling (\ref{261}) and (\ref{108}), we have
\begin{equation}\label{443}
\nabla_{\frac{\partial}{\partial x_m}(x)}\Big(d_x(L_{x\cdot}F(x))\frac{\partial}{\partial x_j}(x)\Big)=\sum_{i=1}^n\nabla_2\nabla_1a_i(x,x)\Big(\frac{\partial}{\partial x_j},\frac{\partial}{\partial x_m}\Big)\frac{\partial}{\partial x_i}(x).
\end{equation}

Clearly, to obtain (\ref{312}), we need to compute $\nabla_2\nabla_1a_i(x,x)\Big(\frac{\partial}{\partial x_j},\frac{\partial}{\partial x_m})$ for $i,j,m=1,\cdots,n$. To this end, we first show some properties of short geodesics. Given $z,y\in O$,
 write $V\equiv \exp_z^{-1}y$ and  $\gamma_{zy}(s)\equiv \exp_zsV$ with $s\in [0,1]$. Then, $\gamma_{zy}(\cdot)$ is the shortest geodesic connecting $z$ and $y$, and satisfies
$
\gamma_{zy}(0)=z,\,\,\gamma_{zy}(1)=y$ and $\dot{\gamma}_{zy}(0)=\exp_z^{-1}y
$.
We assume that the local expression of  $\gamma_{zy}(s)$ with $s\in[0,1]$ in the system of coordinates $(O,x_1,\cdots,x_n)$ is as follows:
$$(\gamma_{zy}^1(s),\cdots,\gamma_{zy}^n(s)).$$
Then, $\dot\gamma_{zy}(s)=\sum_{j=1}^n\dot\gamma_{zy}^j(s)
\frac{\partial}{\partial x_j}(\gamma_{zy}(s))$. In particular, by Lemma \ref{17}, recalling the definition of a geodesic and using the local expression of $\gamma_{zy}(\cdot)$, we have
 \begin{equation}\label{495}
 \dot\gamma_{zy}(1)=\sum_{i=1}^n\dot\gamma_{zy}^i(1)
 \frac{\partial}{\partial x_i}(y)=L_{zy}\dot\gamma_{zy}(0)=L_{zy}\exp_z^{-1}y
 =-\frac{1}{2}L_{zy}\nabla_1\rho^2(z,y)
 =\frac{1}{2}\nabla_2\rho^2(z,y),
 \end{equation} which implies that $\dot\gamma_{zy}^i(1)$ with $i=1,\cdots,n$ are determined by the endpoints of geodesic $\gamma_{zy}(\cdot)$. Thus,
  we may view $\dot\gamma_{zy}^1(1),  \cdots, $ $ \dot\gamma_{zy}^n(1)$ as functions of the first argument $z$ and the second argument $y$.

We claim that (Recall (\ref{476}) and (\ref{108}) for the notations)
  \begin{eqnarray}\label{446}
&\displaystyle\langle\nabla_2\dot\gamma_{xy}^j(1),\frac{\partial}{\partial x_\zeta}(y)\rangle=\delta_j^\zeta+o(1),\qquad&\displaystyle\langle\nabla_1\dot\gamma_{xy}^j(1),\frac{\partial}{\partial x_\zeta}(x)\rangle=-\delta_j^\zeta+o(1),
\\[3mm]\label{447}&\displaystyle\!\!\!\!\!\!\!\!\!\!\!\!\!\!\!\!\!\!\!\!\!\!\!\!\!\!\!\!\nabla_1\nabla_2\dot\gamma_{xy}^j(1)=o(1),\qquad&\displaystyle\nabla_2^2\dot\gamma_{xy}^j(1)=o(1),
\end{eqnarray}
for $\zeta,j=1,\cdots,n$, where $\delta_j^\zeta$ is the usual Kronecker symbol, $o(1)$ is a tensor  of suitable order satisfying $\lim_{y\to x}o(1)=0$.
Indeed,  for any $j=1,\cdots,n$, it follows from (\ref{495}) that
 \begin{equation}\label{448}
 \frac{1}{2}\langle \nabla_2\rho^2(z,y),\frac{\partial}{\partial x_j}(y)\rangle=\sum_{i=1}^n\langle\dot\gamma_{zy}^i(1)
 \frac{\partial}{\partial x_i}(y),\frac{\partial}{\partial x_j}(y)\rangle.
 \end{equation}
 Letting $\frac{\partial}{\partial x_\zeta}(y)$ act on the above identity, by (\ref{99}), we get
 \begin{equation}\label{445}\begin{array}{l}
 \displaystyle\frac{1}{2}\nabla_2^2\rho^2(z,y)
 (\frac{\partial}{\partial x_j},\frac{\partial}{\partial x_\zeta})+\frac{1}{2}\langle \nabla_2\rho^2(z,y),\nabla_{\frac{\partial}{\partial x_\zeta}(y)}\frac{\partial}{\partial x_j}\rangle
 \\[3mm] \displaystyle=\sum_{i=1}^n\Big(\langle \nabla_2\dot\gamma_{zy}^i(1),\frac{\partial}{\partial x_\zeta}(y)\rangle g_{ij}(y)+\dot\gamma_{zy}^i(1)\frac{\partial}{\partial x_\zeta}(y)g_{ij}\Big).
\end{array} \end{equation}
According to (\ref{80}), (\ref{81}), (\ref{472}), (\ref{260}), (\ref{261}), (\ref{440}) and the fact that
\begin{equation}\label{473}\dot\gamma_{xx}^i(1)=0,\quad i=1,\cdots,n,
\end{equation}
 we obtain the first identity of (\ref{446}) by letting $z=x$ in (\ref{445}). In a similar way, we obtain the second identity of (\ref{446}) by letting $\frac{\partial}{\partial x_\zeta}(z)$ act on (\ref{448}) at $z=x$ with $\zeta=1,\cdots,n$.
  For $m=1,\cdots,n$, letting $\frac{\partial}{\partial x_m}(z)$ act on (\ref{445}), via (\ref{99}), we get
  $$\begin{array}{c}
  \displaystyle\frac{1}{2}\nabla_1\nabla_2^2\rho^2(z,y)
  \Big(\frac{\partial}{\partial x_j},\frac{\partial}{\partial x_\zeta},\frac{\partial}{\partial x_m}(z)\Big)+\frac{1}{2}\nabla_1\nabla_2\rho^2
  (z,y)\Big(\nabla_{\frac{\partial}{\partial x_\zeta}(y)}\frac{\partial}{\partial x_j},\frac{\partial}{\partial x_m}(z)\Big)
  \\[3mm]\displaystyle =\sum_{i=1}^n\Big(\nabla_1\nabla_2\dot
  \gamma_{zy}^i(1)(\frac{\partial}{\partial x_\zeta}(y),\frac{\partial}{\partial x_m}(z))g_{ij}(y)+\langle\nabla_1\dot
  \gamma_{zy}^i(1),\frac{\partial}{\partial x_m}(z)\rangle \frac{\partial}{\partial x_\zeta}(y)g_{ij}\Big).
\end{array}  $$
By Lemma \ref{17}, (\ref{261}),  (\ref{440}) and (\ref{446}), we obtain the first equality of (\ref{447}) by taking $z=x$ in the above identity. Similarly, we obtain the second equality of (\ref{447})  by letting $\frac{\partial}{\partial x_m}(y)$ ($m=1,\cdots,n$) act on (\ref{445}) with $z=x$.

Recalling the definition of parallel translation and (\ref{475}), for any $s\in[0,1]$, we have
$$\begin{array}{lll}
0&=&\displaystyle\nabla_{\dot\gamma_{zy}(s)}L_{z\cdot}F(z)
=\sum_{i,\zeta=1}^n\dot\gamma_{zy}^i(s)
\nabla_{\frac{\partial}{\partial x_i}(\gamma_{zy}(s))}\Big(a_\zeta(z,\cdot)\frac{\partial}{\partial x_\zeta}\Big)
\\[3mm] &=&\displaystyle\sum_{\zeta=1}^n\Big(\sum_{i=1}^n\dot
\gamma_{zy}^i(s)\langle\nabla_2a_\zeta(z,\gamma_{zy}(s)),
\frac{\partial}{\partial x_i}\rangle
\\[3mm] &&\displaystyle+\sum_{i,\eta=1}^n\dot\gamma_{zy}^i(s)
a_\eta(z,\gamma_{zy}(s))\Gamma_{i\eta}^\zeta(\gamma_{zy}
(s))\Big)\frac{\partial}{\partial x_\zeta}(\gamma_{zy}(s)),
\end{array}$$
which indicates that
$$
\sum_{i=1}^n\dot\gamma_{zy}^i(1)\langle\nabla_2
a_\zeta(z,y),\frac{\partial}{\partial x_i}\rangle+\sum_{i,\eta=1}^n\dot\gamma_{zy}^i(1)a_\eta(z,
y)\Gamma_{i\eta}^\zeta(y)=0, \qquad\zeta=1,\cdots,n.
$$
Letting  $\frac{\partial }{\partial x_m}(y)$ act on both sides of the above identity, via (\ref{99}), we obtain that
\begin{equation}\label{460}\begin{array}{l}
\displaystyle\sum_{i=1}^n\Big(\langle\nabla_2\dot\gamma_{z
y}^i(1),\frac{\partial}{\partial x_m}(y)\rangle\langle \nabla_2a_\zeta(z,y),\frac{\partial}{\partial x_i}(y)\rangle+\dot\gamma_{zy}^i(1)\nabla_2^2
a_\zeta(z,y)(\frac{\partial}{\partial x_i},\frac{\partial}{\partial x_m})
\\[3mm]\displaystyle\quad +\dot\gamma_{zy}^i(1)\langle \nabla_2a_\zeta(z,y),\nabla_{\frac{\partial}{\partial x_m}(y)}\frac{\partial}{\partial x_i}\rangle\Big)+\sum_{i,\eta=1}^n\Big(\langle\nabla_2
\dot\gamma_{zy}^i(1),\frac{\partial}{\partial x_m}(y)\rangle a_\eta(z,y)\Gamma_{i\eta}^\zeta(y)
\\[3mm]\displaystyle\quad+\dot\gamma_{zy}^i(1)
\langle\nabla_2a_\eta(z,y),\frac{\partial}{\partial x_m}(y)\rangle\Gamma_{i\eta}^\zeta(y)
+\dot\gamma_{zy}^i(1)a_\eta(z,y)\frac{\partial}{\partial x_m}(y)\Gamma_{i\eta}^\zeta\Big)\\[3mm]\displaystyle
=0.
\end{array}\end{equation}
Letting $\frac{\partial}{\partial x_j}(z) $ act on the above identity at $z=x$, we get
$$\begin{array}{l}
\displaystyle\sum_{i=1}^n\langle\nabla_1\dot\gamma_{xy}^i(1),\frac{\partial}{\partial x_j}(x)\rangle\nabla_2^2a_\zeta(x,y)\Big(\frac{\partial}{\partial x_i},\frac{\partial}{\partial x_m}\Big)
 +\sum_{i,\eta=1}^n\langle\nabla_1\dot\gamma_{xy}^i(1),\frac{\partial}{\partial x_j}(x)\rangle a_\eta(x,y)\frac{\partial}{\partial x_m}(y)\Gamma_{i\eta}^\zeta
 \\[3mm]\displaystyle\quad+\sum_{i=1}^n\Big(\langle\nabla_2\dot\gamma_{xy}^i(1),\frac{\partial}{\partial x_m}(y)\rangle\nabla_1\nabla_2a_\zeta(x,y)\Big(\frac{\partial}{\partial x_i}(y),\frac{\partial}{\partial x_j}(x)\Big) +o(1)\\[3mm]\displaystyle=0,
\end{array}$$
where we have used (\ref{99}),  (\ref{108}), (\ref{261}), (\ref{440}), (\ref{447}) and (\ref{473}). Letting $y$ approach to $x$, we obtain
\begin{equation}\label{449}
\nabla_1\nabla_2a_\zeta(x,x)(\frac{\partial}{\partial x_m}(x),\frac{\partial}{\partial x_j}(x))=\nabla_2^2a_\zeta(x,x)(\frac{\partial}{\partial x_j},\frac{\partial}{\partial x_m})+\sum_{\eta=1}^na_\eta(x,x)\frac{\partial}{\partial x_m}(x)\Gamma_{j\eta}^\zeta,
\end{equation}
where we have used  (\ref{446}).
For $j=1,\cdots,n$, letting $\frac{\partial}{\partial x_j}(y)$ act on (\ref{460}) with $z=x$, we get
$$\begin{array}{l}
\displaystyle\sum_{i=1}^n\Big(\!\langle\nabla_2\dot\gamma_{xy}^i(1),\!\frac{\partial}{\partial x_m}(y)\rangle\nabla_2^2a_\zeta(x,y)(\frac{\partial }{\partial x_i},\!\frac{\partial}{\partial x_j})\!+\!\langle\nabla_2\dot\gamma_{xy}^i(1),\!\frac{\partial}{\partial x_j}(y)\rangle\nabla_2^2a_\zeta(x,y)(\frac{\partial}{\partial x_i},\!\frac{\partial}{\partial x_m})\!\Big)
\\[3mm]\displaystyle\quad +\!\sum_{i,\eta=1}^n\!\Big(\!\langle\nabla_2\dot\gamma_{xy}^i(1),\frac{\partial}{\partial x_m}(y)\rangle a_\eta(x,y)\frac{\partial}{\partial x_j}(y)\Gamma_{i\eta}^\zeta\!+\!\langle\nabla_2\dot\gamma_{xy}^i(1),\frac{\partial}{\partial x_j}(y)\rangle a_\eta(x,y)\frac{\partial}{\partial x_m}(y)\Gamma_{i\eta}^\zeta\!\Big)\\[3mm]\quad+ o(1)=0,
\end{array}$$
where we have used (\ref{99}),  (\ref{108}), (\ref{261}), (\ref{440}), (\ref{447}) and (\ref{473}). Letting $y$ approach to $x$, we obtain
$$\begin{array}{l}
\displaystyle\nabla_2^2a_\zeta(x,x)(\frac{\partial}{\partial x_m},\frac{\partial}{\partial x_j})+\nabla_2^2a_\zeta(x,x)(\frac{\partial}{\partial x_j},\frac{\partial}{\partial x_m})
\\[3mm]\displaystyle\quad+\sum_{\eta=1}^na_\eta(x,x)\Big(\frac{\partial}{\partial x_j}(x)\Gamma_{m\eta}^\zeta+\frac{\partial}{\partial x_m}(x)\Gamma_{j\eta}^\zeta\Big)=0,
\end{array}$$where we have used (\ref{446}). Since $\nabla_2^2 a_\eta(x,x)$ is the Hessian of the function $a_\eta(x,\cdot)$ with respect to  the second argument at $x$, it is a symmetric tensor of order $(0,2)$. Therefore, we have
\begin{equation}\label{497}
\nabla_2^2a_\zeta(x,x)(\frac{\partial}{\partial x_j},\frac{\partial}{\partial x_m})=-\frac{1}{2}\sum_{\eta=1}^na_\eta(x,x)\Big(\frac{\partial}{\partial x_j}(x)\Gamma_{m\eta}^\zeta+\frac{\partial}{\partial x_m}(x)\Gamma_{j\eta}^\zeta\Big).
\end{equation}
Recalling (\ref{443}) and inserting the above identity into (\ref{449}), we can get (\ref{312}) with $X=\frac{\partial}{\partial x_m}$ and $V=\frac{\partial}{\partial x_j}$,
where we have used the expression of the curvature tensor in local coordinates (See  \cite[p. 41]{p1}), the property of the curvature tensor (See \cite[p. 33]{p1}), (\ref{261}),  and the definition of $a_j(x,x)$ with $j=1,\cdots,n$.
It is easy to check that $\nabla_{X(x)}\Big(d_x(L_{x\cdot}F(x))V(x)\Big)$ with $X,V\in T M$ is multi-linear with respect to $X$ and $V$. Thus, the proof of (\ref{312}) is completed.

To obtain (\ref{496}), for any $y\in O$, we apply  (\ref{96})  and get
$$\begin{array}{ll}
&\nabla_1\nabla_2\rho^2(x,y)( F(y),\frac{\partial}{\partial x_j}(x))
\\=&-\nabla_1^2\rho^2(x,y)(L_{yx} F(y),\frac{\partial}{\partial x_j}(x))-\Big\langle\nabla_1\rho^2(x,y),\nabla_{\frac{\partial}{\partial x_j}(x)}\Big(L_{y\cdot} F(y)\Big)\Big\rangle.
\end{array}$$
Letting $\frac{\partial}{\partial x_k}(y)$ ($1\leq k\leq n$) act on the above identity, via (\ref{99}), we can get
$$\begin{array}{l}
\nabla_1\nabla_2^2\rho^2(x,y)( F(y),\frac{\partial}{\partial x_k}(y),\frac{\partial}{\partial x_j}(x))+\nabla_1\nabla_2\rho^2(x,y)(\nabla_{\frac{\partial}{\partial x_k}(y)} F,\frac{\partial}{\partial x_j}(x))
\\ \quad+\nabla_2\nabla_1^2\rho^2(x,y)(L_{yx} F(y),\frac{\partial}{\partial x_j}(x),\frac{\partial}{\partial x_k}(y))+\nabla_1^2\rho^2(x,y)\Big(d_y(L_{yx} F(y))\frac{\partial}{\partial x_k}(y),\frac{\partial}{\partial x_j}(x)\Big)
\\ = -\nabla_2\nabla_1\rho^2(x,y)\Big(\nabla_{\frac{\partial}{\partial x_j}(x)}(L_{y\cdot} F(y)),\frac{\partial}{\partial x_k}(y)\Big)
\\ \quad-\Big\langle\nabla_1\rho^2(x,y),\nabla_{\frac{\partial}{\partial x_j}(x)}\Big(d_y(L_{y\cdot} F(y))\frac{\partial}{\partial x_k}(y)\Big)\Big\rangle.
\end{array}$$
Moreover, letting $\frac{\partial}{\partial x_k}(y)$ act on the above identity, as $y$ approaches to $x$, we can get
\begin{equation}\label{42}\begin{array}{l}
\nabla_1\nabla_2^3\rho^2(x,x)( F,\frac{\partial}{\partial x_k},\frac{\partial}{\partial x_k},\frac{\partial}{\partial x_j})+\nabla_1\nabla_2\rho^2(x,x)\Big(\nabla_{\frac{\partial}{\partial x_k}(x)}\nabla_{\frac{\partial}{\partial x_k}} F,\frac{\partial}{\partial x_j}(x)\Big)
\\ \quad +\nabla_2^2\nabla_1^2\rho^2(x,x)( F,\frac{\partial}{\partial x_j},\frac{\partial}{\partial x_k},\frac{\partial}{\partial x_k})
\\ \quad+\nabla_1^2\rho^2(x,x)\Big(d_y\Big(d_y(L_{yx} F(y))\frac{\partial}{\partial x_k}(y)\Big)\frac{\partial}{\partial x_k}(y)|_{y=x},\frac{\partial}{\partial x_j}(x)\Big)
\\ =-2\nabla_2\nabla_1\rho^2(x,x)\Big(\nabla_{\frac{\partial}{\partial x_j}(x)}\Big(d_y(L_{y\cdot} F(y))_{y=x}\frac{\partial}{\partial x_k}(x)\Big),\frac{\partial}{\partial x_k}(x)\Big),
\end{array}\end{equation}
where (\ref{80}), (\ref{81}) and (\ref{260}) are used.

We claim that
\begin{equation}\label{43}\begin{array}{l}
\nabla_1\nabla_2\rho^2(x,x)\Big(\nabla_{\frac{\partial}{\partial x_k}(x)}\nabla_{\frac{\partial}{\partial x_k}} F,\frac{\partial}{\partial x_j}(x)\Big)
\\ \quad+\nabla_1^2\rho^2(x,x)\Big(d_y\Big(d_y(L_{yx} F(y))\frac{\partial}{\partial x_k}(y)\Big)\frac{\partial}{\partial x_k}(y)|_{y=x},\frac{\partial}{\partial x_j}(x)\Big)=0.
\end{array}\end{equation}
In fact, by   (\ref{475}), we may rewrite $F(y)$ as
$
 F(y)=L_{yy} F(y)=\sum_{\eta=1}^na_\eta(y,y)\frac{\partial}{\partial x_\eta}(y),
$
and  compute $\nabla_{\frac{\partial}{\partial x_k}(y)} F$ as follows:
$$\begin{array}{ll}
\nabla_{\frac{\partial}{\partial x_k}(y)} F=&\sum_{\eta=1}^n\Big(\langle \nabla_1 a_\eta(y,y),\frac{\partial}{\partial x_k}(y)\rangle+\langle\nabla_2a_\eta(y,y),\frac{\partial}{\partial x_k}(y)\rangle\Big)\frac{\partial}{\partial x_\eta}(y)
\\ &+\sum_{\eta,\zeta=1}^na_\eta(y,y)\Gamma_{k\eta}^\zeta\frac{\partial}{\partial x_\zeta}(y).
\end{array}$$
Applying Lemma \ref{17}, and by (\ref{261}) and (\ref{440}), we have
\begin{equation}\label{44}\begin{array}{ll}
&\nabla_{\frac{\partial}{\partial x_k}(y)}\nabla_{\frac{\partial}{\partial x_k}} F
\\=&\sum_{\eta=1}^n\Big(\nabla_1^2a_\eta(y,y)(\frac{\partial}{\partial x_k},\frac{\partial}{\partial x_k})+2\nabla_1\nabla_2a_\eta(y,y)(\frac{\partial}{\partial x_k},\frac{\partial}{\partial x_k})
\\ &+\nabla_2^2a_\eta(y,y)(\frac{\partial}{\partial x_k},\frac{\partial}{\partial x_k})\Big)\frac{\partial}{\partial x_\eta}(y)
+\sum_{\eta,\zeta=1}^na_\eta(y,y)\frac{\partial}{\partial x_k}(y)(\Gamma_{k\eta}^\zeta) \frac{\partial}{\partial x_\zeta}(y)+o(1),
\end{array}\end{equation}
where $\lim\limits_{y\to x}o(1)=0$.
On the other hand, to figure out $d_y\Big(d_y(L_{yx} F(y))\frac{\partial}{\partial x_k}(y)\Big)\frac{\partial}{\partial x_k}(y)$ with $y$ closing to $x$, we choose a geodesic  $\gamma(\cdot)$ on $M$ satisfying $\gamma(0)=y$ and $\dot\gamma(0)=\frac{\partial}{\partial x_k}(y)$. Then,
$$
d_y(L_{yx}F(y))\frac{\partial}{\partial x_k}(y)=\frac{\partial}{\partial s}\Big|_0\sum_{i=1}^na_i(\gamma(s),x)\frac{\partial}{\partial x_i}(x)=\sum_{i=1}^n\langle\nabla_1a_i(y,x),\frac{\partial}{\partial x_k}(y)\rangle \frac{\partial}{\partial x_i}(x),
$$
and
\begin{equation}\label{46}\begin{array}{lll}
d_y\Big(d_y(L_{yx}F(y))\frac{\partial}{\partial x_k}(y)\Big)\frac{\partial}{\partial x_k}(y)&=&\frac{\partial}{\partial s}\Big|_0\sum_{i=1}^n\langle\nabla_1 a_i(\gamma(s),x),\frac{\partial}{\partial x_k}(\gamma(s))\rangle\frac{\partial}{\partial x_i}(x)
\\ &=&\sum_{i=1}^n\nabla_1^2 a_i(y,x)(\frac{\partial}{\partial x_k}(y), \frac{\partial}{\partial x_k}(y))\frac{\partial}{\partial x_i}(x)+o(1).
\end{array}\end{equation}
Inserting (\ref{46}) and (\ref{44}) (with $y$ approaching to $x$) into the left hand side of (\ref{43}), via Lemma \ref{17}, (\ref{449}) and (\ref{497}), we  get
 (\ref{43}).

Inserting (\ref{43}) into (\ref{42}), via Lemma \ref{17} and (\ref{312}), we can get (\ref{496}) with $X=\frac{\partial }{\partial x_j}$ and $V=\frac{\partial}{\partial x_k}$.
It is easy to check that $\nabla_2^3\nabla_1\rho^2(x,x)(X,F,V,W)$, $\nabla_2^2\nabla_1^2
\rho^2(x,x)(X,F,V,W)$ and
$R(X,V,F,W)(x)$ with $X,F,V,W\in T M$ are multi-linear with respect to $X,F,V$ and $W$. Hence, the desired identity (\ref{496}) follows.
 $\Box$

In particular, when $M$ is a two dimensional Riemannian manifold, we can employ Gauss-Bonnet Theorem to prove (\ref{312}).
\begin{lem}\label{414}
Assume that $M$ is a two dimensional Riemannian manifold. Let $X,V, F$ $\in T M$. Then, for any  $x\in M$, it holds that
\begin{equation}\label{418}
\nabla_{X(x)}\Big(d_x(L_{x\cdot}F(x))V(x)\Big)=\frac{1}{2} k(x)\Big(\langle F(x),V_X^\bot(x)\rangle X(x)-\langle F(x),X(x)\rangle V_X^\bot(x)\Big),
\end{equation} where $k(x)$ is the Gaussian curvature of  $M$ at the point $x$, and $V^\bot_X(x)$ is given by
$$
V_X^\bot(x)\equiv\cases{V(x)-\Big\langle V(x), \frac{X(x)}{|X(x)|}\Big\rangle  \frac{X(x)}{|X(x)|},\quad \textrm{if}\,\, X(x)\neq 0,\cr V(x),\quad \textrm{if}\,\,X(x)=0.}
$$
\end{lem}
\begin{rem}\label{416} Assume that $M$ is a two dimensional Riemannian manifold. If $X,Y\in T M$ are orthonormal at some point $x\in M$, then one has $k(x)=sec(X(x),Y(x))$. Applying this relation and formula (\ref{415}), one can easily check that (\ref{418}) is consistent with  (\ref{312}).

\end{rem}

\textbf{Proof of Lemma \ref{414}}\quad If $X(x)=0$ or $V(x)=0$, (\ref{418}) holds obviously.

If $X(x)\neq 0$, $V(x)\neq 0$ and $V^\bot_X(x)\neq 0$, set
$
e_1\equiv \frac{X(x)}{|X(x)|}$ and $ e_2=\frac{V^\bot_X(x)}{|V^\bot_X(x)|}.
$
Let us compute $\nabla_{e_1}(d_xL_{x\cdot}F(x)e_2)$ firstly. Set
$
\gamma_i(s)\equiv \exp_{x}(se_i)
$ with $s\geq0$ and  $i=1,2$. For any $y\in M$ closing  enough to $x$, by the definition of the  differential of a map, we have
$$
d_x(L_{xy}F(x))e_2=d_x(L_{xy}F(x))\dot{\gamma_2}(0)=\frac{d}{ds}\Big|_0L_{\gamma_2
(s)y}F(\gamma_2(s)).
$$
Using the definition of the covariant derivative in the sense of parallel translation (\ref{224}) and the above identity, we see that
\begin{equation}\label{319}\begin{array}{ll}
&\nabla_{e_1}(d_x(L_{x\cdot}F(x))e_2)
\\=&\lim\limits_{\tau\to 0^+}\frac{1}{\tau}\Big(L_{\gamma_1(\tau)x}d_x(L_{x\gamma_1(\tau)}F(x))e_2-d_x(L_{xx}F(x))e_2\Big)
\\ =& \lim\limits_{\tau,s\to 0^+}\frac{1}{\tau s}\Bigg(L_{\gamma_1(\tau)x}\Big(L_{\gamma_2(s)\gamma_1(\tau)}F(\gamma_2(s))-L_{x\gamma_1(\tau)}F(x)\Big)-\Big(L_{\gamma_2(s)x}F(
\gamma_2(s))
\\ &-F(x)\Big)\Bigg)
\\ =& \lim\limits_{\tau,s\to 0^+}\frac{1}{\tau s}\Big(L_{\gamma_1(\tau)x}L_{\gamma_2(s)\gamma_1(\tau)}F(\gamma_2(s))-L_{\gamma_2(s)x}F(\gamma_2(s))\Big)
\end{array}\end{equation}
It is well known that the parallel translation  of a vector along a curve conserves its norm and the angle between this vector and the curve. We will use this property to compute
$$L_{\gamma_1(\tau)x}L_{\gamma_2(s)\gamma_1(\tau)}F(\gamma_2(s))-L_{\gamma_2(s)x}F(\gamma_2(s)).$$
 Let $\gamma(\cdot)$ be the shortest geodesic connecting $\gamma_1(\tau)$ and $\gamma_2(s)$ with $\gamma(0)=\gamma_2(s)$ and $\gamma(1)=\gamma_1(\tau)$. Denote by $\Delta_{s\tau}$ the domain inside the geodesic triangle with vertexes  $\gamma_1(\tau)$, $\gamma_2(s)$ and $x$. The above formula represents the difference of parallel translations of $F(\gamma_2(s))$ along different curves.

For  $s\geq 0$, denote by   $\beta(s)$  the angle between $F(\gamma_2(s))$ and $\dot{\gamma}_2(s)$. Then, the angle between $L_{\gamma_2(s)x}F(\gamma_2(s))$ and $\dot{\gamma}_2(0)=e_2$ is $\beta(s)$, while the angle between  $L_{\gamma_2(s)x}F(\gamma_2(s))$ and $e_1$ is $\frac{\pi}{2}-\beta(s)$. Therefore,
\begin{equation}\label{316}
L_{\gamma_2(s)x}F(\gamma_2(s))=|F(\gamma_2(s))|\Big(\cos(\frac{\pi}{2}-\beta(s))e_1+\cos\beta(s)e_2\Big).
\end{equation}
  Denote the internal angles of the geodesic triangle $\Delta_{s\tau}$ at the points $\gamma_1(\tau)$ and $\gamma_2(s)$ by $\alpha_1$ and $\alpha_2$, respectively.  The angle between $F(\gamma_2(s))$ and $\dot{\gamma}(0)$ is $\pi-\alpha_2-\beta(s)$, so is the angle between $L_{\gamma_2(s)\gamma_1(\tau)}F(\gamma_2(s))$ and $\dot{\gamma}(1)$. Then, the angle between $L_{\gamma_2(s)\gamma_1(\tau)}F(\gamma_2(s))$ and $\dot{\gamma}_1(\tau)$ is       $\pi-\alpha_2-\alpha_1-\beta(s)$. Thus, the angle between $L_{\gamma_1(\tau)x}L_{\gamma_2(s)\gamma_1(\tau)}F(\gamma_2(s))$ and $\dot{\gamma}_1(0)=e_1$ is $\pi-\alpha_2-\alpha_1-\beta(s)$. Consequently, the angle between $L_{\gamma_1(\tau)x}L_{\gamma_2(s)\gamma_1(\tau)}F(\gamma_2(s))$ and $e_2$ is $\alpha_1+\alpha_2+\beta(s)-\frac{\pi}{2}$, due to the fact that $e_1$ is perpendicular to $e_2$. Therefore,
  \begin{equation}\label{314}\begin{array}{ll}
  &L_{\gamma_1(\tau)x}L_{\gamma_2(s)\gamma_1(\tau)}F(\gamma_2(s))
  \\ =&|F(\gamma_2(s))|\Big(\cos(\pi-\alpha_1-\alpha_2-\beta(s))e_1+\cos(\alpha_1+\alpha_2+\beta(s)-\frac{\pi}{2})e_2\Big).
 \end{array} \end{equation}
The Gauss-Bonnet formula (ref. \cite[p. 274]{c1}) gives
$$
\alpha_1+\alpha_2=\int\!\!\!\int\limits_{\Delta_{s\tau}}kdA+\frac{\pi}{2},
$$
where $dA$ is the Riemannian volume
element. By using (\ref{316}), (\ref{314}), the above identity and Taylor's theorem,  we have
\begin{equation}\label{320}\begin{array}{ll}
&L_{\gamma_1(\tau)x}L_{\gamma_2(s)\gamma_1(\tau)}F(\gamma_2(s))-L_{\gamma_2(s)x}F(\gamma_2(s))
\\ =&|F(\gamma_2(s))
|\int\!\!\!\int\limits_{\Delta_{s\tau}}kdA \Big(\cos\beta(s)e_1-\sin\beta(s)e_2\Big)+o(\int\!\!\!\int\limits_{\Delta_{s\tau}}kdA ).
\end{array}\end{equation}
We choose $\{s,\tau\}$ as the normal coordinates at the point $x$. When $s$ and $\tau$ are small enough,  the area of $\Delta_{s\tau}$, denoted by $A(\Delta_{s\tau})$, is given by
$$A(\Delta_{s\tau})=\int_0^s\langle \exp_x^{-1}\gamma(\frac{\xi}{s}),e_1\rangle d\xi. $$
Set $v(s,\tau)\equiv \exp_{\gamma_2(s)}^{-1}\gamma_1(\tau)$. Then, $\gamma(\zeta)=\exp_{\gamma_2(s)}(\zeta v(s,\tau))$ with $\zeta\in[0,1]$ and
$$\begin{array}{ll}
&\lim\limits_{s,\tau\to 0^+}\frac{1}{s\tau}A(\Delta_{s\tau})
\\ =&\lim\limits_{s\to 0^+}\frac{1}{s}\int_0^s\Big\langle d\exp_x^{-1}\Big|_{\exp_{\gamma_2(s)}\frac{\xi}{s}v(s,0)}\circ d\exp_{\gamma_2(s)}\Big|_{\frac{\xi}{s}v(s,0)}\frac{\xi}{s}(d\exp_{\gamma_2(s)}^{-1}\Big|_xe_1),e_1\Big\rangle d\xi
\\ \stackrel{\zeta=\xi /s}{=}&\lim\limits_{s\to 0^+}\int_0^1\Big\langle d\exp_x^{-1}\Big|_{\exp_{\gamma_2(s)}\zeta v(s,0)}\circ d\exp_{\gamma_2(s)}\Big|_{\zeta v(s,0)}\zeta(d\exp_{\gamma_2(s)}^{-1}\Big|_xe_1),e_1\Big\rangle d\zeta
\\ =& \frac{1}{2},
\end{array}$$
where we have used Lemma \ref{317}  and the  L'H\^opital's rule. Recalling (\ref{319}) and (\ref{320}), we have
$$\begin{array}{ll}
&\nabla_{e_1}(d_xL_{x\cdot}F(x)e_2)
\\ =& \lim\limits_{s,\tau\to 0^+}\frac{1}{2}\frac{1}{A(\Delta_{s\tau})}\Big(|F(\gamma_2(s))
|\int\!\!\!\int\limits_{\Delta_{s\tau}}kdA \Big(\cos\beta(s)e_1-\cos\Big(\frac{\pi}{2}-\beta(s)\Big)e_2\Big)+o(\int\!\!\!\int\limits_{\Delta_{s\tau}}kdA )\Big)
\\ =& \frac{1}{2}k(x)|F(x)|\Big(\cos\beta(0)e_1-\cos\Big(\frac{\pi}{2}-\beta(0)\Big)e_2\Big)
\\ =& \frac{1}{2}k(x)\Big(\Big\langle F(x), \frac{V^\bot_X}{|V^\bot_X|}\Big\rangle \frac{X}{|X|}-\Big\langle F(x), \frac{X}{|X|}\Big\rangle \frac{V^\bot_X}{|V^\bot_X|}\Big),
\end{array}$$
which implies  (\ref{418}).

If $X\neq 0$, $V\neq 0$ and $V^\bot_X=0$,  $V$ is parallel to $X$. Similar to  (\ref{319}), we can get $\nabla_{e_1}(d_xL_{x\cdot}F(x)e_1)=0$, which implies  $\nabla_X(d_xL_{x\cdot}F(x)V)=0.$ The proof of Lemma \ref{414} is complete. $\Box$

\begin{rem}\label{410}
If $X(x)=e_1$ and $V(x)=e_2$ are orthonormal, the left hand side of (\ref{312}) characterizes the limiting case of parallel translations of vector field $F$ along different curves. More precisely, by setting $\gamma_i(s)=\exp_x (s e_i)$ with $s\geq0$ and $i=1,2$, we can get (\ref{319}). 
\end{rem}

\setcounter{equation}{0}
\hskip\parindent \section{Statement of the main results}
\def\theequation{3.\arabic{equation}}

We begin with the following assumptions (Actually, as we shall see later, for our \textbf{Problem I}, we only need the assumptions $(C1)$ and $(C2)$ below):

\begin{description}

\item[$(C1)$] The maps $f(=f(t,x,u)): [0, T]\times M\times U \to T M$ and $f^0(=f^0(t,x,u)): [0, T]\times M\times U \to \R$  are measurable in $t$, continuous in $u$, and $C^1$ in $x$. Moreover,   there exists a constant $L>1$ such that,
\begin{equation}\begin{array}{c}\label{10}
\,|f^0(s,x_1,u)-f^0(s,x_2,u)|\leq L\rho(x_1,x_2),
\\[3mm] \;\;\;\;\;\;|L_{x_1x_2}f(s,x_1,u)-f(s,x_2,u)|\leq L\rho(x_1,x_2),
\\[3mm] |f^0(s,x_0,u)|\leq L,\qquad |f(s,x_0,u)|\leq L,\end{array}
\end{equation} for all $s\in [0,T]$, $u\in U$, and  $x_1, x_2\in M$ with $\rho(x_1,x_2)\leq \min\{i(x_1),i(x_2)\}$,  where $x_0\in M$ is arbitrarily fixed.

\item[$(C2)$]The maps $f(t,x,u)$ and $f^0(t,x,u)$ are $C^2$  in $x$. Furthermore,
 \begin{equation}\label{116}\begin{array}{l}
|\nabla_{x}f^0(t,x_1,u)-L_{x_2x_1}\nabla_xf^0(t,x_2,u)|\leq L\rho(x_1,x_2),
\\[3mm]|\nabla_{x}f(t,x_1,u)-L_{x_2x_1}\nabla_xf(t,x_2,u)|\leq L\rho(x_1,x_2),
\end{array}\end{equation} for all $x_1,x_2\in M$ with $\rho(x_1,x_2)<\min\{i(x_1),i(x_2)\}$ and $(t,u)\in[0,T]\times U$, where  $\nabla_xf(s,\cdot,u)$ and $\nabla_xf^0(s,\cdot,u)$ are the covariant derivatives of  $f(s,\cdot,u)$ and $f^0(s,\cdot,u)$  with respect to the state  variable, and their  norms are given by (\ref{270}).

\item[$(C3)$] $U\subset\R^m$ is  open. The maps $f(t,x,u)$ and $f^0(t,x,u)$ are $C^2$  in $u$. Furthermore,
\begin{equation}\label{437}\begin{array}{l}
|\nabla_uf^0(t,x,u_1)-\nabla_uf^0(t,x,u_2)|\leq L|u_1-u_2|,
\\[3mm]|\nabla_uf(t,x,u_1)-\nabla_uf(t,x,u_2)|\leq L|u_1-u_2|,
\end{array}\end{equation}
for all $u_1,u_2\in U$ and $(t,x)\in[0,T]\times M$, where $\nabla_uf(s,x,\cdot)$ and $\nabla_uf^0(s,x,\cdot)$ are the derivatives of $f(s,x,\cdot)$ and $f^0(s,x,\cdot)$ with respect to the control variable.

\end{description}

In this paper, for $X\in T M$ (or $X\in T^*M$), we denote by $\tilde X$ the dual covector (vector) of $X$.
Denote by
$H^\nu: [0,T]\times T^*M\times U\to \R$ the Hamiltonian function, defined by
\begin{equation}\label{65}H^\nu(t,x,p,u)\equiv p(f(t,x,u))+\nu f^0(t,x,u),\quad\forall\;(t,x,p,u,\nu)\in[0,T]\times T^*M \times U\times\R^-,
\end{equation}
with $R^- \equiv (-\infty,0]$. In particular, when $\nu=-1$, we denote by
\begin{equation}\label{420}H\equiv H^{-1}\end{equation}
 for abbreviation.

In this section, we fix a control $\bar{u}(\cdot)\in\U_{ad}$, where $\U_{ad}$ is given in (\ref{436}).  Let $\bar y(\cdot)$ be the solution to (\ref{25}) associated to $\bar u(\cdot)$.
For abbreviation, we denote
\begin{equation}\label{434}
[t]\equiv (t,\bar y(t),\bar u(t)),\quad\forall\; t\in[0,T].
\end{equation}
Suppose that $\psi(t)\in T_{\bar{y}(t)}^*M$ is the solution to the following first order dual equation${}^1$:
\footnote{${^1}$Our first order dual equation $(\ref{64})$ which contains the Levi-Civita connection $\nabla$, is the same as the dual equation (12.30) in \cite[p. 168]{as}, if they are rewritten in the local coordinates.  Nevertheless, the explicit form $(\ref{64})$ is convenient for us to derive further the second order necessary condition for optimal controls by introducing a suitable second order dual equation (i.e., (\ref{70}) below).}
 \begin{equation}\label{64}\displaystyle\cases{\nabla_{\dot{
 \bar{y}}(t)}\psi=-\nabla_xf[t]
 (\psi(t),\cdot)-\nu d_xf^0[t],\quad a.e. \,t\in[0,T),\cr
 \psi(T)=\psi_1}
 \end{equation}
  with $d_xf^0$ denoting the exterior derivative of $f^0$ with respect to the state   variable $x$,  and $\nabla_x f[t](\psi(t),\cdot)$ ($t\in[0,T]$) being a tensor given by
$$
\nabla_x f[t]\Big(\psi(t),X(\bar{y}(t))\Big)\equiv \nabla_{X(\bar{y}(t))}f(t,\cdot,\bar{u}(t))(\psi(t)),\qquad \forall\; X\in T M.
$$

\subsection{Optimality conditions for systems without endpoint constraints}

 Let $\psi(\cdot)$ be the solution to (\ref{64}) with $\nu=-1$ and $\psi_1=0$. Recalling (\ref{420}), we
put
$$
\widetilde{U}(t)\equiv\{u\in U;\ \ H(t,\bar{y}(t),\psi(t),\bar{u}(t))=H(t,\bar{y}(t),\psi(t),u)\},\quad t\in [0,T].$$

Let $w(t)\in \T_2^0(\bar y(t))$ with $t\in [0,T]$, and satisfy the following equation
\begin{equation}\label{70}\cases{\nabla_{\dot{\bar{y}}(t)}w+\E^{12}\Big(\nabla_x f[t]\otimes w(t)\Big)+\E^{12}\Big(w(t)\otimes\nabla_xf[t]\Big) \cr \noalign{\medskip}\quad\quad\quad\quad+\nabla_x^2 H(t,\bar{y}(t),\psi(t),\bar{u}(t))
-R\Big(\tilde{\psi}(t),\cdot,f[t],\cdot\Big)
=0,\quad t\in[0,T),
\cr \noalign{\medskip}w(T)=0,}
\end{equation} where $\E^{12}$ is defined in (\ref{76}), the
tensor $R\Big(\tilde\psi(t),\cdot,f[t],\cdot\Big)$  is given by
$$
R\Big(\tilde\psi(t),\cdot,f[t],\cdot\Big)(X,Y)\equiv R\Big(\tilde\psi(t),X,f[t],Y\Big),\quad\forall\; X,Y\in T M,
$$
and the tensor $\nabla_x^2 H(t,\bar{y}(t),\psi(t),\bar{u}(t))$  is given by
$$
\nabla_x^2H(t,\bar{y}(t),\psi(t),\bar{u}(t))(X,Y)\equiv \nabla_x^2f[t](\psi(t),X,Y)-\nabla_x^2f^0[t](X,Y),\quad\forall\; X,Y\in T M.
$$

Recalling the tensor contraction map in Section $2.4$, we introduce two tensors $\Phi(t)\in \T_0^2(\bar y(t))$ and $\Phi_1(t)\in \T_2^0(\bar y(t))$ with $t\in[0,T]$, which  solve the following equations respectively
\begin{equation}\label{459}\cases{
\nabla_{\dot{\bar y}(t)}\Phi=\E^{21}\Big(\nabla_xf[t]\otimes\Phi(t)\Big),
\quad t\in(0,T],
\cr \Phi(0)(X,Y)=\langle X,Y\rangle,\,\,\forall\; X,Y\in T_{y_0}^*M,}\end{equation}
and
\begin{equation}\label{458}\cases{\nabla_{\dot{\bar y}(t)}\Phi_1=-\E^{12}\Big(\Phi_1(t)\otimes \nabla_xf[t] \Big),
\quad t\in(0,T],
\cr \Phi_1(0)(X,Y)=\langle X,Y\rangle,\,\,\forall\; X,Y\in T_{y_0}M.
}\end{equation}

Our second order necessary  conditions for optimal controls without endpoint constraints are stated as follows.

\begin{thm}\label{66}
Assume that the assumptions $(C1)$ and $(C2)$  hold, and $(U,\tilde d)$ is a separable metric space.  Let $(\bar y(\cdot), \bar u(\cdot))$ be an optimal pair for \textbf{Problem I}. Then,
\begin{equation}\label{3}
H(t,\bar y(t),\psi(t),\bar u(t))=\max\limits_{u\in U}H(t,\bar y(t),\psi(t),u),\quad \hbox{a.e. }t\in[0,T],
\end{equation}
where $\psi(\cdot)$ is the solution to the first order dual equation (\ref{64}) with $\nu=-1$ and $\psi_1=0$. Furthermore,
 for any $u(\cdot)\in\U_{ad}$ with $u(t)\in\widetilde U(t)$, a.e. $t\in[0,T]$, it holds that
\begin{equation}\label{457}
\begin{array}{l}
\int_0^T\int_0^t\Big\{\Big\langle\nabla_xH(t,\bar y(t),\psi(t),u(t))
 -\nabla_xH(t,\bar y(t),\psi(t),
 \bar u(t)),
 \\[2mm] \quad {\mathcal E}^{21}\Big({\mathcal E}^{21}\Big(\Phi(t)\otimes L_{\bar y(s)\bar y(t)}^{\bar y(\cdot)}\Phi_1(s)\Big)\otimes
 L_{\bar y(s)\bar y(t)}^{\bar y(\cdot)}(f(s,\bar y(s),u(s))-f[s]) \Big)\Big\rangle
 \\[2mm] \quad+\frac{1}{2}\Big(w(t)+w(t)^\top\Big)\Big(f(t,\bar y(t),u(t))-f[t],
  {\mathcal E}^{21}\Big({\mathcal E}^{21}\Big(\Phi(t)\otimes L_{\bar y(s)\bar y(t)}^{\bar y(\cdot)}\Phi_1(s)\Big)\otimes
\\[2mm] \quad L_{\bar y(s)\bar y(t)}^{\bar y(\cdot)}(f(s,\bar y(s),u(s))-f[s]) \Big)\Big)\Big\}dsdt\\[2mm]
  \leq 0,
\end{array}
\end{equation}
where the $2$-form $w^\top(t)$ is the transpose  of $w(t)$, defined by $
w^\top(t)(X,V)\equiv w(t)(V,X)$ for any $X,V\in T_{\bar{y}(t)}M$, $\E^{21}$ is defined in (\ref{76}), and the tensor $\nabla_xH(t,\bar{y}(t),\psi(t),u)$ (with $u\in U$) is given by
\begin{equation}\label{469}
\nabla_x H(t,\bar{y}(t),\psi(t),u)(X)\equiv \nabla_x f(t,\bar{y}(t),u)(\psi(t),X)-d_xf^0(t,\bar{y}(t),u)(X),\quad\forall\; X\in T M.
\end{equation}
Moreover, if $U$ is a Polish space, then the following pointwise condition holds:
\begin{equation}\label{71}\begin{array}{ll}
\displaystyle\frac{1}{2}\Big(w(t)+w^\top(t)\Big)\Big(f[t]-f(t,\bar{y}(t),v),f[t]-f(t,\bar{y}(t),v)\Big)
\\[3mm] \displaystyle\quad+\Big(\nabla_xH(t,\bar{y}(t),\psi(t),\bar{u}(t))-\nabla_xH(t,\bar{y}(t),\psi(t),v)\Big)\Big(f[t]-f(t,\bar{y}(t),v)\Big)\\[3mm]
\displaystyle\leq 0, \qquad \forall\; v\in \widetilde{U}(t), \quad \hbox{a.e. }t\in[0,T].
\end{array}\end{equation}

\end{thm}

\begin{rem} The explicit dependence of the second order dual variable $w(\cdot)$ (involved  in  equation (\ref{70})) on the curvature tensor $R(\tilde{\psi},\cdot,f,\cdot)$ reveals that the second order necessary condition (\ref{71}) depends on the  curvature.
When $M$ is $\R^n$, the  curvature is zero everywhere, and for this special case,  Theorem \ref{66} coincides with \cite[Theorem 4.3]{l}.
\end{rem}

Denote by $d(\cdot,\cdot)$ the Ekeland metric over the admissible control space $\U_{ad}$, given by
$$
d(u_1(\cdot),u_2(\cdot))=|\{t\in[0,T];u_1(t)\neq u_2(t)\}|,\;\forall\; u_1(\cdot),\,u_2(\cdot)\in\U_{ad}.
$$
It is well known that $(\U_{ad},d)$ is a complete metric space.

Our second order sufficient condition for locally optimal controls for \textbf{Problem I} is given as follows:
\begin{thm}\label{456} Assume that $(C1)$ and $(C2)$  hold, and $(U,\tilde d)$ is a separable metric space. If an admissible pair $(\bar y(\cdot),\bar{u}(\cdot))$ for \textbf{Problem I} satisfies (\ref{3}), and there exist an $\epsilon_0>0$ and a $\beta>0$ satisfying
\begin{equation}\label{20}
\begin{array}{l}
\int_0^T\int_0^t\Big\{\Big\langle\nabla_xH(t,\bar y(t),\psi(t),u(t))
 -\nabla_xH(t,\bar y(t),\psi(t),
 \bar u(t)),
 \\[2mm] \quad {\mathcal E}^{21}\Big({\mathcal E}^{21}\Big(\Phi(t)\otimes L_{\bar y(s)\bar y(t)}^{\bar y(\cdot)}\Phi_1(s)\Big)\otimes
 L_{\bar y(s)\bar y(t)}^{\bar y(\cdot)}(f(s,\bar y(s),u(s))-f[s]) \Big)\Big\rangle
 \\[2mm] \quad +\frac{1}{2}\Big(w(t)+w(t)^\top\Big)\Big(f(t,\bar y(t),u(t))-f[t],
  {\mathcal E}^{21}\Big({\mathcal E}^{21}\Big(\Phi(t)\otimes L_{\bar y(s)\bar y(t)}^{\bar y(\cdot)}\Phi_1(s)\Big)\otimes
 \\[2mm] \quad L_{\bar y(s)\bar y(t)}^{\bar y(\cdot)}(f(s,\bar y(s),u(s))-f[s]) \Big)\Big)\Big\}dsdt\\[2mm]
  \leq -\beta d(u(\cdot),\bar u(\cdot))^2,\qquad\forall\; u(\cdot)\in{\mathcal U}_{ad}\hbox{ with }d(u(\cdot),\bar u(\cdot))<\epsilon_0,
\end{array}
\end{equation}
 then, one can find an $\epsilon_1>0$ such that
$$
J(v(\cdot))\geq J(\bar u(\cdot)),\qquad\forall\;v(\cdot)\in \V\equiv \{u(\cdot)\in{\U}_{ad};\ \ d(u(\cdot),\bar u(\cdot))\leq \epsilon_1\}.
$$
\end{thm}

\subsection{Optimality  conditions for systems with endpoint constraints}
Before stating the main results of this subsection, we introduce some more notations. Recalling (\ref{65}), under the assumptions $(C2)$ and $(C3)$, we introduce the following two tensors:
\begin{eqnarray}\label{468}
&&\nabla_x^2 H^\nu(t,x,p,u)(X,Y)=
 \nabla_x^2f(t,x,u)(p,X,Y)+\nu\nabla_x^2f^0(t,x,u)(X,Y),
\\ \label{524}
&&\nabla_u\nabla_xH^\nu(t,x,p,u)(X,\eta)=\frac{\partial}{\partial u}\Big(\nabla_x f(t,x,u)(p,X)
+\nu d_xf^0(t,x,u)(X)\Big)
\cdot \eta,
\end{eqnarray}
for all $(t,x,p,u,\nu)\in[0,T]\times T^*M \times U\times\R^-$ and $X,Y\in T M$ and $\eta\in\R^m$.

Assume that the condition $(C1)$  hold and $(U,\tilde d)$ is a separable metric space. If $(\bar y(\cdot), \bar u(\cdot))$ is a optimal pair for  \textbf{Problem II}, then there exist a $\nu\leq0$ and a $\psi_1\in T^*_{\bar y(T)}M$ with $|\nu|+|\psi_1|>0$, such that the following Pontryagin's type maximum principle holds (\cite[p. 181, Theorem 12.10]{as}):
\begin{equation}\label{505}
H^\nu(t,\bar y(t),\bar\psi(t),\bar u(t))=\max_{u\in U}H^\nu(t,\bar y(t),\bar\psi(t),u), \,a.e. t\in[0,T],
\end{equation}
where   the covector field $\bar\psi(\cdot)$ along $\bar y(\cdot)$ solves (\ref{64}).
Moreover, when $U\subset \R^m$ is open, and $f$ and $f^0$ are  $C^1$ in $u$, it is easy to see that (\ref{505}) implies that
\begin{equation}\label{547}
\frac{\partial}{\partial u}\Big|_{\bar u(t)}H^\nu(t, \bar y(t),\bar\psi(t),u)=0,\quad a.e.\, t\in[0,T].
\end{equation}
For any $\xi(\cdot)\in L^2(0,T;\R^m)$, we denote by $V(t)\in T_{\bar y(t)}M$ ($t\in[0,T]$) the solution to the following first order variational equation:
\begin{equation}\label{522}
\cases{\nabla_{\dot{\bar y}(t)}V=\nabla_{V(t)}f(t,\cdot,\bar u(t))+\nabla_u f[t](\cdot,\xi(t)),\quad a.e.\,t\in(0,T],\cr V(0)= V(T)=0.}
\end{equation}

Based on the above first order necessary condition, we give below a second order necessary condition for optimal pairs of \textbf{Problem II}.
\begin{thm}\label{504}
  Assume that the conditions  $(C1)$, $(C2)$ and $(C3)$ hold, and $(\bar y(\cdot), \bar u(\cdot))$ is an optimal pair for \textbf{Problem II}. Assume that $U\subset \R^m$ is open, $\bar u (\cdot)\in L^2(0,T;U)$, and for any neighborhood $\O$ of $\bar u(\cdot)$ in $L^2(0,T;\R^m)$, there exists a $v(\cdot)\in\V_{ad}\cap \O$. If $(\nu,\psi_1)$ is a unique pair (up to a positive factor) satisfying (\ref{547}), then the following second order  necessary condition holds:
\begin{description}
\item[(i)] If $\nu<0$ (the normal case), then \begin{equation}\label{507}\begin{array}{l}
\int_0^T\Big\{\nabla_u^2H^\nu(t,\bar y(t),\bar\psi(t),\bar u(t))(\xi(t),\xi(t))+\nabla_x^2H^\nu(t,\bar y(t),\bar\psi(t),\bar u(t))(V(t),V(t))
\\[2mm] \quad+2\nabla_u\nabla_xH^\nu(t,\bar y(t),\bar\psi(t),\bar u(t))(V(t),\xi(t))-R(\tilde{\bar\psi}(t),V(t),\dot{\bar y}(t),V(t))\Big\}dt\\[2mm] \leq 0,
\end{array}\end{equation}
for all $\xi(\cdot)\in L^2(0,T;\R^m)$ and $V(\cdot)\in T_{\bar y(\cdot)}M$ satisfying (\ref{522}), where the tensors $\nabla_x^2H^\nu$ and
$\nabla_u\nabla_xH^\nu$ are given respectively by (\ref{468}) and (\ref{524}).

\item[(ii)] If $\nu=0$ (the abnormal case), the left hand side of (\ref{507}) is sign-definite on the following set:
\begin{equation}\label{544}\begin{array}{rr}
\Big\{\xi(\cdot)\in L^2(0,T;\R^m);&\,\int_0^T\Big\{\langle\nabla_xf^0[t],V(t)\rangle+\nabla_uf^0[t](\xi(t))\Big\}dt=0,
\\&\textrm{where $V(\cdot)\in T_{\bar y(\cdot)}M$ solves (\ref{522})}\Big\}.
\end{array}\end{equation}
\end{description}
 \end{thm}
\begin{rem} In optimal control theory, usually one may use the needle variation technique to deduce a pointwise necessary condition from the integral one, provided that the later still holds  when needle variations are
used as control perturbations (e.g. the proof of \cite[Theorem 1.6, pp. 150--158]{yz} for the first order case, and the proof of \cite[Theorem 4.3, p. 1455]{l} for the second order case).  However, in our case, (\ref{507}) does not enjoy such a property, due to the hard restriction (\ref{522}). Hence, it seems quite difficult to derive a pointwise condition from  (\ref{507}).

\end{rem}

Our second order sufficient conditions for locally optimal controls for \textbf{Problem II} is stated as follows:

\begin{thm}\label{453}  Assume that the conditions  $(C1)$, $(C2)$ and $(C3)$ hold, $U\subset \R^m$ is open, and $(\bar y(\cdot), \bar u(\cdot))$ is an admissible pair for \textbf{Problem II} satisfying (\ref{505}) for some $\nu<0$ and $\psi_1\in T_{y_1}^*M$. If there exist a $\gamma_0>0$ such that
$$
\frac{\partial^2}{\partial u^2}\Big|_{\bar u(t)}H^\nu(t,\bar y(t),\bar \psi(t),u)(v,v)\leq-\gamma_0|v|^2,\quad \forall\; v\in\R^m, \,a.e.\, t\in[0,T],
$$
and a $\gamma_1>0$ so that
 \begin{equation}\label{452}\begin{array}{l}
\int_0^T\Big\{\nabla_u^2H^\nu(t,\bar y(t),\bar\psi(t),\bar u(t))(\xi(t),\xi(t))+\nabla_x^2H^\nu(t,\bar y(t),\bar\psi(t),\bar u(t))(V(t),V(t))
\\[2mm]  \quad+2\nabla_u\nabla_xH^\nu(t,\bar y(t),\bar\psi(t),\bar u(t))(V(t),\xi(t))-R(\tilde{\bar\psi}(t),V(t),\dot{\bar y}(t),V(t))\Big\}dt
\\[2mm]  \leq -\gamma_1\|\xi\|_{L^2(0,T;\R^m)}^2,
\end{array}\end{equation}
for all $\xi\in L^2(0,T;\R^m)$ and $V(\cdot)\in T_{\bar y(\cdot)}M$ satisfying (\ref{522}), then  $\bar u(\cdot)$ is locally optimal in the $L^2(0,T;U)$ topology.
\end{thm}

\setcounter{equation}{0}
\section{Two illustrative examples}
\def\theequation{4.\arabic{equation}}

In this section, we shall give two simple but illustrative examples to show the usefulness of our results in the last section.

The first two lines of (\ref{10}) and (\ref{116}) are essentially Lipschitz  conditions, and they are the key assumptions  for our main results. The following result provides a method which can be employed to  verify them.

\begin{lem}\label{425} Let $\T$ be a tensor on $ M$. Then, the following two conditions are equivalent:
\begin{description}
\item[$(i)$] There exists a positive constant $L$ such that
$
|\nabla \T|\leq L$;
\item[$(ii)$] There exists a positive constant $L$ such that
$
|L_{x_1x_2} \T(x_1)- \T(x_2)|\leq L\rho(x_1,x_2)$,
 for all $x_1, x_2\in  M$ with $\rho(x_1,x_2)<\min\{i(x_1),i(x_2)\}$.
\end{description}
\end{lem}

\textbf{Proof:}  \, Suppose $\T\in T  M$. Assume that ($i$) holds.  Take any $x_1,x_2\in  M$ with $0<\rho(x_1,x_2)<\min\{i(x_1), i(x_2)\}$. Set $V=\frac{1}{\rho(x_1,x_2)}\exp_{x_2}^{-1}x_1$. Then, by Lemma \ref{17}, we have $|V|=1$. Assume that $\{e_1, e_2, \cdots, e_n\}$ (with $e_1\equiv V$) is an orthonormal basis at $T_{x_2}M$. Set $\gamma(s)\equiv \exp_{x_2}sV$ for $s\in[0,\tau]$ with $\tau\equiv \rho(x_1,x_2)$. Then, $\gamma$ is the shortest geodesic connecting $x_1$ and $x_2$, and  parameterized by the  arc length. Let $\{d_1,\cdots,d_n\}$ be the dual basis of $\{e_1,\cdots,e_n\}$ at $x_2$. Set
$$
e_i(s)=L_{x_2\gamma(s)}e_i,\quad d_i(s)=L_{x_2\gamma(s)} d_i,\qquad s\in[0,\tau],\,\,i=1,\cdots,n.
$$
Then, $\{e_1(s),\cdots,e_n(s)\}$ is the orthonormal basis at $T_{\gamma(s)}M$, $\{d_1(s),\cdots,d_n(s)\}$ is the orthonormal basis at $T_{\gamma(s)}^*M$ for $s\in [0,\tau]$, and they are dual to each other. Thus, we can express $\T$ along the curve $\gamma$ as
$$\T(\gamma(s))=\sum\limits_{i=1}^na_{i}(s)e_i(s),\quad \forall\;s\in[0,\tau].$$
Take any $v=\sum\limits_{k=1}^n v_k d_k\in T_{x_2}^*M$. By  ($i$), (\ref{100}), (\ref{82}), the parallel translation of $e_i(\cdot)$ and $d_i(\cdot)$ along $\gamma$ for $i=1,\cdots,n$ and the Mean Value Theorem, one can get
$$\begin{array}{ll}
&\ds\Big|\T(x_2)(v)-L_{x_1x_2}\T(x_1)(v)\Big|
 = \ds\Big|\sum\limits_{i=1}^n a_{i}(0)e_i\Big(\sum\limits_{k=1}^nv_kd_k\Big)-\T(x_1)\Big(L_{x_2x_1}v\Big)\Big|
\\[5mm] =& \ds\Big| \sum\limits_{i=1}^n a_{i}(0)v_i -\sum\limits_{i=1}^na_{i}(\tau)e_i(\tau)\Big(\sum\limits_{k=1}^n v_k d_k(\tau)\Big)\Big|
 = \Big|\sum\limits_{i=1}^n(a_{i}(0)-a_{i}(\tau))v_i \Big|
\\[5mm] = & \ds\Big| \tau \sum\limits_{i=1}^n \dot{a}_{i}(\zeta)v_i \Big|
 = \Big| \rho(x_1,x_2)\sum\limits_{i=1}^n \dot{a}_{i}(\zeta) e_i(\zeta)\Big(\sum\limits_{k=1}^nv_kd_k(\zeta)\Big)\Big|
\\[5mm] = & \ds \Big|\rho(x_1,x_2)\nabla_{\dot{\gamma}(\zeta)}\T\Big(L_{x_2\gamma(\zeta)}v\Big) \Big|
 = \Big|\rho(x_1,x_2)\nabla \T(\gamma(\zeta))\Big(L_{x_2\gamma(\zeta)}v, \dot{\gamma}(\zeta)\Big)\Big|
\\[2mm]  \leq &\ds L |v| \rho(x_1,x_2),
\end{array}$$
 where $\zeta\in[0,1]$, and we get ($ii$).
 Conversely, if ($ii$) holds, we can apply (\ref{224}) to get ($i$).

 When $\T$ is a tensor of general order, one can apply (\ref{82}) and the same method as above to get the equivalence of ($i$) and ($ii$).
 $\Box$

 \medskip

By means of the  above lemma, we give below an  example which satisfies the assumptions $(C1)$ and $(C2)$.

\begin{exl}\label{exl} Consider the hyperbolic surface under the hyperboloid model
 $$H_R^2\equiv\{(x_1,x_2,x_3)\in\R^3; \; x_3^2-(x_1^2+x_2^2)=R^2, x_3>0\},$$ endowed  with the metric $h=\iota^*m$, where $\iota: H_R^2\to \R^3$ is the  inclusion,  $m$ is the Minkowski metric
$
m=dx_1^2+dx_2^2-dx_3^2,
$ and $h$ is the pullback metric.  $H_R^2$ is non-compact, and its  Gaussian curvature is $-1$. For the details of this model, we refer to  \cite[Proposition 3.5, p. 38]{le}. Let $(H_R^2, \Psi)$ be the local coordinates with
$$
\Psi(x_1,x_2)=\Big(x_1,x_2,\sqrt{R^2+x_1^2+x_2^2}\Big),\quad (x_1,x_2)\in \R^2.
$$
 Let $U=\{1,2,3,4\}$ with the discrete metric.
 Denote $$f(x,u)=u^3e^{-(R^2+x_1^2+x_2^2)}(x_2\frac{\partial}{\partial x_1}-x_1\frac{\partial}{\partial x_2}),
 \qquad f^0(x,u)=u^2e^{-(R^2+x_1^2+x_2^2)},
 $$ for $x=(x_1,x_2)\in \R^2$ and $u\in U$.  We consider the optimal control problem $(\ref{25})$ and $(\ref{26})$ on $H_R^2$, with $f$ and $f^0$ given above. Then, the conditions $(C1)$ and $(C2)$ are satisfied.
\end{exl}

\textbf{Proof:} \, It is easy to check that  the third line  of (\ref{10}) hold. We will apply Lemma \ref{425}  to check the first two   inequalities of (\ref{10}) and (\ref{116}).
By a direct computation, the metric $h$ can be expressed  in the local coordinates $(x_1,x_2)$ by
$$\begin{array}{ll}
h=&\ds\Big(1-\frac{x_1^2}{R^2+x_1^2+x_2^2}\Big)dx_1\otimes dx_1-\frac{x_1x_2}{R^2+x_1^2+x_2^2}(dx_1\otimes dx_2+dx_2\otimes dx_1)+
\\[2mm] &\ds\Big(1-\frac{x_2^2}{R^2+x_1^2+x_2^2}\Big)dx_2\otimes dx_2.
\end{array}$$
Hence,  $\nabla f$ and $\nabla^2 f$ in the coordinates $x=(x_1,x_2)\in R^2$ are given respectively as follows:
\begin{eqnarray}
\ds\nabla_x f(x,u)&=&\ds u^3e^{-(R^2+x_1^2+x_2^2)}\Big(  (-2x_1x_2-\frac{x_1x_2}{R^2})\frac{\partial}{\partial x_1}\otimes dx_1+(-2x_2^2+1+\frac{x_1^2}{R^2})\frac{\partial}{\partial x_1}\otimes dx_2\nonumber
\\[2mm] \label{411}  &&\ds+(2x_1^2-\frac{x_2^2}{R^2}-1)\frac{\partial }{\partial x_2}\otimes dx_1+(2x_1x_2+\frac{x_1x_2}{R^2})\frac{\partial}{\partial x_2}\otimes dx_2\Big),
\\[2mm] \label{412}\ds\nabla_x^2 f(x,u)&=&\ds u^3\sum_{i,j,l=1}^2a_{ijl}(x_1,x_2)\frac{\partial}{\partial x_i}\otimes dx_j\otimes dx_l,
\end{eqnarray}
where $a_{ijl}(x_1,x_2)=e^{-(R^2+x_1^2+x_2^2)}(R^2+x_1^2+x_2^2)^{-1}p_{ijl}(x_1,x_2)$ with $p_{ijl}(x_1,x_2)$ being a polynomial of $x_1$ and $ x_2$ for $i,j,l=1,2$.
We  also obtain that, for any tangent  vector field $Y=a_1\frac{\partial}{\partial x_1}+a_2\frac{\partial}{\partial x_2}$ and any cotangent vector  field $\eta=b_1dx_1+b_2dx_2$,
\begin{equation}\label{413}\begin{array}{c}
\ds|Y|\leq 1 \ \ \,\,\,\textrm{if and only if}\,\,\,\ \  a_1^2+a_2^2\leq \frac{R^2+x_1^2+x_2^2}{R^2},
\\[2mm]\ds
|\eta|\leq 1 \ \ \,\,\,\textrm{if and only if}\,\,\,\ \  b_1^2+b_2^2\leq \frac{R^2+x_1^2+x_2^2}{R^2},
\end{array}\end{equation} where we have used the definition (\ref{270}).
Furthermore, for any $x=(x_1,x_2)$, by using (\ref{411}), (\ref{412}) and (\ref{413}),  we get
$$\begin{array}{lll}
\ds|\nabla_x f(x,u)|&\equiv & \ds\sup\{\nabla_x f(x,u)\Big(\eta,Y\Big); \eta\in T_x^*M,Y \in T_xM, |Y|\leq 1, |\eta|\leq 1\}
\\[2mm] &= & \ds u^3\sup \Big\{ e^{-(R^2+x_1^2+x_2^2)}(a_1\,\,a_2)\left(\begin{array}{cc}-2x_1x_2-\frac{x_1x_2}{R^2} \quad & -2x_2^2+1+\frac{x_1^2}{R^2} \\[2mm]\ds  2x_1^2-\frac{x_2^2}{R^2}-1&\ds 2x_1x_2+\frac{x_1x_2}{R^2}
\end{array} \right)\left(\begin{array}{c}b_1  \\[2mm]\ds b_2 \end{array}\right);
\\[3mm] &&\ds \max\{a_1^2+a_2^2,\;\; b_1^2+b_2^2\}\leq \frac{R^2+x_1^2+x_2^2}{R^2} \Big\}
\\[2mm] & \leq & L
\end{array}
$$
and
$$\begin{array}{lll}
\ds  |\nabla_x^2 f(x,u)| & \equiv & \ds\sup \Big\{  \nabla_x^2f(x,u)(\eta, Y, Z);  \eta\in T_x^*M, Y,Z\in T_xM,  |\eta|\leq 1, |Y|\leq 1, |Z|\leq 1\Big\}
\\[2mm]  & = &\ds u^3\sup\Big\{ \sum\limits_{i,j,l=1}^2 a_{ijl} \frac{\partial}{\partial x_i}\otimes dx_j\otimes dx_l\sum\limits_{k,m,n=1}^2( a_kdx_k, b_m\frac{\partial}{\partial x_m}, c_n\frac{\partial}{\partial x_n});
\\[3mm] &&\ds\max\{a_1^2+a_2^2, b_1^2+b_2^2, c_1^2+c_2^2\}\leq \frac{R^2+x_1^2+x_2^2}{R^2}\Big\}
\\[2mm] & =&\ds u^3\sup\Big\{ \sum\limits_{i,j,l=1}^2 a_{ijl} a_ib_j c_l; \max\{a_1^2+a_2^2, b_1^2+b_2^2, c_1^2+c_2^2\}\leq \frac{R^2+x_1^2+x_2^2}{R^2}\Big\}
\\[2mm] &\leq&\ds L
\end{array}$$ for some positive constant $L$. By applying Lemma \ref{425}, we obtain     the second inequalities of (\ref{10}) and (\ref{116}).  Using the same way, one can get the first inequalities of (\ref{10}) and (\ref{116}).
$\Box$

By means of Pontryagin-type maximum principle, Sussmann (\cite[Example 5.10.2]{s2}) gave a new proof of the classical result that the locally shortest curve connecting two fixed points on a Riemannian manifold must be a geodesic. In the following example, we shall apply our second order necessary condition for optimal controls with endpoint constraints (Theorem \ref{504}) to recover the well-known second variation of energy:

\begin{exl}\label{510} Given any two points $y_0,y_1\in M$, assume that $\bar y: [0,T]\to M$ is a smooth curve with $\bar y(0)=y_0$, $\bar y(T)=y_1$ and constant speed $|\dot{\bar y}(t)|\equiv C>0$ (for all $t\in[0,T]$). If $\bar y(\cdot)$ is a locally shortest curve connecting $y_0$ and $y_1$, then it must be a geodesic. Moreover, if $\bar y(\cdot)$ is a geodesic, the following second variation of energy holds (see \cite[p. 159]{p1}):
\begin{equation}\label{451}
\int_0^T\Big\{|\nabla_{\dot{\bar y}(t)}V|^2+R(\dot{\bar y}(t),V(\bar y(t)),\dot{\bar y}(t),V(\bar y(t)))\Big\}dt\geq 0
\end{equation}
for any $V\in T M$ with $V(y_0)=0$ and $V(y_1)=0$.

 \end{exl}

\textbf{Proof}\,\,According to the completeness of $(M,g)$, given any bounded domain $D\subset M$ such that $\bar y(t)\in D$ for all $t\in[0,T]$, one can find $f_1,\cdots, f_m\in T M$ ($m\in \N $) with compact supports, such that
\begin{equation}\label{516}
span\{f_1|_{\overline{D}},\cdots, f_m|_{\overline{D}}\}=\{X|_{\overline{D}}; \ \; X\in T M\}.
\end{equation}
Indeed, by the completeness of  $(M,g)$ and the Hopf-Rinow Theorem (see \cite[Theorem $16$, p. 137]{p1}), we see that  $\overline{D}$ is compact. Hence, by  \cite[Lemma $1$, p. 52]{wsy} or \cite[Theorem $3.7$, p. 72]{c}, there exist a $\delta>0$ and $x_1,\cdots,x_l\in \overline{D}$ ($l\in\N$) such that, for all $x\in \overline{D}$, the map $\exp_x: B(O,\delta)\subset T_x M\to B_x(\delta)$ is diffeomorphic and $\cup_{i=1}^l B_{x_i}(\delta/2)\supset \overline{D}$. Therefore, for each $i=1,\cdots, l$, we can define vector fields $\{f_i^j(\cdot)\}_{j=1}^n$ on $B_{x_i}(\delta)$ as a basis for $T M$ restricted to $B_{x_i}(\delta)$ as follows:
$$
f_i^j(y)=d \exp_{x_i}\Big|_{\exp_{x_i}^{-1}y}e_i^j,\quad \forall\; y\in B_{x_i}(\delta),
$$
where $\{e_i^j\}_{j=1}^n$ is a basis at $x_i$. Then, we extend $\{f_i^j(\cdot)\}_{j=1}^n$ from $B_{x_i}(\delta)$ to $M$ smoothly, which are denoted by $\{\tilde{f}_i^j\}_{j=1}^n$, such that $f_i^j\Big|_{B_{x_i}(\delta/2)}=\tilde{f}_i^j\Big|_{B_{x_i}(\delta/2)}$ and $\tilde{f}_i^j\Big|_{M\setminus B_{x_i}(\delta)}=0$, for $j=1,\cdots,n$. Thus, $\{\tilde{f}_i^j\}$ ($i=1,\cdots,l, j=1,\cdots,n$) can express linearly any vector field restricted to $\overline{D}$.

According to (\ref{516}), we see that
\begin{equation}\label{511}
\cases{\dot y(t)=\sum_{i=1}^mu_i(t)f_i(y(t)),\,a.e.t\in[0,T],\cr
y(0)=y_0,\quad y(T)=y_1,}
\end{equation}
with $u(t)=(u_1(t),\cdots,u_m(t))\in\R^m$  a.e. $t\in[0,T]$, describes all absolutely continuous curves  contained in $D$, with endpoints $y_0$ and $y_1$. In order to apply our optimal control results to the above shortest curve problem, we need to seek an appropriate cost functional.   To this end, set
$$
\Omega=\{y: [0,T]\to M; \ \; y(0)=y_0,\,y(T)=y_1,\,\textrm{and }\, y(\cdot)\,\textrm{ is absolutely continuous}\}.
$$
Denote by $L(y)$ the length of $y\in\Omega$. \cite[Proposition 17, p. 126]{p1} says that, if $\bar y(\cdot)$ is a minimizer of  $L$ in $\Omega$, then it also minimizes $E(y)\equiv \frac{1}{2}\int_0^T|\dot y(t)|^2dt$ within $\Omega$. Therefore, we define the cost functional as follows:
\begin{equation}\label{512}
J(u(\cdot))=\frac{1}{2}\int_0^T|\dot y_u(t)|^2dt=\frac{1}{2}\int_0^T|\sum_{i=1}^mu_i(t)f_i(y_u(t
))
|^2dt,
\end{equation}
where $y_u(\cdot)$ is the solution to (\ref{511}) corresponding
 to the control $u(\cdot)\in \U\equiv\{u:[0,T]\to \R^m ; \,u(\cdot)\,\textrm{is measurable, and }y_u(0)=y_0, y_u(T)=y_1\} $.

Hence, we choose (\ref{512}) to be the desired cost functional.
The corresponding Hamiltonion function is
\begin{equation}\label{421}
H^\nu(t,y,p,u)\equiv p(\sum_{i=1}^m u_i f_i(y))+\frac{1}{2}\nu|\sum_{i=1}^m u_if_i(y)|^2,\end{equation}
for $(t,y,p,u,\nu)\in [0,T]\times T^*M\times \R^m\times \R^-$.

Suppose that $(\bar y(\cdot), \bar u(\cdot))$ is optimal for problem (\ref{511}) and (\ref{400}) with $J$ given in (\ref{512}). Then $\bar y(\cdot)$ is the locally shortest curve connecting $y_0$ and $y_1$. Applying Theorem \ref{504} to this specific problem, we obtain that
\begin{equation}\begin{array}{l}
\bar\psi(t)\Big(\sum_{i=1}^m\bar u_i(t)f_i(\bar y(t))\Big)+\frac{\nu}{2}|\sum_{i=1}^m\bar u_i(t)f_i(\bar y(t))|^2
\\
\label{518}=\max\limits_{
u=(u_1,\cdots,u_m)\in \R^m}\{\bar\psi(t)\Big(\sum_{i=1}^m
u_if_i(\bar y(t))\Big)+\frac{\nu}{2}|\sum_{i=1}^m u_if_i(\bar y(t))|^2\}, \,a.e.\, t\in[0,T],
\end{array}\end{equation}
where  the first order dual variable $\bar\psi$ solves
\begin{equation}\label{514}\cases{
\nabla_{\dot{\bar y}(t)}\bar\psi=-\sum_{i=1}^m \bar u_i(t)\nabla f_i(\bar y(t))(\bar\psi(t),\cdot)-\nu\sum_{i,j=1}^m\bar u_i(t)\bar u_j(t)\nabla f_i(\bar y(t))(\tilde f_j(\bar y(t)),\cdot),\cr
\quad\quad\quad\quad\quad\quad\quad\quad\quad\quad\quad\quad\quad\quad\quad\quad\quad\quad\quad\quad\quad\quad\quad\quad\quad\quad
\quad\quad\quad\quad a.e.\, t\in[0,T],\cr
\bar\psi(T)=\psi_1,}\end{equation}
with
\begin{equation}\label{513}
|\nu|+|\psi_1|>0,\quad \nu\leq 0,
\end{equation}
and $\tilde f_j$ ($j=1,\cdots,m$) being the dual covector of $f_j$. Clearly, (\ref{518}) is the first order necessary condition for $(\bar y(\cdot), \bar u(\cdot))$.

By differentiating (\ref{518}) with respect to the control variable at $\bar u(t)$ for almost every $t\in[0,T]$, we have
\begin{equation}\label{515}
\bar\psi(t)(f_i(\bar y(t)))+\nu\langle f_i(\bar y(t)),\dot{\bar y}(t)\rangle=0, \quad a.e.\, t\in[0,T]
\end{equation}
 for $i=1,\cdots,m$.
According to (\ref{516}), for any $V\in T M$ such that $V(y_0)=0$ and $V(y_1)=0$, we can find suitable real valued functions on $[0,T]$: $\xi_1(\cdot),\cdots, \xi_m(\cdot)$, such that
\begin{equation}\label{550}
\nabla_{\dot{\bar y}(t)}V-\sum_{i=1}^m\bar u_i(t)\nabla_{V(\bar y(t))}f_i=\sum_{i=1}^m \xi_i(t)f_i(\bar y(t)),\quad t\in[0,T].
\end{equation}
Combining  (\ref{515}) and the above identity, we deduce that
$$
\int_0^T\Big\{\bar\psi(t)\Big(\nabla_{\dot{\bar y}(t)}V-\sum_{i=1}^m\bar u_i(t)\nabla_{V(\bar y(t))}f_i\Big)+\nu\langle\nabla_{\dot{\bar y}(t)}V-\sum_{i=1}^m\bar u_i(t)\nabla_{V(\bar y(t))}f_i,\dot{\bar y}(t)\rangle\Big\}dt=0.
$$
Using integration by parts to the above identity, via the dual equation (\ref{514}), we can get
\begin{equation}\label{517}
-\nu\int_0^T\langle V(\bar y(t)),\nabla_{\dot{\bar y}(t)}\dot{\bar y}\rangle dt=0.
\end{equation}

If $\nu=0$, (\ref{515}) implies $\bar\psi(t)=0$ for all $t\in[0,T]$, which contradicts (\ref{513}). Hence, $\nu<0$. By the choice of $V$, (\ref{517}) implies $\nabla_{\dot{\bar y}(t)}\dot{\bar y}=0$ for all $t\in[0,T]$, which means that $\bar y(\cdot)$ is  a geodesic.

From (\ref{515}) we can get
\begin{equation}\label{525}
\dot{\bar y}(t)=-\frac{1}{\nu}\tilde{\bar\psi}(t),\quad\forall\; t\in[0,T],
\end{equation} which implies $(\nu,\psi_1)$
 is unique up to a positive factor. In order to deduce the second order necessary condition, we need to compute $\nabla_u^2H^\nu$, $\nabla_u\nabla_xH^\nu$ and $\nabla_x^2H^\nu$ at $(t,\bar y(t),\bar\psi(t),\bar u(t))$ respectively.

Recalling (\ref{421}), for any $1\leq i, j\leq n$, we have
\begin{equation}\label{419}\frac{\partial^2}{\partial u_i\partial u_j}\Big|_{\bar u(t)}H^\nu(t,\bar y(t),\bar\psi(t),u)=\nu\langle f_i(\bar y(t), f_j(\bar y(t)))\rangle,\quad \forall\; t\in[0,T].\end{equation}
Hence, for any $\xi\in L^2(0,T;\R^m)$ satisfying (\ref{522}) (or (\ref{550})), we get
 \begin{equation}\label{523}
 \nabla_u^2H^\nu(t,\bar y(t),\bar \psi(t),\bar u(t))(\xi(t),\xi(t))=\nu|\sum_{i=1}^m\xi_i(t)f_i
 (\bar y(t))|^2.
 \end{equation}
Recalling (\ref{468}), (\ref{524}) and (\ref{525}), we obtain that, via the first order variational equation (\ref{550}),
$$\begin{array}{l}
\nabla_x^2H^\nu(t,\bar y(t),\bar \psi(t),\bar u(t))(V(t),V(t))
\\=\nu\Big(|\nabla_{\dot{\bar y}(t)}V|^2+|\sum_{i=1}^m\xi_i(t)f_i(\bar y(t))|^2
-2\langle\nabla_{\dot{\bar y}(t)}V,\sum_{i=1}^m\xi_i(t)f_i(\bar y(t))\rangle\Big);
\\ \nabla_u\nabla_xH^\nu(t,\bar y(t),\bar \psi(t),\bar u(t))(V(t),\xi(t))=\nu\langle\nabla_{V(t)}\Big(\sum_{j=1}^m\bar u_j(t)f_j\Big),\sum_{i=1}^m\xi_i(t)f_i(\bar y(t))\rangle.
\end{array}$$
Inserting (\ref{523}) and the above two identities into (\ref{507}), and via (\ref{550}), (\ref{525}) and (\ref{522}), we end up with
$$\begin{array}{lll}
0&\leq&\int_0^T\Big\{|\nabla_{\dot{\bar y}(t)}V|^2+R(\dot{\bar y}(t),V(t),\dot{\bar y}(t),V(t))+2\langle \nabla_{V(t)}\Big(\sum_{j=1}^m\bar u_j(t)f_j\Big)
\\ &&-\nabla_{\dot{\bar y}(t)}V,\sum_{i=1}^m\xi_i(t)f_i(\bar y(t))\rangle +2|\sum_{i=1}^m\xi_i(t)f_i(\bar y(t))|^2\Big\}dt
\\  & = &\int_0^T\Big\{|\nabla_{\dot{\bar y}(t)}V|^2+R(\dot{\bar y}(t),V(t),\dot{\bar y}(t),V(t))\Big\}dt,
\end{array}$$
which implies (\ref{451}).
 $\Box$

\begin{rem}\label{rem4.1} The Legendre condition (\cite[Theorem 20.6, p. 300 and Proposition 20.11, p. 310]{as})
 says that the second order necessary condition of an optimal control $\bar u(\cdot)$ is :
 $$\frac{\partial^2}{\partial u^2}\Big|_{\bar u(t)}H^\nu(t,\bar y(t),\bar\psi(t),u)(v,v)\leq 0,\quad \forall\; v\in\R^m,\,a.e.t\in[0,T].$$
 Applying this condition to deal with the concrete problem in Example \ref{510}, instead of (\ref{451}), one can only  get
 \begin{equation}\label{5zhqiu23}
\nu|\sum_{i=1}^m v_if_i(\bar y(t))|^2\leq 0,\,\,\forall\; v=(v_1,\cdots,v_m)^\top\in\R^m,
 \end{equation}
 a.e. $t\in[0,T]$, where (\ref{419}) is used. Clearly, the condition (\ref{5zhqiu23}) is trivially correct. Hence, our second order necessary condition in Theorem \ref{504} provides more information than that in \cite{as}. In addition, from the inequality (\ref{5zhqiu23}), it follows that the optimal control is not totally singular, i.e. $\bar u(\cdot)$ does not fulfills $\frac{\partial^2}{\partial u^2}\Big|_{\bar u(t)}H^\nu(t,\bar y(t),\bar\psi(t),u)=0$ $a.e.\; t\in[0,T]$, and therefore, the Goh condition and the generalized Legendre condition do not work.
\end{rem}

\setcounter{equation}{0}
\hskip\parindent \section{Variations of Trajectories}
\def\theequation{5.\arabic{equation}}
In this section,  we shall give the first and second order variations of a trajectory of the control system (\ref{25}), by  employing respectively the needle variation in the case that the control set $U$ is a  metric space,  and the classical variation in the case that $U$ is an open subset in $\R^m$.

\subsection{Needle Variation }

In this subsection, we assume that $U$ is a metric space.
\begin{pro}\label{31}
Assume that the conditions $(C1)$ and $(C2)$ hold. Fix any $\bar u(\cdot)\in\U_{ad}$, and let $\bar y(\cdot)$ be the corresponding solution to  (\ref{25}). For any $u(\cdot)\in\U_{ad}$, denote by $X_u(t),Y_u(t)\in T_{\bar y(t)}M$ respectively  solutions to the following
 first and second order variational equations:
\begin{equation}\label{227}\cases{
\nabla_{\dot{\bar y}(t)}X_u=\nabla_{X_u(t)}f(t,\cdot,\bar{u}(t))+
f(t,\bar y (t),u(t))-f[t],\quad a.e.\,t\in(0,T],\cr X_u(0)=0,
}
\end{equation}
and
\begin{equation}\label{477}\cases{
\langle\nabla_{\dot{\bar y}(t)}Y_u,Z\rangle=\langle \nabla_{Y_u(t)}f(t,\cdot,\bar{u}(t)),Z\rangle-\frac{1}{2}R(Z,X_u(t),f[t],X_u(t))
\cr\,\,\,\quad\quad\quad\quad\quad\quad+\langle \nabla_{X_u(t)}f(t,\cdot,u(t))-\nabla_{X_u(t)}f(t,\cdot,\bar u(t)),Z\rangle
 \cr \quad\quad\quad\quad\quad\quad\,\,+\frac{1}{2}\nabla_x^2f[t](\tilde Z,X_u(t),X_u(t)),\quad a.e.\,t\in(0,T],\cr
Y_u(0)=0,
}\end{equation}
for any $Z\in T M$, where $[t]$ is given in (\ref{434}).
Then,  as $d(u(\cdot),\bar u(\cdot))\equiv\epsilon\to 0$,
\begin{eqnarray}\label{32}
|X_u(t)|=O(\epsilon),&|Y_u(t)|=O(\epsilon^2),
\\ \label{37}
|V_u(t)-X_u(t)|=o(\epsilon),& \;\;|V_u(t)-X_u(t)-Y_u(t)|=o(\epsilon^2),
\end{eqnarray}
where
\begin{equation}\label{33}
V_u(t)=\exp_{\bar y(t)}^{-1}y_u(t),\quad t\in[0,T],
 \end{equation}
and  $y_u(\cdot)$ is the corresponding solution to  (\ref{25}) with the control $u(\cdot)$.
\end{pro}

To prove Proposition \ref{31}, we need the following result.
\begin{lem}\label{480}
Under the assumptions in Proposition \ref{31}, for any $u(\cdot)\in\U_{ad}$, it holds that \begin{eqnarray}\label{6}&\rho(y_u(\hat{t}),y_u(t))\leq (1+\rho(x_0,y_0))(e^{L\hat{t}}-e^{Lt}),\qquad\forall\; 0\leq t\leq \hat{t}\leq T,\end{eqnarray}
 where $x_0\in M$ and $L>0$ are given in assumption $(C1)$.
 Furthermore, there exists a  constant  $\epsilon_1>0$ such that
\begin{equation}\label{478}
\rho(\bar y(t),y_u(t))\leq 2L(1+\rho(x_0,y_0))e^{Lt}\Big|[0,t]\cap \{\tau\in[0,T];u(\tau)\neq\bar u(\tau)\}\Big|,
\end{equation}
for any $ t\in[0,T]$ and $u(\cdot)\in\U_{ad}$ with $d(\bar u(\cdot),u(\cdot))<\epsilon_1$.
\end{lem}
\textbf{Proof:} \, The proof is divided into two parts.

\textbf{Part I}. In this part, we prove (\ref{6}). Firstly, we claim that
\begin{equation}\label{240}
|f(t,y,v)|\leq L\Big(\rho(x_0,y)+1\Big)
\end{equation} holds for any $t\in[0,T]$, $y\in M$ and $v\in U$.

In fact, by  the completeness of the manifold $M$, we can find the shortest geodesic $\gamma$ connecting $x_0$ and $y$ with $\gamma(0)=x_0$ and $\gamma(1)=y$. Let $0=s_0<s_1<s_2<\cdots<s_{N-1}<s_N=1$ (for some $N\in \N$) be such that $\rho(\gamma(s_j),\gamma(s_{j-1})
)<\min\{ i(\gamma(s_j)),i(\gamma(s_{j-1}))\}$ for $j=1,2,\cdots,N$. By means of  (\ref{268}), the condition $(C1)$, and the triangle inequality of the norm $|\cdot|$ yields
$$\begin{array}{ll}
&|f(t,y,v)|
\\ \leq& |L_{\gamma(s_{N-1})\gamma(s_N)}f(t,\gamma(s_{N-1}),v)-f(t,y,v)|+|f(t,\gamma(s_{N-1}),v)|\\
\leq &|L_{\gamma(s_{N-1})\gamma(s_N)}f(t,\gamma(s_{N-1}),v)-f(t,y,v)|
\\ &+|f(t,\gamma(s_{N-1}),v)-L_{\gamma(s_{N-2})\gamma(s_{N-1})}f(t,\gamma(s_{N-2}),v)|+|f(t,\gamma(s_{N-2}),v)|
 \\ \leq &\cdots
 \\ \leq & \sum\limits_{j=1}^N |L_{\gamma(s_{j-1})\gamma(s_j)}f(t,\gamma(s_{j-1}),v)-f(t,\gamma(s_j),v)|+|f(t,x_0,v)|
 \\ \leq & L \sum\limits_{j=1}^N  \rho\Big(\gamma(s_{j-1}),\gamma(s_j)\Big)+L
 =L\Big(\rho(x_0,y)+1\Big),
\end{array}$$ which implies (\ref{240}).

Secondly, we estimate $|f(t, y_u(t),v)|$ for any  $v\in U$ and $t\in[0,T]$. Let $\beta: [0,1]\to
M$ be the shortest geodesic connecting $y_u(t)$ and $y_0$ and satisfy $\beta(0)=y_0$ and $\beta(1)=y_u(t)$. By the compactness, we can choose $\tau_0=0,\tau_1,\tau_2,\cdots,\tau_l=1$ for some $l\in\N$, such that  $\tau_0<\tau_1<\cdots<\tau_l$ and $\rho(\beta(\tau_i),\beta(\tau_{i+1}))<\min\{i(\beta(\tau_i)),i(\beta(\tau_{i+1}))$ for $i=0,\cdots,l-1$. Applying Lemma \ref{17}, the condition $(C1)$, (\ref{240}), and the triangle inequality of the norm $|\cdot|$, we obtain that
\begin{equation}\label{27}\begin{array}{ll}
&|f(t,y_u(t),v)|
\\ \leq& |f(t,y_u(t),v)-L_{\beta(\tau_{l-1})\beta(\tau_l)}f(t,\beta(\tau_{l-1}),v)|
\\ &+|f(t,\beta(\tau_{l-1}),v)
-L_{\beta(\tau_{l-2})\beta(\tau_{l-1})}f(t,\beta(\tau_{l-2}),v)|+\cdots+|f(t,y_0,v)|
\\ \leq & L\sum\limits_{i=1}^l\rho(\beta(\tau_i),\beta(\tau_{i-1})) +L(\rho(x_0,y_0)+1)
\\ = & L\rho(y_u(t),y_0)+L\rho(x_0,y_0)+L.
\end{array}
\end{equation}

Thirdly, we estimate $\rho(y_u(t),y_0)$ for any $t\in[0,T]$. Applying (\ref{27}), we get
$$\begin{array}{ll}
\rho(y_u(t),y_0)\leq \int_0^t|\dot{y}_u(s)|ds
\leq  L\int_0^t\rho(y_u(s),y_0)ds+Lt(\rho(y_0,x_0)+1).\end{array}
$$
Applying the Gronwall inequality, we obtain
\begin{equation}\label{28}
\rho(y_u(t),y_0)\leq (1+\rho(x_0,y_0))(e^{Lt}-1),\,\forall\;t\in[0,T].
\end{equation}

Finally, for any $t,\hat{t}\in[0,T]$ with $t\leq \hat{t}$, applying (\ref{27}) and (\ref{28}), we have
$$\begin{array}{ll}
\rho(y_u(t),y_u(\hat{t}))&\leq\int_t^{\hat{t}}|\dot{y}_u(s)|ds
\\ &\leq L\int_t^{\hat{t}}\Big(\rho(y_u(s),y_0)+\rho(x_0,y_0)+1\Big)ds
\\ & \leq(1+\rho(x_0,y_0))(e^{L\hat{t}}-e^{Lt}),
\end{array}$$ which yields (\ref{6}).

\textbf{Part II}. In this part, we prove (\ref{478}). We will firstly show that (\ref{478}) holds on a small interval, and then we extend it to the whole  interval $[0,T]$. From (\ref{6}), it follows that, for any $u(\cdot)\in\U_{ad}$, the corresponding trajectory $y_u(\cdot)$ is contained in the closed ball $\overline{B_{y_0}( R)}$ with $R\equiv (1+\rho(x_0,y_0))(e^{LT}-1)$, where $B_{y_0}( R)$ is defined in (\ref{409}).Clearly $\overline{B_{y_0}( R)}$  is compact, because of  the completeness of the manifold $M$ and the Hopf-Rinow Theorem (see   \cite[Theorem $16$, p. 137]{p1}). Therefore, there exists a $\delta>0$ such that, for any $x\in \overline{B_{y_0}( R )}$, the map $\exp_x:B(O,\delta)\subset T_xM\to B_x(\delta)$ is diffeomorphic, due to   \cite[Lemma $1$, p. 52]{wsy} or a direct consequence of  \cite[Theorem $3.7$, p. 72]{c}.

We claim that (\ref{478}) holds on the interval $[0,t_1]$, where
$$
t_1=\frac{1}{L}\ln\frac{\delta+2(1+\rho(x_0,y_0))}{2(1+\rho(x_0,y_0))}.
$$
In fact, by the triangle inequality of the norm $|\cdot|$  and (\ref{6}), for any $t\in[0,t_1)$, we have
$$\begin{array}{ll}
\rho(y_u(t),\bar y(t))&\leq\rho(y_u(t),y_0)+\rho(y_0,\bar y(t))
\\ & \leq 2(e^{Lt}-1)(1+\rho(x_0,y_0))
\\&<\delta.
\end{array}$$
Denote by $I_{u\neq \bar u}$ the indicator function of the set $\{\tau\in[0,T];\ \;u(\tau)\neq \bar u(\tau)\}$. For $t\in[0,t_1)$, using the condition $(C1)$, Lemma \ref{17}, (\ref{476}), (\ref{205}), (\ref{268}), (\ref{6}), (\ref{27}), and (\ref{28}), we get
$$
\begin{array}{ll}
&\frac{d}{dt}\rho^2(y_u(t), \bar y(t))
\\ =& \langle  \nabla_1\rho^2(y_u(t), \bar y(t)),\dot{y}_u(t)\rangle +\langle  \nabla_2\rho^2(y_u(t), \bar y(t)),\dot{\bar y}(t)\rangle
\\ = & -2\langle  \exp_{y_u(t)}^{-1}\bar y(t),f(t,y_u(t),u(t))\rangle -2\langle  \exp_{\bar y(t)}^{-1}y_u(t),f[t]\rangle
\\=&-2\langle \exp_{y_u(t)}^{-1}\bar y(t),f(t,y_u(t),u(t))\rangle-2\langle\exp_{y_u(t)}^{-1}\bar y(t),f(t,
y_u(t),\bar u(t))-f(t,
y_u(t),\bar u(t))\rangle
\\  &-2\langle\exp_{\bar y(t)}^{-1}y_u(t),f[t]\rangle
\\  =& 2\langle \exp_{\bar y(t)}^{-1}y_u(t),L_{y_u(t)\bar y(t)}f(t,y_u(t),\bar u(t))-f[t]\rangle
\\ &-2\langle \exp_{y_u(t)}^{-1}\bar y(t),f(t,y_u(t),u(t))-f(t,y_u(t),\bar u(t))\rangle I_{u\neq \bar u}(t)
\\ \leq & 2L \rho(\bar y(t),y_u(t))\Big(\rho(\bar y(t),y_u(t))+2e^{Lt}(1+\rho(x_0,y_0))I_{u\neq \bar u}(t)\Big),
\end{array}
$$
which leads to
\begin{equation}\label{426}
\frac{d}{dt}\rho(y_u(t),y(t))\leq L\Big(\rho(\bar y(t),y_u(t))+2e^{Lt}(1+\rho(x_0,y_0))I_{u\neq \bar u}(t)\Big).
\end{equation}
Applying Gronwall's inequality to the above inequality, we get     (\ref{478}) on $[0,t_1)$. By the inequality (\ref{6}) and the continuity of $\rho(\cdot,\cdot)$,  (\ref{478}) actually holds on the closed interval  $[0,t_1]$.

Secondly, we show that (\ref{478}) holds on the interval $[0,t_2]$ with
$$
t_2=\frac{1}{L}\ln\Big\{\frac{\delta}{2(1+\rho(x_0,y_0))}
+e^{Lt_1}(1-L\epsilon_1)\Big\},
$$
 where
\begin{equation}\label{51}
\epsilon_1\equiv \min\Big\{T,\frac{\delta}{3(1+\rho(x_0,y_0))Le^{LT}}\Big\}.
\end{equation}
We already showed that (\ref{478}) holds on $[0,t_1]$. It remains to show that, for $u(\cdot)\in\U_{ad}$ with $d(u(\cdot),\bar u(\cdot))<\epsilon_1$,  (\ref{478}) still holds on $[t_1,t_2]$. For this purpose, for any $s\in[t_1,t_2)$, applying (\ref{478}) with $t=t_1$ and (\ref{6}), we see that
$$\begin{array}{ll}
\rho(y_u(s), \bar y(s))&\leq \rho(y_u(s),y_u(t_1))+\rho(y_u(t_1), \bar y(t_1))+\rho( \bar y(t_1), \bar y(s))
\\ & \leq 2(e^{Ls}-e^{Lt_1})(1+\rho(x_0,y_0))+2L(1+\rho(x_0,y_0))
e^{Lt_1}d(u(\cdot),\bar u(\cdot))
\\ & <\delta.
\end{array}
$$
 Then, similarly to (\ref{426}),  by Gronwall's inequality and (\ref{478}) with $t=t_1$, we can get (\ref{478}) on $[t_1,t_2)$. Furthermore,
by the continuity of $\rho^2(\cdot,\cdot)$ and (\ref{6}), one can easily get (\ref{478}) on $[0,t_2]$.

Using the same technique as above, one can prove by induction that (\ref{478}) holds on $[0,t_i]$ for $i=2,3,\dots$ with
\begin{equation}\label{34}
t_i=\frac{1}{L}\ln\Big\{\frac{\delta}{2(1+\rho(x_0,y_0))}+
e^{Lt_{i-1}}(1-L\epsilon_1)\Big\}.
\end{equation}

We claim that there exists an integer $I>0$ being large enough such that  $[0,T]\subset [0,t_I]$. Let us use the contradiction argument and assume that
\begin{equation}\label{441}
t_i<T,\quad \forall\;i\in\{0\} \cup\N.
\end{equation}
Then,  there would exist a $\tilde{t}>0$ and  a subsequence of $\{t_i\}_{i=0}^{+\infty}$ (still denoted  by $\{t_i\}_{i=0}^{+\infty}$),  such that  $\lim_{i\to +\infty}t_i=\tilde{t}$. Letting $i\to+\infty$ in (\ref{34}) yields
$\tilde{t}=\frac{1}{L}\ln \frac{\delta}{2(1+\rho(x_0,y_0))L\epsilon_1}$. Recalling (\ref{51}), one has $\tilde{t}>T$, which contradicts (\ref{441}). Hence we get (\ref{478}) on $[0,T]$. $\Box$

\medskip

We are now in a position to prove Proposition \ref{31}.

\textbf{Proof of Proposition \ref{31}}\,
First, we will prove (\ref{32}). Multiplying both sides of (\ref{227}) by $X_u(t)$, using $(C1)$ and Lemma \ref{425}, we get
$$\begin{array}{ll}
&\frac{1}{2}|X_u(t)|\frac{\partial}{\partial t}|X_u(t)|
\\ =&\nabla_x f[t](X_u(t),X_u(t))+\langle f(t,\bar y(t),u(t))-f[t],X_u(t)\rangle
\\ \leq & L |X_u(t)|^2+| f(t,\bar y(t),u(t))-f[t]||X_u(t)|,\end{array}$$
which leads to
$$
|X_u(t)|\leq 2L \int_0^t|X_u(s)|ds+2\int_0^t|f(s,\bar y(s),u(s))-f[s]|ds.
$$
Applying Gronwall's inequality to the above inequality, we can get, via $(C1)$, the first estimate in (\ref{32}).

To prove the second estimate in (\ref{32}), we take $Z=Y_u$ in (\ref{477}), and get
$$\begin{array}{lll}
\frac{1}{2}|Y_u(t)|\frac{\partial}{\partial t}|Y_u(t)|
&=&\nabla_x f[t](Y_u(t),Y_u(t))+\nabla_x(f(t,\bar y(t),u(t))-f[t])(Y_u(t),X_u(t))
\\ & &-\frac{1}{2}R(Y_u(t),X_u(t),f[t],X_u(t))
+\frac{1}{2}\nabla_x^2f[t](\tilde Y_u(t),X_u(t),X_u(t)).
\end{array}$$
Applying $(C1)$, $(C2)$,  Lemma \ref{425} and the first relation  of (\ref{32}), we can find a positive constant  $C>0$ such that
$$
\frac{\partial}{\partial t}|Y_u(t)|\leq C(|Y_u(t)|+|\nabla_xf(t,\bar y(t),u(t))-\nabla_xf[t]||X_u(t)|+O(\epsilon^2)).
$$
Integrating the above inequality on $[0,t]$, by Gronwall's inequality, and applying  the first estimate in (\ref{32}), the condition ($C2$) and Lemma \ref{425}, we  obtain the second estimate in (\ref{32}).

Next, we are going to prove (\ref{37}).

By Lemma \ref{480}, there exists an $\epsilon_0>0$ such that, for all $u(\cdot)\in\U_{ad}$ with $d(u(\cdot),\bar u(\cdot))=\epsilon<\epsilon_0$, we have $\rho(y_u(t),\bar y(t))<i(\bar y(t))$ for all $t\in[0,T]$. Then, we can define $V_u(\cdot)$ as in (\ref{33}).
For $t\in[0,T]$, let  $\tilde V_u(t)$ be a vector at $T_{\bar y(t)}M$ satisfying
\begin{equation}\label{481}
\tilde V_u(t)\equiv\cases{ \frac{1}{\rho(\bar y(t),y_u(t))}V_u(t),\quad\textrm{if}\;|V_u(t)|\neq 0;
\cr 0,\quad\quad\quad\quad\quad\quad\quad\;\textrm{if}\,|V_u(t)|= 0.}
\end{equation}
Then, we define a geodesic connecting $\bar y(t)$ and $y_u(t)$ as follows:
\begin{equation}\label{482}
\beta(\theta;t)\equiv \exp_{\bar y(t)}\theta \tilde V_u(t),\quad \theta\in[0,\rho(\bar y(t),y_u(t))].
\end{equation}
The inequality (\ref{478}) indicates that $\beta(\cdot;t)$ is the shortest geodesic connecting $\bar y(t)$ and $y_u(t)$, provided $d(\bar u(\cdot),u(\cdot))=\epsilon$ is  small  enough.
 In particular, in the case that $|V_u(t)|\neq0$, $\beta(\cdot;t)$ is parameterized by the arc length. For any $Z\in T M$,
applying Lemma \ref{17}, Lemma \ref{480} and Taylor's expansion,  we get
$$\begin{array}{ll}
&\langle \nabla_{\dot{\bar y}(t)}V_u,Z\rangle
\\=&-\frac{1}{2}\nabla_2\nabla_1\rho^2\Big(\bar y(t),y_u(t)\Big)(Z,f(t,y_u(t),u(t)))
+\frac{1}{2}\nabla_2\nabla_1\rho^2\Big(\bar y(t),\bar y(t)\Big)(Z,f(t,\bar y(t),u(t)))
\\ &-\frac{1}{2}\nabla_1^2\rho^2\Big(\bar y(t),y_u(t)\Big)(Z,f[t])
+\frac{1}{2}\nabla_1^2\rho^2(\bar y(t),\bar y(t))(Z,f[t])
+\langle Z,f(t,\bar y(t),u(t))-f[t]\rangle
\\ =&-\frac{1}{2}\frac{\partial}{\partial \theta}|_{\theta=0}\nabla_2\nabla_1\rho^2(\bar y(t),\beta(\theta;t))\Big(Z,f(t,\beta(\theta;t),u(t))\Big)
\rho(\bar y(t),y_u(t))
\\ & -\frac{1}{4}\frac{\partial^2}{\partial \theta^2}|_{\theta=0}\nabla_2\nabla_1\rho^2(\bar y(t),\beta(\theta;t))\Big(Z,f(t,\beta(\theta;t),u(t))\Big)
\rho^2(\bar y(t),y_u(t))
\\ &-\frac{1}{2}\frac{\partial}{\partial \theta}|_{\theta=0}\nabla_1^2\rho^2(\bar y(t),\beta(\theta;t))\Big(Z,f[t]\Big)\rho(\bar y(t),y_u(t))+\langle Z,f(t,\bar y(t),u(t))-f[t]\rangle
\\ &-\frac{1}{4}\frac{\partial^2}{\partial\theta^2}|_{\theta=0}
\nabla_1^2\rho^2(\bar y(t),\beta(\theta;t))\Big(Z,f[t]\Big)\rho^2(\bar y(t),y_u(t))
+o(\epsilon^2)
\\ =&\langle Z,\nabla_{V_u(t)}f(t,\cdot,\bar u(t))\rangle +\langle Z,f(t,\bar y(t),u(t))-f[t]\rangle
+\langle Z,\nabla_{V_u(t)}f(t,\cdot,u(t))
\\ &-\nabla_{V_u(t)}f(t,\cdot,
\bar u(t))\rangle
-\frac{1}{4}\nabla_2^3\nabla_1\rho^2(\bar y(t),\bar y(t))\Big(Z,f(t,\bar y(t), u(t)),V_u(t),V_u(t)\Big)
\\ &-\frac{1}{4}\nabla_2^2\nabla_1^2\rho^2(\bar y(t),\bar y(t))\Big(Z,f[t],V_u(t),V_u(t)\Big)
\\ &+\frac{1}{2}\langle Z,\nabla_{V_u(t)}\nabla_{\frac{\partial}{\partial \theta}\beta(\theta;t)}f(t,\cdot,u(t))|_{\theta=0}\rangle\rho(\bar y(t),y_u(t))+o(\epsilon^2).\end{array}
$$
By Lemma \ref{311}, we obtain
\begin{equation}\label{35}\begin{array}{ll}
&\langle \nabla_{\dot{\bar y}(t)}V_u,Z\rangle
\\ =&\langle Z,\nabla_{V_u(t)}f(t,\cdot,\bar u(t))\rangle +\langle Z,f(t,\bar y(t),u(t))-f[t]\rangle
+\langle Z,\nabla_{V_u(t)}f(t,\cdot,u(t))
\\ &-\nabla_{V_u(t)}f(t,\cdot,
\bar u(t))\rangle
-\frac{1}{2}R(f[t],V_u(t),Z,V_u(t))
\\ &+\frac{1}{2}\langle Z,\nabla_{V_u(t)}\nabla_{\frac{\partial}{\partial \theta}\beta(\theta;t)}f(t,\cdot, u(t))|_{\theta=0}\rangle\rho(\bar y(t),y_u(t))
\\ & -\frac{1}{4}\nabla_2^3\nabla_1\rho^2(\bar y(t),\bar y(t))\Big(Z, f(t,\bar y(t),u(t))-f[t],V_u(t),V_u(t)\Big)+o(\epsilon^2).
\end{array}\end{equation}
We observe that, for any $v\in U$,
\begin{equation}\label{38}\begin{array}{l}
\langle Z,\nabla_{V_u(t)}\nabla_{\frac{\partial}{\partial \theta}\beta(\theta;t)}f(t,\cdot,v)|_{\theta=0}\rangle\rho(\bar y(t),y_u(t))
=\nabla_x^2f(t,\bar y(t),v)(\tilde Z,V_u(t),V_u(t)).
\end{array}\end{equation}
In fact, recalling the definition of $\tilde{V}_u(t)$, (\ref{99}), (\ref{100}) and (\ref{482}), we have
$$\begin{array}{ll}
&\textrm{The left hand side of (\ref{38})}
\\ = & V_u(t)\langle \nabla_{\frac{\partial}{\partial \theta}\beta(\theta;t)}f(t,\cdot,v),Z\rangle \rho(y_u(t),\bar y(t))-\langle\nabla_{\tilde{V}_u(t)}f(t,\cdot,v),\nabla_{V_u(t)}Z\rangle \rho(y_u(t),\bar y(t))
\\ = &V_u(t)\Big(\nabla_x f(t,\cdot,v)(Z,\frac{\partial}{\partial \theta}\beta(\theta;t))\Big)\rho(y_u(t),\bar y(t))-\nabla_x f(t,\bar y(t),v)(\nabla_{V_u(t)}Z,V_u(t))
\\ =&\Big(\nabla_{V_u(t)}(\nabla_x f(t,\cdot,v))\Big)(Z,\frac{\partial}{\partial\theta}\beta(\theta;t)|_{\theta=0})\rho(y_u(t),\bar y(t))
\\ =& \nabla_x^2f(t,\bar y(t),v)(Z,V_u(t),V_u(t)),
\end{array}$$
which implies (\ref{38}). Inserting (\ref{38}) into (\ref{35}),
and subtracting (\ref{227}) from (\ref{35}), we get that
\begin{equation}\label{39}\cases{
\langle\nabla_{\dot{\bar y}(t)}(V_u-X_u),Z\rangle=\langle Z ,\nabla_{V_u(t)-X_u(t)}f(t,\cdot,\bar u(t))\rangle-\frac{1}{2}R(f[t],V_u(t),Z,V_u(t))
\cr \quad\quad\quad\quad+\nabla_x\Big( f(t,
\bar y(t), u(t))-f[t]\Big)(\tilde Z,V_u(t))+\frac{1}{2}\nabla_x^2f[t](\tilde Z,V_u(t),V_u(t))
\cr \quad\quad\quad\quad
-\frac{1}{4}\nabla_2^3\nabla_1\rho^2(\bar y(t),\bar y(t))\Big(Z, f(t,\bar y(t),u(t))-f[t],V_u(t),V_u(t)\Big)
\cr\quad\quad\quad\quad+\frac{1}{2}\nabla_x^2(f(t,\bar y(t),u(t))-f[t])(\tilde Z,V_u(t),V_u(t))
 +o(\epsilon^2),\quad t\in(0,T],
\cr V_u(0)-X_u(0)=0.
}\end{equation}
Taking $Z=V_u(t)-X_u(t)$ in the above system, and using $(C2)$, $(C3)$, (\ref{32}) and Lemma \ref{17}, we can find a positive constant $C>0$ such that
$$\frac{\partial}{\partial t}|V_u(t)-X_u(t)|^2
\leq C|V_u(t)-X_u(t)|^2
+CI_{\{\bar u\neq u\}}(t)|V_u(t)|^2+o(\epsilon^2),
$$
where $I_{\{\bar u\neq u\}}$
 is the indicator function of the set $\{\tau\in[0,T];\ \;\bar u(\tau)\neq u(\tau)\}$, and the boundness of curvature tensor $R$ along $\bar y(\cdot)$ is used. By integrating the above inequality over $[0,t]$, applying Gronwall's inequality of integral form and (\ref{32}), we can get the first estimate in (\ref{37}).

We subtract (\ref{477}) from (\ref{39}), and get
$$\cases{\langle\nabla_{\dot{\bar y}(t)}(V_u-X_u-Y_u),Z\rangle=\langle Z ,\nabla_{V_u(t)-X_u(t)-Y_u(t)}f(t,\cdot,\bar u(t))\rangle
\cr \quad\quad\quad\quad+\nabla_x\Big(f(t,\bar y(t),u(t))-f[t]\Big)(\tilde Z,V_u(t)-X_u(t))
\cr \quad\quad\quad\quad-\frac{1}{2}R(Z,V_u(t)-X_u(t),f[t],X_u(t))-\frac{1}{2}R(f[t],V_u(t)-X_u(t),Z,V_u(t))
\cr \quad\quad\quad\quad+\frac{1}{2}\nabla^2_xf[t](\tilde Z,V_u(t)-X_u(t),V_u(t))
+\frac{1}{2}\nabla^2_xf[t](\tilde Z,X_u(t),V_u(t)-X_u(t))
\cr \quad\quad\quad\quad-\frac{1}{4}\nabla_2^3\nabla_1\rho^2(\bar y(t),\bar y(t))\Big(Z, f(t,\bar y(t),u(t))-f[t],V_u(t),V_u(t)\Big)
\cr\quad\quad\quad\quad+\frac{1}{2}\nabla_x^2(f(t,\bar y(t),u(t))-f[t])(\tilde Z, V_u(t),V_u(t))+o(\epsilon^2),\quad t\in(0,T],
\cr V_u(0)-X_u(0)-Y_u(0)=0,
}$$
which implies
$$\begin{array}{ll}
&|\nabla_{\dot{\bar y}(t)}(V_u-X_u-Y_u)|
\\ \leq & C|V_u(t)-X_u(t)-Y_u(t)|+CI_{\{u\neq \bar u\}}(t)\Big(|V_u(t)-X_u(t)|+|V_u(t)|^2\Big)
\\ &+C|V_u(t)-X_u(t)|(|V_u(t)|+|X_u(t)|)+o(\epsilon^2),
\end{array}$$
where we have used  $(C1)$, $(C2)$, Lemma \ref{425}, the boundness of curvature tensor $R$ along $\bar y(\cdot)$. Moreover, by applying  (\ref{32}), the first estimate in (\ref{37}), Lemma \ref{480} and (\ref{80}), we can get
$$\begin{array}{ll}
&\frac{\partial }{\partial t}|V_u(t)-X_u(t)-Y_u(t)|^2
\\ \leq &2|\nabla_{\dot{\bar y}(t)}(V_u-X_u-Y_u)||V_u(t)-X_u(t)-Y_u(t)|
\\ \leq & C|V_u(t)-X_u(t)-Y_u(t)|^2+ C|V_u(t)-X_u(t)-Y_u(t)||V_u(t)-X_u(t)|(|V_u(t)|+|X_u(t)|
)
\\&+CI_{\{u\neq \bar u\}}(t)\Big(|V_u(t)-X_u(t)|+|V_u(t)|^2\Big)|V_u(t)-X_u(t)-Y_u(t)|
\\ &+o(\epsilon^2)|V_u(t)-X_u(t)-Y_u(t)|
\\ \leq & C|V_u(t)-X_u(t)-Y_u(t)|^2+CI_{\{u\neq \bar u\}}(t)(|Y_u(t)|^2+|V_u(t)|^4)
+o(\epsilon^4).
\end{array}$$
By integrating the above inequality over $[0,t]$ and applying Gronwall's inequality and (\ref{32}), we can get the second estimate in (\ref{37}).  $\Box$

\subsection{Classical Variation}

In this subsection, we assume that $U$ is an open subset in $\R^m$. Similarly to Proposition \ref{31}, we have the following simple result.

\begin{pro}\label{439} Assume that the conditions $(C1)$, $(C2)$ and $(C3)$ hold.
For any $\bar u(\cdot)\in L^2(0,T;U)$ and $v(\cdot)\in L^2(0,T;R^m)$, set
$$
u^\epsilon(\cdot)=\bar u(\cdot)+\epsilon v(\cdot),\quad \epsilon\geq 0.
$$
Denote by $y^\epsilon(\cdot)$ the solution to (\ref{25}) with control $u^\epsilon(\cdot)$. In particular, we denote by $\bar y(\cdot)$ the solution to (\ref{25}) with control $\bar u(\cdot)$. For any $\epsilon>0$ being small enough, we define a vector field along $\bar y(\cdot)$ as follows:
\begin{equation}\label{438}
V_\epsilon(t)=\exp_{\bar y(t)}^{-1}y^\epsilon(t),\quad t\in[0,T].
\end{equation}
Let $V(\cdot)$ and $Y(\cdot)$ be respectively the vector fields along $\bar y(\cdot)$  solving:
\begin{equation}\label{527}
\cases{\nabla_{\dot{\bar y}(t)}V(Z)=\nabla_xf[t](Z,V(t))+ \nabla_uf[t](Z,v(t)),\quad a.e. \,t\in(0,T],\;\forall\; Z\in T^*M,\cr V(0)=0,
}
\end{equation}
and
\begin{equation}\label{528}
\cases{\nabla_{\dot{\bar y}(t)}Y(Z)=\nabla_xf[t](Z,Y(t))+\nabla_x\nabla_uf[t](Z,v(t),V(t))-\frac{1}{2}R(\tilde Z,V(t),\dot{\bar y}(t),V(t))
\cr \quad\quad\quad+\frac{1}{2}\nabla_x^2f[t](Z,V(t),V(t))+\frac{1}{2}
\nabla_u^2f[t](Z, v(t),v(t)),\quad a.e. \,t\in(0,T], \forall\; Z\in T^*M,\cr Y(0)=0.
}\end{equation}
Then, we have
\begin{equation}\label{427}
V_\epsilon(t)=\epsilon V(t)+\epsilon^2Y(t)+o(\epsilon^2),\quad\forall\; t\in[0,T].
\end{equation}
\end{pro}

\textbf{Proof}\,\,\, Similarly to Lemma \ref{480}, when $\epsilon> 0$ is small enough, $\rho(y^\epsilon(t),\bar y(t))=O(\epsilon)$ for any $t\in[0,T]$. Thus, one can define the vector field (\ref{438}).

Following a similar argument as that in the proof of (\ref{35}), we can obtain,
for any $Z\in T_{\bar y(t)}M$,
$$\begin{array}{ll}
&\langle \nabla_{\dot{\bar y}(t)}V_\epsilon,Z\rangle
\\=&\nabla_x f[t](Z,V_\epsilon(t))+\nabla_uf[t](Z,v(t))\epsilon+\nabla_x
\nabla_uf[t](Z,v(t),V_\epsilon(t))\epsilon-\frac{1}{2}R(Z,
V_\epsilon(t),\dot{\bar y}(t),V_\epsilon(t))
\\ &+\frac{1}{2}\nabla_x^2f[t](Z,V_\epsilon(t),V_\epsilon(t))+
\frac{1}{2}\nabla_u^2f[t](Z,v(t),v(t))\epsilon^2+
o(\epsilon^2).
\end{array}$$
Dividing the above identity by $\epsilon$, we obtain that
\begin{equation}\label{537}\lim_{\epsilon\to 0^+}\frac{V_\epsilon(t)}{\epsilon}=V(t)\end{equation}
 uniformly in $t\in[0,T]$, where
$V(\cdot)$ is the solution to (\ref{527}).

Similarly, one can check that
$$
\lim_{\epsilon\to 0^+}\frac{V_\epsilon(t)-\epsilon V(t)}{\epsilon^2}=Y(t)
$$
uniformly in $t\in[0,T]$, where $Y(\cdot)$ is a vector field along $\bar y(\cdot)$ solving (\ref{528}).
Thus, the proof is concluded. $\Box$

\setcounter{equation}{0}
\section{Proof of the Main Results}
\def\theequation{6.\arabic{equation}}

This section is addressed to proving our main results, i.e., Theorems \ref{66}--\ref{453}.

\subsection{Proof of Theorem \ref{66}}
The proof of Pontryagin-type maximum principle (\ref{3}) is given in \cite[p. 183, Theorem 12.13]{as}. This subsection aims at  proving the second order necessary condition. We divide the proof into $3$ steps.

\medskip

\textbf{Step 1}  In this step, we will give the second order expansion of $J(\cdot)$ around the  optimal control $\bar{u}(\cdot)$.
Fix any $u(\cdot)\in\U_{ad}$ and  $\epsilon>0$. Define
\begin{equation}\label{492}
u^\epsilon(t)=\bar u(t)I_{[0,T]\setminus E_\epsilon}(t)+u(t)I_{E_\epsilon}(t),
\end{equation}
where $E_\epsilon\subset [0,T]$ is a measurable subset  of $[0,T]$ with $|E_\epsilon|=\epsilon$. Let $y^\epsilon(\cdot)$ be the solution to (\ref{25}) corresponding to the control $u^\epsilon(\cdot)$.
Denote by $X^\epsilon(\cdot)$ and $Y^\epsilon(\cdot)$ the solutions to (\ref{227}) and (\ref{477}) with $u(\cdot)=u^\epsilon(\cdot)$, respectively.
By Proposition \ref{31}, there exists an $\epsilon_0>0$ such that, for all $\epsilon<\epsilon_0$,  we can define
$$
V^\epsilon(t)\equiv \exp_{\bar y(t)}^{-1}y^\epsilon(t),\quad \forall\; t\in[0,T],
$$
and
\begin{equation}\label{485}
X^\epsilon(t)=O(\epsilon),\quad Y^\epsilon(t)=O(\epsilon^2),\quad
V^\epsilon(t)=X^\epsilon(t)+Y^\epsilon(t)+o(\epsilon^2),\quad
\forall\; t\in[0,T].
\end{equation}
Similar to (\ref{481}) and (\ref{482}), for $t\in[0,T]$, let $\tilde V^\epsilon(t)$ be a vector at $\bar y(t)$
satisfying
\begin{equation}\label{501}
\tilde V^\epsilon(t)\equiv\cases{\frac{1}{\rho(\bar y(t),y^\epsilon(t))}V^\epsilon(t),\quad\textrm{if}\,
|V^\epsilon(t)|\neq0;
\cr 0,\quad\quad\quad\quad\quad\quad\quad\textrm{if}\,
|V^\epsilon(t)|=0,}
\end{equation}
 and let $\gamma(\cdot;t)$ be the shortest geodesic connecting $\bar y(t)$ and $y^\epsilon(t)$, given by
 \begin{equation}\label{502}
 \gamma(\theta;t)=\exp_{\bar y(t)}(\theta \tilde V^\epsilon(t)),\quad \theta\in[0,\rho(\bar y(t),y^\epsilon(t))].
 \end{equation}
 In particular, in the case that $|V^\epsilon(t)|\neq0$, $\gamma(\cdot;t)$ is parameterized by the arc length. Recalling (\ref{434}),
 and applying (\ref{203}), Tayor's Theorem, Lemma \ref{480}, (\ref{501}) and (\ref{502}), we obtain that
$$\label{53}
\begin{array}{ll}
 & J(u^\epsilon(\cdot))-J(\bar u(\cdot))
\\ =&\int_0^T\Big(f^0(t,y^\epsilon(t),u^\epsilon(t))-f^0(t,\bar y(t),u^\epsilon(t))+f^0(t,\bar y(t),u^\epsilon(t))-f^0[t]\Big)dt
\\ =&\int_0^T\Big(f^0(t,\gamma(\rho(\bar y(t),y^\epsilon(t));t),u^\epsilon(t))-f^0(t,\gamma(0;t),
u^\epsilon(t))
\\ &+f^0(t,\bar y(t),u^\epsilon(t))-f^0[t]\Big)dt
\\ =&\int_0^T\Big(\frac{\partial}{\partial s}\Big|_0f^0(t,\gamma(s;t),u^\epsilon(t))\rho(\bar y(t),y^\epsilon(t))
+\frac{1}{2}\frac{\partial^2}{\partial s^2}\Big|_0f^0(t,\gamma(s;t),u^\epsilon(t))\rho^2(\bar y(t),y^\epsilon(t))
\\ &+f^0(t,\bar y(t),u^\epsilon(t))-f^0[t]+o\Big(\rho^2(\bar y(t),y^\epsilon(t))\Big)\Big)dt
\\ =&\int_0^T\Big(\langle\nabla_xf^0(t,\bar y(t),u^\epsilon(t)),\tilde V^\epsilon(t)\rangle\rho(\bar y(t),y^\epsilon(t))
\\ &+\frac{1}{2}\nabla_x^2f^0(t,\bar y(t),u^\epsilon(t))(\tilde V^\epsilon(t),\tilde V^\epsilon(t))\rho^2(\bar y(t),y^\epsilon(t))
+f^0(t,\bar y(t),u^\epsilon(t))
\\ &-f^0[t]+o\Big(\rho^2(\bar y(t),y^\epsilon(t))\Big)\Big)dt
\\ =&\int_0^T\Big(\langle\nabla_xf^0[t],V^\epsilon(t)\rangle+\frac{1}{2}\nabla_x^2f^0[t](V^\epsilon(t),V^\epsilon(t))
+f^0(t,\bar y(t),u^\epsilon(t))-f^0[t]
\\ &+\frac{1}{2}\nabla_x^2\Big(f^0(t,
\bar y(t),u^\epsilon(t))
-f^0(t,\bar y(t),\bar u(t))\Big)(V^\epsilon(t),V^\epsilon(t))
\\ & +\langle \nabla_xf^0(t,\bar y(t),u^\epsilon(t))
- \nabla_x f^0[t],V^\epsilon(t)\rangle+o(\rho^2(\bar y(t),y^\epsilon(t)))\Big)dt.
\end{array}
$$
Recalling (\ref{478}) and (\ref{485}), we further get
\begin{equation}\label{54}
J(u^\epsilon(\cdot))-J(\bar u(\cdot))=I_1+I_2+o(\epsilon^2),\end{equation}
where
\begin{eqnarray}\label{498}I_1&=& \int_0^T\Big(\langle\nabla_xf^0[t],X^\epsilon(t)\rangle+f^0(t,\bar y(t),u^\epsilon(t))-f^0[t]\Big)dt,
\\ \label{499}I_2&=&\int_0^T\Big(\langle\nabla_xf^0[t],Y^\epsilon(t)\rangle+\langle\nabla_xf^0(t,\bar y(t),u^\epsilon(t))-\nabla_xf^0[t],X^\epsilon(t)\rangle\nonumber
\\&& +\frac{1}{2}\nabla_x^2f^0[t](X^\epsilon(t),X^\epsilon(t))\Big)dt.
\end{eqnarray}

 In what follows, we will rewrite $I_1$ and $I_2$ by the dual variables and the perturbed control $u^\epsilon(\cdot)$. Recalling the first order variational equation (\ref{227}), the first order dual equation (\ref{64}) and (\ref{65}) with $\nu=-1$ and $\psi_1=0$, we have, via integration by parts over $[0,T]$,
\begin{equation}\label{55}\begin{array}{lll}
I_1&=&\int_0^T\Big(\nabla_{\dot{\bar y}(t)}\psi(X^\epsilon(t))+\nabla_{X^\epsilon(t)}f(t,\cdot,
\bar u(t))(\psi(t))
+f^0(t,\bar y(t),u^\epsilon(t))-f^0[t]\Big)dt
\\ &=& \int_0^T\Big(\psi(t)\Big(-\nabla_{\dot{\bar y}(t)}X^\epsilon+\nabla_{X^\epsilon(t)}f(t,\cdot,\bar u(t))\Big)+f^0(t,\bar y(t),u^\epsilon(t))-f^0[t]\Big)dt
\\ &=& \int_0^T\Big(\psi(t)\Big(f[t]-f(t,\bar y(t),u^\epsilon(t))\Big)+f^0(t,\bar y(t),u^\epsilon(t))-f^0[t]\Big)dt
\\ &=&\int_0^T\Big(H(t,\bar y(t),\psi(t),\bar u(t))-H(t,\bar y(t),\psi(t),u^\epsilon(t))\Big)dt.
\end{array}\end{equation}

Using (\ref{64}) with $\nu=-1$ and $\psi_1=0$, the second order dual equation (\ref{70}), the second order variational equation (\ref{477}) and integration by parts over $[0,T]$, via   (\ref{469}) and (\ref{468}), we can get
\begin{equation}\label{57}\begin{array}{lll}
I_2&=&\int_0^T\Big(\nabla_{\dot{\bar y}(t)}\psi(Y^\epsilon(t))+\psi(t)(\nabla_{Y^\epsilon(t)}
f(t,\cdot,\bar u(t)))+\langle\nabla_xf^0(t,\bar y(t),u^\epsilon(t))
\\ &&-\nabla_xf^0[t],X^\epsilon(t)\rangle+\frac{1}{2}\nabla_x^2f^0[t](X^\epsilon(t),X^\epsilon(t))\Big)dt
\\ &=& \int_0^T\Big(\psi(t)\Big(-\nabla_{\dot{\bar y}(t)}Y^\epsilon+\nabla_{Y^\epsilon(t)}f(t,\cdot,\bar u(t))\Big)+\langle\nabla_xf^0(t,\bar y(t),u^\epsilon(t))
\\ &&-\nabla_xf^0[t],X^\epsilon(t)\rangle+\frac{1}{2}\nabla_x^2f^0[t](X^\epsilon(t),X^\epsilon(t))\Big)dt
\\&=&\int_0^T\Big(\frac{1}{2}R(\tilde\psi(t),X^\epsilon(t),
f[t],X^\epsilon(t))
-\langle\nabla_{X^\epsilon(t)}(f(t,\cdot,u^\epsilon(t))
-f(t,\cdot,\bar u(t))),\tilde\psi(t)\rangle
\\ &&-\frac{1}{2}\nabla_x^2f[t](\psi(t),X^\epsilon(t),X^\epsilon(t))
+\langle\nabla_xf^0(t,\bar y(t),u^\epsilon(t))-\nabla_xf^0[t],X^\epsilon(t)\rangle
\\ &&+\frac{1}{2}\nabla_x^2f^0[t](X^\epsilon(t),X^\epsilon(t)) \Big)dt
\\ &=&I_3+I_4,
\end{array}\end{equation}
where $\tilde\psi$ is the dual vector of $\psi$,
\begin{equation}\label{85}\begin{array}{ll}
I_3
=&\frac{1}{2}\int_0^T\Big(R(\tilde\psi(t),X^\epsilon(t),
f[t],X^\epsilon(t))
-\nabla_x^2H(t,\bar y(t),\psi(t),\bar u(t))(X^\epsilon(t),X^\epsilon(t))\Big)dt,
\end{array}\end{equation}
and
\begin{equation}\label{486}
I_4=\int_0^T\langle\nabla_xH(t,\bar y(t),\psi(t),\bar u(t))
-\nabla_xH(t,\bar y(t),\psi(t),u^\epsilon(t)),X^\epsilon(t)\rangle dt.
\end{equation}

Next, we shall re-write $I_3$ and $I_4$ in terms of local coordinates along $\bar y(\cdot)$.
For this purpose, let $\{e_1,e_2,\cdots,e_n\}$ be an orthonormal basis at $\bar{y}(0)=y_0$, and $\{d_1,d_2,\cdots,d_n\}$ be the dual basis of $\{e_1,e_2,\cdots,e_n\}$. For $t\in[0,T]$, we define
\begin{equation}\label{277}
e_i(t)\equiv L^{\bar{y}(\cdot)}_{y_0\bar{y}(t)}e_i,
\quad d_i(t)\equiv L^{\bar{y}(\cdot)}_{y_0\bar{y}(t)}d_i,
\quad i=1,2,\cdots,n.
\end{equation}
By (\ref{205}) and (\ref{82}), we have
\begin{equation}\label{89}
\langle e_i(t),e_j(t)\rangle=
d_i(t)(e_j(t))=
\delta_i^j,\quad i,j=1,2,\cdots,n,\,\,\,t\in[0,T].
\end{equation}
Set
\begin{equation}\label{293}
X^\epsilon(t)=\sum_{i=1}^nX^\epsilon_i(t)e_i(t),\quad \vec{X}^\epsilon(t)=(X^\epsilon_1(t),\cdots,
X^\epsilon_n(t))^\top.
\end{equation}
Inserting the above equalities into
(\ref{85}), we get
\begin{equation}\label{77}\begin{array}{ll}
I_3=&\frac{1}{2}tr\int_0^T (M(t)-\vec H(t))\vec{X}^\epsilon(t)\vec{X}^{\epsilon\top}(t)dt,
\end{array}\end{equation}
where
\begin{eqnarray}\label{494}
M(t)=(M_{ik}(t))\equiv\Big(R(\tilde\psi(t),e_i(t),f[t],e_k(t))\Big),\;t\in[0,T],\quad\quad
\\\label{72}\vec H(t)=(H_{ik}(t))\equiv\Big(\nabla_x^2H(t,\bar y(t),\psi(t),\bar u(t))(e_i(t),e_k(t))\Big),\;t\in[0,T].
\end{eqnarray}

Let us derive the differential equation satisfied by the matrix $\vec X^\epsilon(t)\vec X^{\epsilon\top}(t)$.
Inserting (\ref{293}) into (\ref{227}), via the parallel translation of $e_i(\cdot)$ ($i=1,\cdots,n$) along $\bar y(\cdot)$, we can get
\begin{equation}\label{84}\cases{
\dot{\vec {X}^\epsilon}(t)=F(t)\vec X^\epsilon(t)+F_1(t,u^\epsilon(t)),\quad t\in(0,T],\cr
\vec X^\epsilon(0)=0,
}\end{equation}
where
\begin{eqnarray}\label{286}
&
\ds
 F(t)=\Big(\Big\langle\nabla_{e_j(t)}f(t,\cdot,\bar{u}
 (t)),e_i(t)\Big\rangle\Big)=\nabla_xf[t](d_i(t),e_j(t))\equiv \Big(F_{ij}(t)\Big),
\\
[3mm]
&
\ds
\label{287}F_1(t,u)=\left(\begin{array}{c}\langle f(t,\bar{y}(t),u)-f[t],e_1(t)\rangle
\\
[3mm]
 \vdots
\\
[3mm]
\ds
  \langle f(t,\bar{y}(t),u)-f[t],e_n(t)\rangle\end{array}\right)\equiv (F_1^1(t,u),\cdots, F_1^n(t,u))^\top.\qquad
\end{eqnarray}
Hence, we get the equation for $\vec X^\epsilon(t)\vec X^{\epsilon\top}(t)$ as follows:
\begin{equation}\label{74}
\cases{\ds\frac{d}{dt}(\vec{X}^\epsilon
\vec{X}^{\epsilon\top})(t)=F(t)(\vec{X}^\epsilon
\vec X^{\epsilon\top})(t)+(\vec{X^\epsilon}
\vec{X}^{\epsilon\top})(t)F(t)^\top+F_1(t,
u^\epsilon(t))\vec{X^\epsilon}(t)^\top \cr
\ns\quad\quad\quad\quad\quad\quad+\vec{X^\epsilon}(t)
F_1(t,u^\epsilon(t))^\top,\qquad t\in(0,T],
\cr \ns(\vec{X}^\epsilon\vec{X}^{\epsilon\top})(0)=0.}
\end{equation}

We now derive the local form of the second order dual equation (\ref{70}).
Recalling (\ref{89}), we can express the solution $w$ to (\ref{70}) by
\begin{equation}\label{75}
w(t)=\sum_{i,j=1}^nw_{ij}(t)d_i(t)\otimes d_j(t),
\end{equation}
where $w_{ij}(t)\equiv w(t)(e_i(t),e_j(t))$ for any $t\in[0,T]$. The transpose of $w(t)$ can be rewritten as
\begin{equation}\label{88}
w^\top(t)=\sum_{i,j=1}^nw_{ji}(t)d_i(t)\otimes d_j(t),\quad \forall\; t\in[0,T].
\end{equation}
Recalling (\ref{286}), we have
\begin{equation}\label{73}
\nabla_xf[t]=\sum_{i,j=1}^nF_{ij}(t)e_i(t)\otimes d_j(t).
\end{equation}
Thanks to (\ref{494}) and (\ref{72}), we have
\begin{eqnarray}\label{483}
R(\tilde\psi(t),\cdot,f[t],\cdot)=\sum_{i,j=1}^nM_{ij}(t)d_i(t)\otimes d_j(t),
\\ \label{484}\nabla_x^2H(t,\bar y(t),\psi(t),\bar u(t))=\sum_{i,j=1}^nH_{ij}(t)d_i(t)\otimes d_j(t).
\end{eqnarray}
Set
\begin{equation}\label{86}
W(t)=(w_{ij}(t)),\quad\forall\; t\in[0,T].
\end{equation}
Recalling (\ref{277}) and (\ref{89}), we insert (\ref{75}), (\ref{73}), (\ref{483}) and (\ref{484}) into (\ref{70}), and get
\begin{equation}\label{83}
\cases{\dot{W}(t)+F(t)^\top W(t)+W(t)F(t)-M(t)+\vec{H}(t)=0,\,\,\,t\in[0,T),
\cr W(T)=0.
}\end{equation}
Applying  (\ref{74}), (\ref{83}) to (\ref{77}), via integration by parts over $[0,T]$, we get
\begin{equation}\label{490}\begin{array}{ll}
&I_3
\\ =& \frac{1}{2}tr\int_0^T\Big(\dot W(t)+F(t)^\top W(t)+W(t)F(t)\Big)\vec X^\epsilon(t)\vec X^{\epsilon\top}(t)dt
\\ =&\frac{1}{2}tr\int_0^T\Big\{W(t)\Big(-\frac{d}{dt}(\vec X^\epsilon\vec X^{\epsilon\top})(t)+F(t)(\vec X^\epsilon\vec X^{\epsilon\top})(t)\Big)+F(t)^\top W(t)(\vec X^\epsilon\vec X^{\epsilon\top})(t)\Big\}dt
\\ =&\frac{1}{2}tr\int_0^T\Big\{-W(t)(\vec X^\epsilon\vec X^{\epsilon\top})(t)F(t)^\top-W(t)F_1(t,u^\epsilon(t))\vec X^{\epsilon\top}(t)-W(t)\vec X^\epsilon(t)F_1(t,u^\epsilon(t))^\top
\\ &+F(t)^\top W(t)\vec X^\epsilon(t)\vec X^{\epsilon\top}(t)\Big\}dt
\\=&-\frac{1}{2}tr\int_0^T(W(t)+W^\top(t))F_1(t,u^\epsilon
(t))\vec X^{\epsilon\top}(t)dt
\\ =& -\frac{1}{2}\int_0^T\Big((W(t)+W^\top(t))F_1(t,u^\epsilon
(t))\cdot \vec X^\epsilon(t)dt,
\end{array}\end{equation}
where
we have used the  properties:
$
tr (AB)=tr (BA)$ and $ tr A=tr (A^\top)
$ for any $n\times n$ matrixes $A$ and $B$.

In order to express $I_4$ in the local form, for any $u\in U$ and $i=1,
\cdots,n$, set
\begin{eqnarray}\label{487}
&\ds\delta\partial_i H(t,u)=\Big(\nabla_xH(t,\bar{y}(t),\psi(t),\bar{u}(t))-
\nabla_xH(t,\bar{y}(t),\psi(t),u)\Big)(e_i(t)),
\\[2mm] \label{489} &\ds\delta\partial H(t,u)=\Big(\delta\partial_1H(t,u),\cdots,\delta\partial_n
H(t,u)\Big)^\top.
\end{eqnarray}
Recalling (\ref{293}), we have
\begin{equation}\label{491}
I_4=\int_0^T\delta\partial H(t,u^\epsilon(t))\cdot \vec{X}^\epsilon(t)dt.
\end{equation}

Now, we re-write $\vec X^\epsilon$ in terms of $u^\epsilon(\cdot)$.
 Assume that $\tilde\Phi(\cdot)=(\phi_{ij}(\cdot))$ is an $\R^{n\times n}$-valued function satisfying
\begin{equation}\label{455}
\cases{\dot{\tilde\Phi}(t)=F(t)\tilde\Phi(t),\,\,\,t\in(0,T],\cr
\tilde\Phi(0)=I,}
\end{equation}
where $I\in \R^{n\times n}$ is the $n\times n$ identity matrix. Recalling (\ref{84}), we have
\begin{equation}\label{PPP455}
\vec X^\epsilon(t)=\int_0^t\tilde\Phi(t){\tilde\Phi}^{-1}(s)F_1(s,
u^\epsilon(s))ds,
\end{equation}
where ${\tilde\Phi}^{-1}(\cdot)=(\psi_{ij}(\cdot))$ is the inverse of $\tilde \Phi(\cdot)$.

 Inserting (\ref{490}) and (\ref{491}) into (\ref{57}) via (\ref{PPP455}), and recalling (\ref{492}), (\ref{55}), (\ref{287}), (\ref{487}) and (\ref{489}), we can rewrite (\ref{54}) by
 \begin{equation}\label{548}\begin{array}{ll}
& J(u^\epsilon(\cdot))-J(\bar u(\cdot))
\\ =& \int_{E_\epsilon}\Big(H(t,\bar y(t),\psi(t),\bar u(t))-H(t,\bar y(t),\psi(t),u(t))\Big)dt+\int_{E_\epsilon}\Big\{-\frac{1}{2}(W(t)
 \\ &+W^\top(t))F_1(t,u(t))+
 \delta\partial H(t,u(t))\Big\}\cdot\tilde\Phi(t)\int_0^t{\tilde\Phi}^{-1}(s)F_1(s,u^\epsilon(s))dsdt+o(\epsilon^2).
\end{array}\end{equation}

\medskip

\textbf{Step 2}\,\,\, In order to get the integral form of second order necessary condition, we apply \cite[Corollary 3.8, p. 144]{yz} to choose suitable measurable set $E_\epsilon$ in (\ref{548}). Let $u(\cdot)$ be such that $u(t)\in\widetilde U(t)$, a.e. $t\in [0,T]$. Let $E_\epsilon\subset[0,T]$ with $|E_\epsilon|=\epsilon$ be such that
$$\begin{array}{l}
\int_{[0,t]\cap E_\epsilon}{\tilde\Phi}^{-1}(s)F_1(s,u(s))ds=\epsilon\int_0^t{\tilde\Phi}^{-1}(s)F_1(s,u(s))ds+\eta_1(t),
\\[2mm] \int_{[0,t]\cap E_\epsilon}\Big\{-\frac{1}{2}(W(\tau)+W^\top(\tau))F_1(\tau,u(\tau))\\[2mm] \quad+\delta\partial H(\tau,u(\tau))\Big\}\cdot\tilde\Phi(\tau)
 \int_0^\tau{\tilde\Phi}^{-1}(s)F_1(s,u(s))ds d\tau
\\[2mm]  =\epsilon\int_0^t\Big\{-\frac{1}{2}(W(\tau)+W^\top(\tau))F_1(\tau,u(\tau))
 \\[2mm] \quad +\delta\partial H(\tau,u(\tau))\Big\}\cdot\tilde\Phi(\tau)\int_0^\tau{\tilde\Phi}^{-1}(s)F_1(s,u(s))ds d\tau+\eta_2(t),
\end{array}$$
 where
 $$
 |\eta_i(t)|\leq \epsilon^3,\quad\forall\; t\in[0,T],\,i=1,2.
 $$
 Then, dividing (\ref{548}) by $\epsilon^2$ and taking the limit as $\epsilon\to0$, one can get, via the optimality of $\bar u(\cdot)$,
 \begin{equation}\label{45}
 0\leq \int_0^T\Big\{-\frac{1}{2}(W(t)+W^\top(t))F_1(t,u(t))+\delta\partial H(t,u(t))\Big\}\cdot \tilde \Phi(t)\int_0^t\tilde \Phi^{-1}(s)F_1(s,u(s))dsdt.
 \end{equation}
 Set
 \begin{equation}\label{454}
 \Phi(t)=\sum_{i,j=1}^n\phi_{ij}(t)e_i(t)\otimes e_j(t),\quad \Phi_1(t)=\sum_{i,j=1}^n\psi_{ij}(t)d_i(t)\otimes d_j(t),\quad t\in[0,T].
 \end{equation}
 Recalling (\ref{76}), (\ref{73}) and (\ref{455}), one can check, via (\ref{277}) and (\ref{82}), that $\Phi(\cdot)$ and $\Phi_1(\cdot)$ solve (\ref{459}) and (\ref{458}) respectively.

Recalling  (\ref{287}), (\ref{75}), (\ref{88}), (\ref{86}), (\ref{487}), (\ref{489}) and (\ref{454}), one can get (\ref{457}) from (\ref{45}), via (\ref{277}) and (\ref{89}).

\medskip

\textbf{Step 3} In this step, we shall get the pointwise form of the second order necessary condition. Assume that $U$ is a Polish space. Recalling (\ref{45}), define
 \begin{equation}\label{431}\begin{array}{l}
 G_1(t,u)=\tilde\Phi(t)^\top\Big\{-\frac{1}{2}(W(t)+W^\top(t))F_1(t,u)+\delta \partial H(t,u)\Big\},
 \\ G_2(t,u)=\tilde\Phi^{-1}(t)F_1(t,u),\quad (t,u)\in[0,T]\times U.
 \end{array}\end{equation}
 Then, by (\ref{45}), for any $u(\cdot)\in \U_{ad}$ and $u(t)\in\widetilde U(t)$,  a.e. $t\in [0,T]$, one has
 \begin{equation}\label{433}
 0\leq \int_0^T G_1(t,u(t)) \cdot \int_0^t G_2(s,u(s))ds  dt.
 \end{equation}

Applying the same argument as that in the proof of \cite[Theorem 4.3]{l}, we can get, from the above inequality and the conditions (C1) and (C2),
$$
G_1(t,u(t))\cdot G_2(t,u(t))\geq 0\quad a.e.\,t\in[0,T].
$$
Following the same argument as that in the proof of \cite[Theorem 4.3]{l} again, we can get, from the above inequality and the assumption that
 $U$ is a Polish space,
$$
G_1(t,v)\cdot G_2(t,v)\geq 0\quad \forall\; v\in\widetilde U(t),
$$
for a.e. $t\in[0,T]$. Recalling (\ref{431}), the above inequality  implies
\begin{equation}\label{432}
0\leq \Big\{-\frac{1}{2}\Big(W(t)+W^\top(t)\Big)F_1(t,v)+\delta
\partial H(t,v)\Big\}\cdot F_1(t,v),\quad\forall\; v\in \widetilde U(t),\;\, a.e.\,t\in[0,T].
\end{equation}
Combining the above inequality with (\ref{89}), (\ref{287}), (\ref{75}), (\ref{88}), (\ref{487}) and (\ref{489}), we can easily obtain (\ref{71}). $\Box$

\subsection{Proof of Theorem \ref{456}}
For any $v(\cdot)\in\U_{ad}$ with $d(\bar u(\cdot), v(\cdot))=\epsilon$ being small enough,  from (\ref{548}), we have,
$$\begin{array}{ll}
&J(v(\cdot))-J(\bar u(\cdot))
\\ =& \int_0^T\Big(H(t,\bar y(t),\psi(t),\bar u(t))-H(t,\bar y(t),\psi(t), v(t))\Big)dt+\int_0^T\Big\{-\frac{1}{2}(W(t)+W^\top(t))
\\ &F_1(t,v(t))+\delta\partial H(t,v(t))\Big\}\cdot\tilde\Phi(t)\int_0^t\tilde\Phi^{-1}(s)F_1(s,v(s))dsdt+o(\epsilon^2).
\end{array}$$
 Recalling  (\ref{287}), (\ref{75}), (\ref{88}), (\ref{86}), (\ref{487}), (\ref{489}) and (\ref{454}), one can get, via (\ref{277}) and (\ref{89}),
$$\begin{array}{ll}
&J(v(\cdot))-J(\bar u(\cdot))
\\ =& \int_0^T\Big(H(t,\bar y(t),\psi(t),\bar u(t))-H(t,\bar y(t),\psi(t), v(t))\Big)dt-\int_0^T\int_0^t\Big\{\frac{1}{2}(w(t)+w^\top(t))
\\ &\Big(f(t,\bar y(t),v(t))-f[t],
  {\mathcal E}^{21}\Big({\mathcal E}^{21}\Big(\Phi(t)\otimes L_{\bar y(s)\bar y(t)}^{\bar y(\cdot)}\Phi_1(s)\Big)\otimes
 \\ & L_{\bar y(s)\bar y(t)}^{\bar y(\cdot)}(f(s,\bar y(s),v(s))-f[s]) \Big)\Big)
+\Big\langle\nabla_xH(t,\bar y(t),\psi(t),v(t))
\\ & -\nabla_xH(t,\bar y(t),\psi(t),
 \bar u(t)),
  {\mathcal E}^{21}\Big({\mathcal E}^{21}\Big(\Phi(t)\otimes L_{\bar y(s)\bar y(t)}^{\bar y(\cdot)}\Phi_1(s)\Big)\otimes
 L_{\bar y(s)\bar y(t)}^{\bar y(\cdot)}(f(s,\bar y(s),v(s))
 \\ &-f(s,\bar y(s),\bar u(s))) \Big)\Big\rangle\Big\}dsdt+o(\epsilon^2).
\end{array}
$$
Hence, according to (\ref{3}) and (\ref{20}), it is easy to see that, there exists an $\epsilon_1>0$ such that
$J(u^\epsilon(\cdot))-J(\bar u(\cdot))\geq 0$, for all $u^\epsilon(\cdot)\in\U_{ad}$ with $d(u^\epsilon(\cdot),\bar u(\cdot))=\epsilon\leq\epsilon_1$. This completes the proof of Theorem \ref{456}. $\Box$

\subsection{Proofs of Theorems \ref{504} and \ref{453}}
In this subsection, we shall prove Theorems \ref{504} and \ref{453}.

To begin with, let us recall the following known result (\cite[Theorem 20.3, p. 297]{as}):
\begin{lem}\label{533}Let $F:\U\to {\cal M}$ be a continuous mapping having smooth restrictions to finite-dimensional submanifolds of $\U$, where $\U$ is an open subset of a Banach space, and ${\cal M}$ is an $n$-dimensional, smooth differential manifold.
Let $\hat u\in\U$ be a corank one critical point of $F$, i.e. the codimension of $\Im D_{\hat u}F$ (the image of the differential of $F$ at $\hat u$) is equal to $1$. Let $\lambda\in (\Im D_{\hat u}F)^\bot\equiv \{\eta\in T^*_{F(\hat u)}{\cal M}; \ \;\eta(X)=0, \forall\; X\in \Im D_{\hat u}F\subset T_{F (\hat u)}{\cal M}\} $, $\lambda\neq 0$. If the quadratic form $\lambda Hess_{\hat u}F: \ker D_{\hat u}F\times \ker D_{\hat u}F\to \R$ is sign-indefinite,
 which is defined by
 \begin{equation}\label{539}
\lambda Hess_{\hat u}F(v,v)\equiv \lambda\Big(\frac{d^2}{d\epsilon^2}\Big|_{0}F(\varphi
(\epsilon))\Big),  \quad \forall\; v\in Ker D_{\hat u}F\subset T_{\hat u}\U,
\end{equation}
with $\varphi: [0,\epsilon_0)\to \U$ ($\epsilon_0>0$) satisfying $\varphi(0)=\hat u$ and $\frac{d}{d\epsilon}\Big|_{0}\varphi(\epsilon)=v$,
then $F$ is locally open at $\hat u$, i.e.  $F(\hat u)\in int F(O_{\hat u})$ for any neighborhood $O_{\hat u}\subset \U$ of $\hat u$.
\end{lem}

We are now in a position to prove Theorem \ref{504}.

\textbf{Proof of Theorem \ref{504}}\quad
Similarly to \cite{as}, we introduce an extended system associated to problem (\ref{25}) and  (\ref{26})  as follows:
\begin{equation}\label{531}
\cases{\frac{d}{dt}\left(y_0(t)\atop y(t)\right)=\left(f^0(t,y(t),u(t))\atop f(t,y(t),u(t))\right),\quad u(t)\in U,\; a.e.\;t\in(0,T),\cr \left(y_0(0)\atop y(0)\right)
=\left(0\atop y_0\right).}
\end{equation}
Define the  endpoint mapping $E:\;L^2(0,T;U)\to \R\times M$ for the system (\ref{531}) by
\begin{equation}\label{540}
E(u(\cdot))\equiv \left(J(u(\cdot))\atop y(T;u(\cdot))\right),
\end{equation}
where $y(\cdot;u(\cdot))$ is the solution to (\ref{25}) associated to the control $u(\cdot)$.
Then, we define the attainable set of (\ref{531}) at time $T$ by
$$
\A\equiv\{ E(u(\cdot));\ \; u(\cdot)\in L^2(0,T;U)\}.
$$
Since $\bar u(\cdot)$ is optimal for \textbf{Problem II}, we have (\cite[Section 12.4, p. 179]{as})
\begin{equation}\label{532}E(\bar u(\cdot))\in\partial \A.\end{equation}

For any $v(\cdot)\in L^2(0,T;\R^m)$, let $y^\epsilon(\cdot)$ be the solution to (\ref{25}) corresponding to the control $\bar u(\cdot)+\epsilon v(\cdot)$ with small $\epsilon\geq 0$.

Firstly, we claim that
\begin{equation}\label{534}
\frac{\partial }{\partial\epsilon}\Big|_0y^\epsilon(t)=V(t),\quad\frac{\partial^2}{\partial \epsilon^2}\Big|_0
y^\epsilon(t)= 2Y(t),\quad t\in[0,T],
\end{equation}
where $V(\cdot)$ and $Y(\cdot)$ are the solutions to (\ref{527}) and (\ref{528}), respectively.

In fact, for any $t\in[0,T]$, set
\begin{equation}\label{535}
\tilde V_\epsilon(t)\equiv \cases{  \frac{V_\epsilon(t)}{\rho(\bar y(t),y^\epsilon(t))},\quad \textrm{if}\ \rho(\bar y(t),y^\epsilon(t))>0,\cr 0,\quad\quad\quad\quad\quad \textrm{if}\ \rho(\bar y(t),y^\epsilon(t))=0,}
\end{equation}
where $V_\epsilon(t)\equiv \exp_{\bar y(t)}^{-1}y^\epsilon(t)$ with $\epsilon>0$ being small enough.
Define
\begin{equation}\label{536}\eta(\theta;\epsilon)\equiv \exp_{\bar y(t)}(\theta \tilde V_\epsilon(t)),\quad \theta\in[0,\rho(\bar y(t),y^\epsilon(t)].
\end{equation} Then, we have $\eta(0;\epsilon)=\bar y(t)$ and $\eta(\rho(\bar y(t),y^\epsilon(t));\epsilon)=y^\epsilon(t)$.
For any $h\in C^\infty(M)$, applying Taylor's expansion, Lemma \ref{317}, (\ref{80}), (\ref{427}) and (\ref{537}),
we obtain that
$$\begin{array}{lll}
\frac{\partial}{\partial \epsilon}\Big|_0y^\epsilon(t)(h)&=&\frac{\partial }{\partial \epsilon}\Big|_0h\Big(y^\epsilon(t)\Big)=\lim\limits_{\epsilon\to 0}\frac{h\Big(y^\epsilon(t)\Big)-h(\bar y(t))}{\epsilon}
\\  &=&\lim\limits_{\epsilon\to 0}\frac{h\Big(\eta(\rho(\bar y(t),y^\epsilon(t));\epsilon)\Big)-h\Big(\eta(0;\epsilon)\Big)}{\epsilon}
\\ &=&  \lim\limits_{\epsilon\to 0}\frac{1}{\epsilon}\Big(\langle \nabla h(\bar y(t)),\frac{\partial}{\partial \theta}\Big|_0\eta(0;\epsilon)\rangle\rho(\bar y(t),y^\epsilon(t))+o(\rho(\bar y(t),y^\epsilon(t)))\Big)
\\ &=& \lim\limits_{\epsilon\to 0}\frac{1}{\epsilon}\langle \nabla h(\bar y(t)),V_\epsilon(t)\rangle=V(t)(h),
\end{array}$$
which implies the first equality of (\ref{534}).

To prove the second equality of (\ref{534}), we employ two ways to compute $\frac{\partial^2}{\partial \epsilon^2}\Big|_0h(y^\epsilon(t))$.
On one hand,
\begin{equation}\label{529}\begin{array}{lll}
\frac{\partial^2}{\partial \epsilon^2}\Big|_0h(y^\epsilon(t))&=&\lim\limits_{\epsilon\to0}\frac{\partial}{\partial\epsilon}\Big|_0\langle\nabla h(y^\epsilon(t)),\frac{\partial}{\partial\epsilon}y^\epsilon(t)\rangle
\\  &=&  \langle\nabla_{V(t)}\nabla h,V(t)\rangle+\langle\nabla h(\bar y(t)),\frac{\partial^2}{\partial\epsilon^2}\Big|_0y^\epsilon(t)\rangle,
\end{array}\end{equation}
where the first equality of (\ref{534}) is used.

On the other hand, recalling (\ref{535}) and (\ref{536}), we can get, via Taylor's expansion, (\ref{80}) and (\ref{427}),
$$\begin{array}{lll}
&&\frac{\partial^2}{\partial \epsilon^2}\Big|_0h(y^\epsilon(t))
\\&=&\lim\limits_{\epsilon\to 0}\frac{1}{\epsilon^2}\Big(h(y^{2\epsilon}(t))-2h( y^\epsilon(t))+h(\bar y(t))\Big)
\\  &=& \lim\limits_{\epsilon\to0}\frac{1}{\epsilon^2}\Big\{h\Big(\eta(\rho(\bar y(t),y^{2\epsilon}(t));2\epsilon)
\Big)-h\Big(\eta(0;2\epsilon)\Big)-2\Big[h\Big(\eta(\rho(\bar y(t),y^\epsilon(t));\epsilon)\Big)
\\ &&-h\Big(\eta(0;\epsilon)\Big)\Big]\Big\}
\\ &=& \lim\limits_{\epsilon\to 0}\frac{1}{\epsilon^2}\Big\{\langle \nabla h(\bar y(t)),\frac{\partial}{\partial\theta}\Big|_0\eta(\theta;2\epsilon)\rangle\rho(\bar y(t),y^{2\epsilon}(t))+\frac{1}{2}\frac{\partial^2}{\partial \theta^2}\Big|_0h(\eta(\theta;2\epsilon))\rho^2(\bar y(t),y^{2\epsilon}(t))
\\&& -2\langle\nabla h(\bar y(t)),\frac{\partial}{\partial\theta}\Big|_0\eta(\theta;\epsilon)\rangle\rho(\bar y(t),y^\epsilon(t))-\frac{\partial^2}{\partial\theta^2}\Big|_0h(\eta(\theta;\epsilon))\rho^2(\bar y(t),y^\epsilon(t))+o(\epsilon^2)\Big\}.
\end{array}$$
Since
$$
\frac{\partial^2}{\partial\theta^2}\Big|_0h(\eta(\theta;\epsilon))=\frac{\partial}{\partial\theta}\Big|_0\langle\nabla h(\eta(\theta;\epsilon)),\frac{\partial}{\partial\theta}\eta(\theta;\epsilon)\rangle=\frac{1}{\rho^2(\bar y(t),y^\epsilon(t))}\langle\nabla_{V_\epsilon(t)}\nabla h,V_\epsilon(t)\rangle,
$$
where we have used Lemma \ref{317} and the fact that $\beta(\cdot;\epsilon)$ is a geodesic.
Thus, by using the above equality and (\ref{427}), we can get
$$\begin{array}{lll}
&&\frac{\partial^2}{\partial\epsilon^2}\Big|_0h(y^\epsilon(t))
\\ &=& \lim\limits_{\epsilon\to0}\frac{1}{\epsilon^2}\Big(\langle\nabla h(\bar y(t)),V_{2\epsilon}(t)\rangle+\frac{1}{2}\langle\nabla_{V_{2\epsilon}(t)}\nabla h,V_{2\epsilon}(t)\rangle-2\langle\nabla h(\bar y(t)),V_\epsilon(t)\rangle
\\ &&-\langle\nabla_{V_\epsilon(t)}\nabla h,V_\epsilon(t)\rangle
+o(\epsilon^2)\Big)
\\  &=&\lim\limits_{\epsilon\to 0}\frac{1}{\epsilon^2}\Big(\langle \nabla h(\bar y(t)),2\epsilon V(t)+4\epsilon^2 Y(t)\rangle+\frac{1}{2}\langle\nabla_{2\epsilon V(t)+4\epsilon^2 Y(t)}\nabla h,2\epsilon V(t)+4\epsilon^2 Y(t)\rangle
\\ &&-2\langle\nabla h(\bar y(t)),\epsilon V(t)+\epsilon^2Y(t)\rangle-\langle\nabla_{\epsilon V(t)+\epsilon^2 Y(t)}\nabla h,\epsilon V(t)+\epsilon^2 Y(t)\rangle+o(\epsilon^2)\Big)
\\ &=& 2\langle\nabla h(\bar y(t)), Y(t)\rangle+\langle\nabla_{V(t)}\nabla h,V(t)\rangle.
\end{array}$$
Combining this equality with (\ref{529}), one can get the second equality of (\ref{534}).

Secondly,  we claim that $(\nu,\psi_1)^\top\in (-\infty,0]\times T^*_{y_1}M$ belongs to $(\Im D_{\bar u(\cdot)}E)^\bot$ if and only if it satisfies (\ref{547}), where $\bar \psi(\cdot)$ is the solution to (\ref{64}) corresponding to the pair $(\nu, \psi_1)$. Thus, by the conditions of Theorem \ref{504}, the codimension of $\Im D_{\bar u(\cdot)}E$ is $1$.

In fact, assume that $(\nu,\psi_1)\in(-\infty, 0]\times T^*_{y_1}M$ satisfies (\ref{547}).
Recalling (\ref{65}), one can further get
\begin{equation}\label{530}
\nu\nabla_uf^0[t](v(t))+\nabla_uf[t](\bar\psi(t),v(t))=0,\,\,a.e.\,t\in[0,T],\quad \forall\; v(\cdot)\in L^2(0,T;\R^m).
\end{equation}
Integrating (\ref{530}) over $[0,T]$, we can obtain, via (\ref{527}), (\ref{64}), (\ref{534}) and integration by parts,
\begin{equation}\label{538}\begin{array}{lll}
0&=&\int_0^T\Big\{\nu\nabla_u f^0[t](v(t))+\nabla_{\dot{\bar y}(t)}V(\bar\psi(t))-\nabla_x f[t](\bar\psi(t),V(t))\Big\}dt
\\ & =&\psi_1(V(T))+\int_0^T\Big\{\nu\nabla_u f^0[t](v(t))-\nabla_{\dot{\bar y}(t)}\bar\psi(V)-\nabla_x f[t](\bar\psi(t),V(t))\Big\}dt
\\ &=&\nu\int_0^T\Big(\nabla_u f^0[t](v(t))+d_xf^0[t](V(t))\Big)dt+\psi_1(V(T))
\\ &=& \nu \frac{d}{d\epsilon}\Big|_{0}J(\bar u(\cdot)+\epsilon v(\cdot))+\psi_1(\frac{\partial}{\partial\epsilon}\Big|_0y^\epsilon(T))
\\ &=& (\nu, \psi_1)^\top \Big(D_{\bar u(\cdot)}E(v(\cdot))\Big).
\end{array}\end{equation}
Since $v(\cdot)\in L^2(0,T;\R^m)$ is arbitrarily chosen, we have
\begin{equation}\label{545}(\nu,\psi_1)^\top\in (\Im D_{\bar u(\cdot)}E)^\bot.\end{equation}

Conversely, if $(\tilde\nu,\tilde\psi_1)^\top\in (\Im D_{\bar u(\cdot)}E)^\bot$ with $\tilde\nu\leq 0$, one can follow the above argument in an inverse way, and get (\ref{547}) with $\nu=\tilde\nu$ and $\bar \psi(\cdot)=\tilde\psi(\cdot)$, where $\tilde\psi(\cdot)$ is the solution to (\ref{64}) with $\nu=\tilde\nu$ and $\psi_1=\tilde\psi_1$.

Thirdly,
according to (\ref{532}),   the above argument  and Lemma \ref{533},  $(\nu,\psi_1)^\top Hess_{\bar u(\cdot)}E$ is sign-definite on $Ker D_{\bar u(\cdot)}E\times Ker D_{\bar u(\cdot)}E$. Now, we need to compute explicitly both $(\nu,\psi_1)^\top Hess_{\bar u(\cdot)}E$ and $Ker D_{\bar u(\cdot)}E$.

Recalling (\ref{539}), (\ref{540})  and (\ref{534}), we have
\begin{equation}\label{541}\begin{array}{lll}
(\nu,\psi_1)^\top Hess_{\bar u(\cdot)}E(v(\cdot),v(\cdot))
&=&(\nu,\psi_1)^\top\Big( \frac{\partial^2}{\partial \epsilon^2}\Big|_0 E(\bar u(\cdot)+\epsilon v(\cdot))\Big)
\\ &=& \nu\frac{\partial^2}{\partial \epsilon^2}\Big|_0 J(\bar u(\cdot)+\epsilon v(\cdot))+2\psi_1\Big(Y(T)\Big).
\end{array}\end{equation}
Applying (\ref{534}), we can get
\begin{equation}\label{542}\begin{array}{lll}
&&\nu\frac{\partial^2}{\partial \epsilon^2}\Big|_0 J(\bar u(\cdot))
\\ &=&\nu\frac{\partial}{\partial\epsilon}\Big|_0\int_0^T\Big\{\langle\nabla_x f^0(t,y^\epsilon(t),\bar u(t)+\epsilon v(t)),\frac{\partial}{\partial\epsilon}y^\epsilon(t)\rangle+\nabla_u f^0(t,y^\epsilon(t),\bar u(t)+\epsilon v(t))(v(t))\Big\}dt
\\ &=& \nu\int_0^T \Big\{\nabla_x^2f^0[t](V(t),V(t))+2\nabla_x\nabla_uf^0[t](v(t),V(t))+\nabla_u^2f^0[t](v(t),v(t))
\\ &&+\langle\nabla_xf^0[t],2Y(t)\rangle
\Big\}dt.
\end{array}\end{equation}
Using integration by parts and employing the second order variational equation (\ref{528}), we obtain that
$$\begin{array}{lll}
&&2\psi_1(Y(T))
\\ &=& 2\int_0^T\Big[\nabla_{\dot{\bar y}(t)}\bar \psi(Y(t))+\bar \psi(t)(\nabla_{\dot{\bar y}(t)}Y)\Big] dt
\\ &=& 2\int_0^T\Big\{\nabla_{\dot{\bar y}(t)}\bar \psi(Y(t))+\nabla_xf[t](\bar \psi(t),Y(t))+\nabla_u\nabla_xf[t](\bar \psi(t),V(t),v(t))
\\ &&-\frac{1}{2}R(\tilde{\bar\psi}(t),V(t),f[t],V(t))+\frac{1}{2}\nabla_x^2f[t](\bar\psi(t),V(t),V(t))+\frac{1}{2}\nabla_u^2f[t](\bar
\psi(t),v(t),v(t))\Big\}dt.
\end{array}$$

Inserting the above identity and (\ref{542}) into (\ref{541}), one can get, via the first order dual equation (\ref{64}),
$$\begin{array}{lll}
&&(\nu, \psi_1)^\top Hess_{\bar u(\cdot)}E(v(\cdot),v(\cdot))
\\&=&\int_0^T\Big\{\nabla_x^2 H^\nu(t,\bar y(t),\bar \psi(t),\bar u(t))(V(t),V(t))+\nabla_u^2 H^\nu(t,\bar y(t),\bar \psi(t),\bar u(t))(v(t),v(t))
\\  &&+2\nabla_u\nabla_x H^\nu(t,\bar y(t),\bar \psi(t),\bar u(t))(V(t),v(t))
-R(\tilde{\bar\psi}(t),V(t),f[t],V(t))\Big\}dt.
\end{array}$$

By the definition of the kernel of $D_{\bar u(\cdot)}E$, for any  $v(\cdot)\in Ker D_{\bar u(\cdot)}E$, one has
\begin{equation}\label{543}\begin{array}{lll}
0=\frac{d}{d\epsilon}\Big|_0\left( \int_0^Tf^0(t,y^\epsilon(t),\bar u(t)+\epsilon v(t))dt\atop y^\epsilon(T)\right)
=\left(\int_0^T\left\{\langle\nabla_xf^0[t],V(t)\rangle+\nabla_uf^0[t](v(t))\right\}dt\atop V(T)\right).
\end{array}\end{equation}

 In the normal case, i.e. $\nu<0$, (\ref{543}) implies, together with the first order dual equation (\ref{64}) and (\ref{527}),
$V(T)=0$ and $\int_0^T\nabla_uH^\nu[t](v(t))dt=0$. By (\ref{547}), it follows that
 \begin{equation}\label{546}\begin{array}{lll}Ker D_{\bar u(\cdot)}E&=&\{v(\cdot)\in L^2(0,T;\R^m);  \,\textrm{  (\ref{522}) with $\xi(\cdot)=v(\cdot)$  admits a solution}\}.\end{array}\end{equation}

In order to determine the sign of the quadratic $(\nu,\psi_1)^\top Hess_{\bar u(\cdot)}E$, from the assumptions,
one can find a $u_\eta(\cdot)\in \V_{ad}\setminus\{\bar u(\cdot)\}$ such that $\|u_\eta(\cdot)-\bar u(\cdot)\|_{L^2(0,T; \R^m)}$ tends to zero as $\eta\to +\infty$. Set
$$
v_\eta(\cdot)=\frac{u_\eta(\cdot)-\bar u(\cdot)}{\|u_\eta(\cdot)-\bar u(\cdot)\|_{L^2(0,T; \R^m)}},\quad\epsilon_\eta=\|u_\eta(\cdot)-\bar u(\cdot)\|_{L^2(0,T; \R^m)}.
$$
Since $D_{\bar u(\cdot)}E: L^2(0,T;\R^m)\to \R\times T_{y_1}M$ is linear, one can find a subspace $N\subset L^2(0,T;\R^m)$ such that $L^2(0,T;\R^m)=Ker D_{\bar u(\cdot)}E\oplus N$, $\Im D_{\bar u(\cdot)}E=(D_{\bar u(\cdot)} E) N$, and $D_{\bar u(\cdot)}: N\to \R\times  T_{y_1}M $ is surjective. Thus, we have $\dim N<\infty$. For each $v_\eta(\cdot)$, there exist a $v_\eta^1(\cdot)\in Ker D_{\bar u(\cdot)}E$ and a $v_\eta^2(\cdot)\in N$ such that $v_\eta(\cdot)=v_\eta^1(\cdot)+v_\eta^2(\cdot)$. In what follows, we will consider two cases.

Case I.\quad If there exists an $n_1>0$ such that $v_\eta^2(\cdot)\in N\setminus\{0\}$ for all $n>n_1$, then one can find a subsequence of $\{v_\eta^2(\cdot)\}$ (still denoted by $\{v_\eta^2(\cdot)\}$) and a $v^2(\cdot)\in N$ such that \begin{equation}\label{330}\lim_{\eta\to+\infty}\|v_\eta^2(\cdot)-v^2(\cdot)\|_{L^2(0,T;\R^m)}=0.
\end{equation}  We claim that $v^2(\cdot)=0$.

Indeed, denote by $V_\eta^i(\cdot)$ and $V^2(\cdot)$ the solutions to the first order variational equation (\ref{527}) with $v(\cdot)= v_\eta^i(\cdot)$ for $i=1,2$ and $v(\cdot)=v^2(\cdot)$, respectively. Especially, applying (\ref{330}), we have
\begin{equation}\label{329}\lim_{\eta\to +\infty}\max_{t\in[0,T]}|V_\eta^2(t)-V^2(t)|=0.
\end{equation}
Denote by $V_\eta(\cdot)=V_\eta^1(\cdot)+V_\eta^2(\cdot)$. One can check that $V_\eta(\cdot)$ is the solution to (\ref{527}) with $v(\cdot)=v_\eta(\cdot)=v_\eta^1(\cdot)+v_\eta^2(\cdot)$. Denote by $y_\eta^\epsilon(\cdot)$ the solution to (\ref{25}) with cotrol $u(\cdot)=\bar u(\cdot)+\epsilon v_\eta(\cdot)$ for $\epsilon\geq 0$. Set $V_{\eta,\epsilon}(t)=\exp_{\bar y(t)}^{-1}(y_\eta^\epsilon(t))$ for $t\in [0,T]$. Applying (\ref{427}), one can obtain
$$
V_{\eta,\epsilon}(t)=\epsilon V_\eta(t)+\epsilon^2 Y_\eta(t)+o(\epsilon^2),\quad\forall\; t\in[0,T],
$$
where $Y_\eta(\cdot )$ is the solution to (\ref{528}) with $V(\cdot)=V_\eta(\cdot)$ and $v(\cdot)=v_\eta(\cdot)$.
Recalling the proof of Proposition \ref{439}, one can check that $\epsilon^2 Y_\eta(t)+o(\epsilon^2)$ in the above formula is an infinitesimal $o(\epsilon)$, which is uniform with respect to $v_\eta(\cdot)$, due to $\|v_\eta(\cdot)\|_{L^2(0,T;\R^m)}=1$.
Note that $u_\eta(\cdot)=\bar u(\cdot)+\epsilon_\eta v_\eta(\cdot)\in\V_{ad}$, $\;\;v_\eta^1(\cdot)\in Ker D_{\bar u(\cdot)}E$ and (\ref{546}). The above formula holds for $\epsilon=\epsilon_\eta$ and
$$
0=V_{\eta,\epsilon_\eta}(T)=\epsilon_\eta V_\eta^2(T)+o(\epsilon_\eta),\quad \eta=n_1+1, n_1+2, \cdots.
$$
From the above formula and (\ref{329}), one can get $V^2(T)=0$, and hence $v^2(\cdot)\in Ker D_{\bar u(\cdot)}E\cap N$, which implies the claim.

Applying  Taylor's expansion, and noting (\ref{534}) and (\ref{545}) and the boundness of $\{v_\eta(\cdot)\}$ in $L^2(0,T; \R^m)$, one can get
\begin{equation}\label{328}\begin{array}{ll}
&\nu\Big(J(\bar u(\cdot)+\epsilon v_\eta(\cdot))-J(\bar u(\cdot))\Big)+\psi_1(V_{\eta,\epsilon}(T))
\\=&\epsilon\Big(\nu \Big(D_{\bar u(\cdot)}J\Big)v_\eta(\cdot)+\psi_1(V_\eta(T))\Big)+\frac{1}{2}\epsilon^2\Big(\nu Hess_{\bar u(\cdot)}J(v_\eta(\cdot),v_\eta(\cdot))+\psi_1(\frac{\partial^2}{\partial\epsilon^2}\Big|_0y^\epsilon_\eta(T))\Big)
\\ &+o(\epsilon^2)
\\ =& \epsilon (\nu,\psi_1)^\top\Big(D_{\bar u(\cdot)}E v_\eta(\cdot)\Big)+\frac{1}{2}\epsilon^2 (\nu,\psi_1)^\top Hess_{\bar u(\cdot)}E(v_\eta(\cdot),v_\eta(\cdot))+o(\epsilon^2)
\\ =& \frac{1}{2}\epsilon^2(\nu,\psi_1)^\top Hess_{\bar u(\cdot)}E(v_\eta^1(\cdot), v_\eta^1(\cdot))+\epsilon^2(\nu,\psi_1)^\top Hess_{\bar u(\cdot)} E(v_\eta^1(\cdot), v_\eta^2(\cdot))
\\&+\frac{1}{2}\epsilon^2(\nu,\psi_1)^\top Hess_{\bar u(\cdot)}E(v_\eta^2(\cdot),v_\eta^2(\cdot))+o(\epsilon^2).
\end{array}\end{equation}
Taking $\epsilon=\epsilon_\eta$ in the above formula, we have, via $\nu<0$, $y_\eta^{\epsilon_\eta}(T)=y_1$, (\ref{545}) and the optimality of $\bar u(\cdot)$,
\begin{equation}\label{331}\begin{array}{ll}
0\geq&\frac{1}{2}(\nu,\psi_1)^\top Hess_{\bar u(\cdot)}E(v_\eta^1(\cdot),v_\eta^1(\cdot))+(\nu,\psi_1)^\top Hess_{\bar u(\cdot)}E(v_\eta^1(\cdot),v_\eta^2(\cdot))
\\ &+\frac{1}{2}(\nu,\psi_1)^\top Hess_{\bar u(\cdot)}E(v_\eta^2(\cdot),v_\eta^2(\cdot))+o(1).
\end{array}\end{equation}
By means of the boundedness of  $\{v_\eta^2(\cdot)\}$ in $L^2(0,T;\R^m)$, there exists a subsequence of
$\{(\nu,\psi_1)^\top Hess_{\bar u(\cdot)}E(v_\eta^1(\cdot),v_\eta^1(\cdot))\}$ (still denoted by $\{(\nu,\psi_1)^\top Hess_{\bar u(\cdot)}E(v_\eta^1(\cdot),v_\eta^1(\cdot))\}$ ) such that $\lim_{\eta\to +\infty}(\nu,\psi_1)^\top Hess_{\bar u(\cdot)}E(v_\eta^1(\cdot),v_\eta^1(\cdot))=\xi$ exists. Taking limit in (\ref{331}) as $\eta\to +\infty$, one can get, via (\ref{330}), $0\geq\frac{1}{2}\xi$. Since $(\nu,\psi_1)^\top Hess_{\bar u(\cdot)}E$ is sign definite on $Ker D_{\bar u(\cdot)}E$, we get $(\nu,\psi_1)^\top Hess_{\bar u(\cdot)}E$ is negative definite on $Ker D_{\bar u(\cdot)}E$.

Case II.\quad There exists a subsequence of $\{v_\eta(\cdot)\}$ (still denoted by $\{v_\eta(\cdot)\}$ ) such that $v_\eta(\cdot)\in Ker D_{\bar u(\cdot)}E$. Applying the same argument as that in (\ref{328}) and (\ref{331}), one can get $(\nu,\psi_1)^\top Hess_{\bar u(\cdot)}E$ is negative definite on $Ker D_{\bar u(\cdot)}E$.

In the abnormal case, i.e. $\nu=0$, $Ker D_{\bar u(\cdot)}E$ is equal to (\ref{544}). Lemma \ref{533} concludes that the left hand of (\ref{507}) is sign-definite on (\ref{544}).   $\Box$

\medskip

Finally, let us prove Theorem \ref{453}.

\textbf{Proof of Theorem \ref{453}}\quad
Theorem \ref{453} is a direct consequence of  \cite[Theorem 21.8, p. 347]{as}. $\Box$

{}

\end{document}